\documentclass[twoside]{article}
 \usepackage{amsthm,amsmath,amssymb}
{\newtheorem{theorem}{Theorem}[section]                     
        
   \newtheorem{lemma}[theorem]{Lemma}
   \newtheorem*{claim}{Claim}
   \newtheorem{corollary}[theorem]{Corollary}
   
   \newtheorem{question}[theorem]{Question}
}
{\theoremstyle{definition}
   \newtheorem{exercise}[theorem]{Exercise}
   \newtheorem{example}[theorem]{Example}
   \newtheorem{definition}[theorem]{Definition}
   \newtheorem{remark}[theorem]{Remark}
}
\newcommand{\n}{\noindent}
\newcommand{\vn}{\vspace{5mm} \noindent}

\newcommand{\mylabel}[1]{\label{#1}}

\newcommand{\Spec}{\operatorname{Spec}}
\newcommand{\Isom}{\operatorname{Isom}}
\newcommand{\Sing}{\operatorname{Sing}}
\newcommand{\Star}{\operatorname{Star}}
\newcommand{\Jac}{{\operatorname{Jac }}}
\newcommand{\Pic}{{\operatorname{\mathbf Pic }}}
\newcommand{\Hom}{{\operatorname{Hom}}}
\newcommand{\Div}{{\operatorname{Div}}}
\newcommand{\Ext}{{\operatorname{Ext}}}
\newcommand{\Mor}{{\operatorname{Mor}}}
\newcommand{\End}{{\operatorname{End}}}
\newcommand{\Mat}{{\operatorname{Mat}}}
\newcommand{\Aut}{{\operatorname{Aut}}}
\newcommand{\Sym}{{\operatorname{Sym}}}
\newcommand{\Def}{{\operatorname{Def}}}

\newcommand{\Sm}{{\operatorname{Sm}}}
\newcommand{\Tor}{\operatorname{Tor}}
\newcommand{\Hl}{\operatorname{H}}
\newcommand{\length}{\operatorname{length}}
\newcommand{\codim}{\operatorname{codim}}
\newcommand{\chara}{\operatorname{char}}
\newcommand{\das}{\dashrightarrow}
\newcommand{\dar}{\downarrow}
\newcommand{\bfA}{{\mathbb A}}
\newcommand{\Syp}{{\operatorname{Sp}}}
\newcommand{\PGL}{{\operatorname{PGL}}}
\newcommand{\hardexercise}{{\LARGE${\mathbf *}$}}
\newcommand{\RR}{{\mathbb{R}}}
\newcommand{\CC}{{\mathbb{C}}}
\newcommand{\QQ}{{\mathbb{Q}}}

\newcommand{\PP}{{\mathbb{P}}}
\newcommand{\ZZ}{{\mathbb{Z}}}

\newcommand{\GG}{{\mathbb{G}}}
\newcommand{\AAA}{{\mathbb{A}}}
\newcommand{\cC}{{\cal C}}
\newcommand{\cD}{{\cal D}}
\newcommand{\cF}{{\cal F}}
\newcommand{\cG}{{\cal G}}
\newcommand{\cH}{{\cal H}}
\newcommand{\cI}{{\cal I}}
\newcommand{\cL}{{\cal L}}
\newcommand{\cO}{{\cal O}}
\newcommand{\cP}{{\cal P}}
\newcommand{\cQ}{{\cal Q}}
\newcommand{\cX}{{\cal X}}
\newcommand{\cY}{{\cal Y}}
\newcommand{\cZ}{{\cal Z}}
\newcommand{\M}{M_g}
\newcommand{\oM}{\overline{M_g}}
\newcommand{\mM}{M^{(m)}_g}
\newcommand{\omM}{\overline{M^{(m)}_g}}
\newcommand{\nM}{M_{g,n}}
\newcommand{\onM}{\overline{M_{g,n}}}
\newcommand{\mnM}{{M^{(m)}_{g,n}}}
\newcommand{\omnM}{\overline{M^{(m)}_{g,n}}}
\newcommand{\cM}{{\cal M}_g}

\newcommand{\conM}{\overline{{\cal M}_{g,n}}}

\newcommand{\cnM}{{{\cal M}_{g,n}}}
\newcommand{\cmnM}{{{\cal M}^{(m)}_{g,n}}}
\newcommand{\setmin}{%
\begin{picture}(6,10)\qbezier(1,5.)(3,4.)(5,3.)\end{picture}}

\title{{\bf Alterations and resolution of singularities}}
\author{Dan Abramovich and Frans Oort}
\date{}
\begin{document}
\maketitle 
 
\tableofcontents

\addtocounter{section}{-1}
\newpage
            \section{Introduction}
%
\noindent
H. Hironaka, 1964: 
\begin{center}
{\it In characteristic zero \\ any variety can be modified into a
nonsingular variety.}  
\end{center}
A. J. de Jong, 1995:
\begin{center} 
{\it Any variety can be altered into a nonsingular variety.}
\end{center}
On July 26, 1995, at the University of California, Santa Cruz, a young Dutch
mathematician by the name Aise Johan de Jong made a revolution in the study of
the 
arithmetic, geometry and cohomology theory of varieties in positive or mixed
characteristic. The talk he delivered, first in a series of three entitled 
``Dominating Varieties by Smooth Varieties'', had a central theme: a
systematic application of fibrations by nodal curves. Among the hundreds of awe
struck members of the audience, participants of the American Mathematical
Society Summer Research Institute on Algebraic Geometry, many recognized the
great potential of Johan de Jong's ideas even for {\em complex} algebraic
varieties, and indeed soon more results along these lines began to form.

\subsection{The alteration paradigm}
A.J. de Jong's main result was, that for any variety $X$, there is a
nonsingular variety $Y$ and  an
{\em alteration,} namely a proper, surjetive and generically finite
morphism, $Y\to 
X$  (see Theorem \ref{Th:deJongMain} for a precise statement). This is
in contrast with Hironaka's result, which uses only a {\em modification,}
namely a proper birational morphism.

Here is the basic structure of the proof by De Jong:
\begin{itemize}
\item {\bf Projection.} For a given variety $X$ of dimension $d$ we
produce a morphism $f: X 
\to P$ with $\dim  P = d-1$, and all fibers of $f$ are curves (we may first
have to apply a modification to $X$).
\item {\bf Desingularization of fibers.} After an alteration of the
base $P$, we arrive at a new morphism 
$f: X \to P$ where all fibers are curves with only ordinary
nodes as singularities.  The main tool here is the theory of moduli of
curves. 
\item {\bf Desingularization of base.}  After a further alteration on
the base $P$, we arrive at a new morphism 
$f: X \to P$ as above, where $P$ is regular.  Here we use induction,
i.e. supposing that the theorem is already true for varieties of dimension
$d-1$. So here we ``desingularize the base''.
\item  {\bf Desingularization of total space.} Given the last two
steps, an explicit and easy method of resolution of 
singularities finishes the job. 
\end{itemize} 

\subsection{The purpose of this paper}

This paper is an outgrowth of our course material prepared for the  Working
Week on Resolution of Singularities, which was held during September 7-14, 1997
in  Obergurgl, Tirol, Austria. As we did in the workshop, we intend to explain
Johan de Jong's results in some detail, and give some other results
following the same paradigm, as well as a few
applications,
 both arithmetic
and in characteristic zero. We hope that  the 
reader will come to share some of the excitement we felt on that
beautiful July day in Santa Cruz. 

\vn
In the rest of this introduction we give an overview 
of the proof and the material involved. We hope that this introduction will
give most readers a general feeling of what the results are about. The
body of the paper is divided in two parts. We begin part I by expanding on some
of the preliminary material necessary for understanding the proofs by
any student of 
algebraic geometry. Then we go back to the proof of de Jong's main theorem, as
well as some generalizations.  Proofs of some variants and
generalizations of de Jong's 
theorems are indicated in the form of exercises, with sufficient hints and
references, which we hope will  enable the reader to appreciate de Jong's
work. Part II is an introduction to an ingredient of
the proof - the theory of moduli of curves. We aim to indicate the
main ideas behind the proofs of the main theorems about existence and
properties of moduli spaces, again accompanied with a  collection of
exercises.

\vn
As a result, this account is mostly expository. The only point where some  
novelty  appears is in Section \ref{Sec:Tautological-pointed}, where we show
the existence of  tautological families  of stable curves over the moduli
spaces of {\em stable pointed curves with level structure.} This has
been ``well known 
to the experts'' for years, and can be collected from the
literature. However a complete account under one roof has not been published. For the definition of a ``tautological curve" we refer to \ref{taut}.

\subsection{Historical context}

There are many cases in geometry in which one wants  to transform a singular
variety into a non-singular one: once arrived in such a situation, various
technical steps can be performed, not possible on singular varieties.  

Since the beginning of the century,  partial results in this direction
appeared, crowned by  
Hironaka's theorem on resolution of singularities in 
characteristic zero, in 1964. 

 Hironaka's ingenuous proof had many applications, but it was not easy to
understand the fine details of his proof. Generalizing that method to varieties
in positive characteristic has failed up to now. Indeed, resolution of
singularities  in positive characteristic has been a topic to which many years
of intensive research have been devoted, and up to now the status is not yet
clear: for the general question of resolution of singularities in positive
characteristic we have neither a fully verified theorem nor a
counterexample. In addition, the algorithms involved  in Hironaka's theory were
difficult  to generalize, even in characteristic 0, to some important
more complicated situations.

It seemed that  a lull in development of this subject had been
reached, until a 
totally new idea came about. In 1995 Johan de Jong approached the problem
above, of transforming a variety into a nonsingular one, from
a different angle. The idea of the proof is surprisingly easy, and for
many applications his result is sufficient. His approach is very geometric,
and hence it works in a wide range of situations.  The alteration
paradigm automatically 
works in all characteristics, and a suitable version works in
mixed characteristic as well. It easily gives rise to some new
``semistable reduction'' type results which are new even over the
complex numbers. 
Moreover, without much effort it give birth to new, ``conceptually easy''
proofs of a weaker form of Hironaka's theorem.

\subsection{Comparison of approaches}

Let us take a moment to make a qualitative comparison of Hironaka's result and
de Jong's result. 

In the approach taken by Hironaka,  singularities   of a variety are studied
closely, invariants measuring the difficulty of the singularities are defined,
and a somewhat explicit algorithm is applied  
 in order to improve the singularities, in the sense that the
given invariants get ``better''. One needs to show that the algorithm 
 terminates (and indeed  in characteristic zero it does), resulting in the
construction of a 
regular variety. A big advantage of this process developed by Hironaka (and by
many others) is the fact that usually it is very explicit, it is canonical in a
certain sense and once it works, the result is in its strongest form,
see \cite{Villa} and \cite{Bie-Mil}, as well as \cite{Enc-Vil} in this
volume. 

In the approach by Johan de Jong, the singularities are, at first, completely
ignored. The idea is to first bring the variety to a special form: a fibration
by nodal curves. Here one pays a big price: in order to arrive at this special
form one needs to use an operation - called alteration - which extends the
function field of the variety. However, once we arrive at this form, we can use
induction on the dimension for the base space of the fibration, and
automatically arrive at a situation where the variety has very mild
singularities. Only then, finally, attention is paid to the singularities. But
these are so mild that an easy and explicit blowing up finishes the job.

\subsection{A sketch of the construction of an alteration giving a regular
variety}  Here we give a much simplified form of the proof of A.J. de Jong's
main Theorem  (Theorem \ref{Th:deJongMain} in this text). We break up
the proof in steps. A star attached to a step means 
that in that phase of the proof a finite extension of the function field might
be involved, i.e. the alteration constructed might not be a modification. In
steps without a star only modifications are used.

Before starting, a small technical point is necessary. In the course of the
proof we use induction on the dimension of the variety 
$X$, and it turns out 
that for the induction to work we need the statement of the theorem to involve
a closed subset $Z\subset X$ as well. Our final goal will be to find an
alteration 
$f:Y \to X$ such that $f^{-1}Z$ is a normal crossings divisor.

{\it We start with a field $k$, a variety $X$ and a closed subset $Z\subset X$,
over the field  $k$.}

\begin{description} \large
\item{\sc Step 0.} 

{\it We can reduce to the case where $k$ is algebraically
closed, the variety $X$ is projective and normal, and the closed
subset $Z$ a divisor.}   
\end{description}

\n
We intend to say: if we prove the theorem with this new additional data, then
the theorem in the original, more general form follows. Reducing to the
algebraically closed field case is standard - in the main body of the
paper we avoid it, assuming $k$ is algebraically closed. The main
ingredient for projectivity is
Chow's Lemma (see \cite{RedBook}, pp. 85-89, or \cite{HAG}, Exercise
II.4.10): 
{\em for a variety $X$ over $K$, there exists a modification $X' \rightarrow
X$, such that $X'$ is quasi-projective.} To make $Z$ into a divisor we simply
blow it up inside $X$.

\noindent{\em Replacement convention.} From now on, in each step,  we shall
replace $X$ by a new 
variety $X'$ over $k$ which admits a modification or an alteration $X'
\rightarrow X$, arriving finally at a regular variety and an alteration of the
variety produced in Step 0.

\begin{description}\large
\item{\sc  Step 1.} 

{\em After modifying  $X$, construct a morphism
$f: X 
\rightarrow P$ of projective varieties 
whose generic fiber
is an irreducible, complete, non-singular curve.}
\end{description}

\vn
Note: $\dim(P) = \dim(X) -1$, which suggests using induction later. 

Actually we need a little more, but the technical details will be discussed in
the main text.

This step follows a classical, geometric idea. 
Set  $\dim(X)=d$, and  assume $X \subset \PP^N$. Using Bertini's theorem we see
that we can find a 
linear subvariety $L \subset \PP^N$ ``in general position" with $\dim(L) = N-d$
such that the projection with center $L$ gives a rational map $X
\dashrightarrow \PP^{d-1}$ where the generic fiber is a regular curve. After
modifying $X$ we can make this rational map into a morphism.

\vn
{\em The strict transform.} We will use
an operation which  de Jong called the ``strict transform''. (In 
\cite{Berthelot}, 815-12 the 
terminology ``strict alteration'' is used). Consider a morphism $X \to S$, and 
a base change $T \to S$. Assume $T$ to be integral, and let $\eta \in T$ be its
generic point. Then define $X' \subset T \times _S X$ as the closure of the
generic fiber $(T\times_S X)_{\eta}$ in  $(T\times_S X)$. A more thorough
discussion of this operation will follow in Section
\ref{Sec:Strict-Transform}.  

In our situation $X \to P$, we will often replace $P$ by an alteration, and
then simply replace $X$ by its strict transform.

\begin{description}\large
\item{\sc  Step 2$^*$.} 

{\em After applying alterations to $X$
and to 
$P$ we can arrive at a morphism $f: X \rightarrow P$ as in Step 1, and 
disjoint sections $\sigma_1, \cdots , \sigma_n:P\to X$, such that every
geometric component $C'$ of every geometric fiber of $f$ meets at least three
of these sections in  the smooth locus of $f$, i.e. in $C'\cap \Sm(f)$.}
\end{description}

\n
There is a ``multi-section" in the situation of Step 1 having this
property. After an alteration on $Y$ and on $X$ this becomes a union of
sections. 

\noindent{\em Stable pointed curves.} 
Here we follow Deligne \& Mumford, and
Knudsen. 
An algebraic curve is called {\it nodal} if it is complete, connected
and  
if the singularities of $C$ are not worse than ordinary double points. Its
arithmetic genus is given by $g = \dim_k \ \Hl^1(C, {\cal O}_C).$ 

Suppose $C$ is a nodal curve of genus $g$ over a field $k$, and let $P_1,
\cdots , P_n \in C(k)$  with $2g-2+n > 0$; we write ${\cal P} = \{P_1, \cdots ,
P_n\}$; this is called a {\it stable $n$-pointed curve} if: 
\begin{itemize}
\item the points are mutually different, $i<j \Longrightarrow P_i \not= P_j$, 
\item none of these marked points is singular, $P_i \not\in \Sing(C)$, \item
and $\Aut(C,{\cal P})$ is a finite group; under the previous conditions (and
$k$ algebraically closed) this amounts to the condition that for every regular
rational irreducible component 
$$ \PP_1 \cong C' \subset C, \ \  \mbox{then} \ \ \ \#(C' \cap ({\cal P} \cup
\Sing(C)) )\geq 3.$$ 
\end{itemize} 
A flat family of curves is called ``a family of  stable $n$-pointed curves'' if
all geometric fibers 
are stable $n$-pointed curves in the sense just defined, the markings given by
sections.

\vn
Historically, stable curves and stable pointed curves were introduced in order
to construct, in a natural way, compactifications of moduli spaces
(see \cite{D.M}). Certainly the  
following names should be mentioned: Zariski, A. Mayer, Deligne, Mumford,
Grothendieck, Knudsen, and many more.  It came a bit as a surprise when de
Jong used these for a desingularization-type problem!

\begin{description}\large
\item{\sc Step 3$^{\ast}$.} 

{\em After an alteration on the base $P$, we  can assume 
that $X \to P$ is a projective family  of stable
$n$-pointed curves.}
\end{description}

\vn
We briefly sketch the heart of the proof of this step - it will be discussed in
detail later.

\noindent{\em Extending families of curves.} We need the following
fundamental 
fact: suppose we are given a variety $P$, an open dense subset $U \subset
P$, and a family of stable  curves  
$C_U \to U$:
$$\begin{array}{ccc} C_U &\subset & ? \\
			\dar & & \dar \\
			U &\subset& P.
\end{array} $$
Then there is an alteration $a: P_1 \to P$  such that the pullback
family $C_{U_1} \to U_1$ over the open set $U_1 = a^{-1}U$  can be 
extended to a family of stable curves $C_1 \to P_1$:
$$\begin{array}{ccc} C_{U_1} &\subset & C_1 \\
			\dar & & \dar \\
			U_1 &\subset& P_1.
\end{array} $$

\n
The first result behind this is the existence of a  moduli space of stable
curves (\cite{Knudsen}, see also Section
\ref{Sec:Moduli-construction}). Then one needs the fact that a
finite cover $M \to \onM$ of the  
moduli space admits a ``tautological  family'' - namely, a family $C
\to M$ such that the associated morphism $M \to \onM$ is the given
finite cover.  
 One could consult \cite{Deligne} (the precise statement we need 
follows from that paper), or use \cite{vG.FO}, where a tautological
family of nodal curves is constructed over a moduli space of stable
 curves with a level structure.

 The markings of the family $X' \to P$
correspond to those of  the stable $n$-pointed curve ${\cal C}
\rightarrow P$, 
under the birational transformation thus defined.  We want to show this
extends to a morphism ${\cal C} \rightarrow X'$. 

\vn
{\em Flattening of the graph.}  We  take the closure $T \subset X
\times _P {\cal C}$ of the graph of 
$\beta_0:{\cal C}_U \to X_U$, and  apply the "Flattening Lemma", see
\ref{Lem:flattening} 
below. We arrive at new $X, T$, and ${\cal C}$ flat over $P$. All we have to
show 
(modulo some technicalities) is that no point of a fiber of ${\cal C} \to P$
is blown up to a component of  a fiber of $X \to P$. 

\vn{\em The Three Point Lemma.} Using the markings, and
studying carefully the geometry we show that indeed $\beta_0$ extends to a
morphism $\beta$. The crucial point here was that every component of every
fiber of $X$ over $P$ has at least three nonsingular points marked by the
sections $\sigma_i$  (see 4.18 - 4.20 of \cite{dj}).  

\begin{description}\large
\item{\sc Step 4$^*$.} 

{\em After an alteration of $P$, we may assume that $P$ is nonsingular.}
\end{description}

\n
We simply  apply induction on the dimension of the base: we
suppose that the 
theorem we want to prove is valid for all varieties having dimension less than
$\dim \ X$. Thus after an alteration of the base $P$ we can suppose $P$ is
regular and the strict transform of $X$ has all the previous properties.

\noindent{\em Following $Z$.} The argument for the previous two steps should
be carried through with a 
proper care given to the divisor $Z$. At the end, we can guarantee that $Z$
is contained in the union of two types of sets:
\begin{itemize}
\item the images of the sections $\sigma_i$, and
\item the inverse image of a normal crossings divisor  $\Delta\subset P$.
\end{itemize}
Moreover, in the induction hypothesis we can guarantee that the final family of
curves $X \to P$ degenerates only over the normal crossings divisor $\delta$.

\begin{description}\large
\item{\sc Step 5.} 

{\em The singularities of the resulting family $X \to P$
are so mild that it is very easy to resolve them explicitly}.
\end{description}

Indeed, each singular point can be described in formal coordinates be the
equation $xy = t_1^{k_1}\cdots t_r^{k_r}$. It is a fairly straightforward
exercise to resolve these singularities.

\newpage\part{The alteration theorem}
            \section{Some preliminaries and generalities on varieties}
%
\subsection{Varieties} To fix notation, we use the following definition of a
variety: 
\begin{definition} By a {\em variety} defined over $k$ we
mean a 
geometrically integral scheme of finite type over $k$. If $k\subset k_1$ we
write $X_{k_1}$ for $X\times_{\Spec k}\Spec k_1$.
\end{definition}

\n
In more down to earth terms this means: an {\it affine} variety defined over
$k$ is 
given as 
a closed subvariety of an affine space $\bfA_k^n$ defined by an ideal $I
\subset k[T_1, \cdots , T_n] = k[T]$ such that  $k_1{\cdot}I \subset k_1[T]$ is
a prime ideal for every (equivalently, for some) algebraically closed field
$k_1$ 
containing $k$. In general, a variety 
then is defined by gluing a finite number of affine varieties in a separated
way.  See \cite{RedBook}, I.5, Definition 1 (p. 35) and I.6,
Definition 2 (p. 52). 
 
\begin{remark} This definition differs slightly from that in \cite{dj}. De
Jong requires the algebraic scheme to be integral, and we require that the
schemes stay integral after extending the field. 
For example for any finite field extension $k \subset K$, the scheme
$\Spec(K)$ is called a $k$-variety by de Jong, but we only say it is a variety
defined over $k$ if $k=K$. For most geometric situations the differences will
not be important. 
\end{remark}

\subsection{Operations on varieties}
\begin{definition} A morphism of varieties  $Y \to X$ is called
a {\em modification} if it is proper and birational.
\end{definition}

\n
A modification is the type of ``surgery operation'' usually associated with
resolution of singularities. Johan de Jong introduced the following important
variant: 

\begin{definition}[de Jong] A morphism of varieties  $Y \to X$ is called
an {\em alteration} if it is proper, surjective and generically finite.  This
notion of alteration will also be used for integral schemes. 

See \cite{dj}, 2.20.
\end{definition}

\begin{remark} A modification is a birational alteration.
\end{remark}

\begin{exercise} Show that an alteration $\varphi: Y \to X$ can be
factored as 
$$Y \stackrel{\pi}\longrightarrow Z  \stackrel{f}\longrightarrow X,$$
where $\pi$ is a modification, and $f$ is a finite morphism.
\end{exercise}

\begin{exercise} Suppose moreover that a finite group $G$ acts on $Y$ by
automorphisms such that 
the field of invariants $K(Y)^G$ contains the function field $K(X)$. Formulate
other factorizations of $\varphi$.  
\end{exercise}
 
\begin{remark} Given a variety $X$ and a nonzero coherent ideal sheaf
${\cal I} \subset {\cal O}_X$, the blowing up $Bl_{\cal I}X =
Proj_X (\oplus_{j\geq 0} {\cal I}^j) \to X$ gives naturally a modification,
such 
that the inverse image of ${\cal I}$ becomes invertible. If $Z\subset X$ is a
subscheme with ideal sheaf ${\cal I}$, the blowing up $Bl_Z(X)$ of $X$ along
$Z$ is defined to be the blowing up $Bl_{\cal I}X$.

 See \cite{HAG}, II.7, p. 163.
\end{remark}

\subsection{Smooth morphisms and regular varieties}
The terminology ``smooth'' will only be used in a relative situation. Thus a
morphism can be smooth. The terminology ``regular", or ``non-singular", will be
used in the absolute sense. Thus a variety can be regular. This means that for
every point 
$P$ in the variety the local ring at $P$ is a regular local ring. If a morphism
$X \to \Spec(K)$ is smooth, then $X$ is regular. It is not recommended to use
the terminology ``a smooth variety", which can be misleading and confusing. 

\subsection{Resolution, weak and strong} We state what we mean by a
resolution of singularities. There are two 
variants we will use: 

\begin{definition} Let $X$ be a variety. A {\em resolution of
singularities in the weak sense} is a modification $Y \to X$ such that $Y$ is
nonsingular.
\end{definition}
 
\begin{definition} Let $X$ be a variety. A {\em resolution of
singularities in the strong sense} is a modification $Y \to X$, which is an
isomorphism over the nonsingular locus $X_{reg}$, such that $Y$ is nonsingular.
\end{definition}

\subsection{Normal crossings}
 The following type of ``nice subschemes'' of a
variety are quite useful in desingularization problems and applications:

\begin{definition} Let $X$ be a variety. A subscheme $Z\subset X$ is
called a {\em strict normal crossings divisor} if for each point $x \in Z$,
there 
is a regular system of parameters $y_1,\ldots, y_k$ for $x$ in $X$ (in
particular the point $x \in X$ is supposed to be a regular point on $X$), such
that 
$Z$ is given on a Zariski neighborhood of $x$ by the equation $y_1\cdots y_l =
0$. 

Suppose furthermore we have a finite group acting on $Z$ and $X$ equivariantly:
$G\subset \Aut(Z\subset X)$. We say that $Z$ is a {\em $G$-strict normal
crossings divisor} if it has normal crossings, and for any irreducible
component $Z'\subset Z$, the  orbit $\cup_{g\in G}\, g(Z')$ is normal.
 
See \cite{dj}, 7.1.
\end{definition}

\n
Strict normal crossings divisors have played an important role in resolution of
singularities, and are essential in the proof of de Jong's result.

\subsection{Flatness} A crucial idea for studying ``families of
schemes'' is Serre's notion of {\em flatness} (see \cite{HAG}, III.9). 
\begin{definition}
Let $A$ be a ring and $M$ and $A$-module. Recall that  $M$ is said to be a {\em flat}
$A$-module if the functor $N \mapsto M\otimes_A N$ is exact.

A morphism of schemes $X \to Y$ is {\em flat} if at any point $x\in X$,
whose image is $y\in Y$, the local ring $\cO_{X,x}$ is a flat
$\cO_{Y,y}$-module. 
\end{definition}

\n
There are many important examples of flat morphisms which we will discuss
later. The reader is advised to consult \cite{HAG} or \cite{Matsumura} for a
more detailed discussion. The general picture should be that in a proper flat
morphism, many essential numerical invariants (e.g. dimension, degree...) are
``constant'' from fiber to fiber, so we should really think about it as a
``family''. 

\vn
Here are some instructive examples of morphisms which are {\em not} flat:

\begin{example} (See \cite{HAG}, III 9.7.1.) Let $Y$ be a curve with a node
(say, the locus $xy=0$ in the affine plane). Let $X\to Y$ be the normalization
(in the specific example, the disjoint union of two lines mapping onto  the
locus $xy=0$). Then $f:X \to Y$ is not flat. The idea one should have in mind
is that since over a general point in 
$Y$ we have one point in $X$, and over the node we have two points in $X$, this
is not really a nice family - it jumps in degree.
\end{example}

\n
The same reasoning gives a more general example:
\begin{example} Let $f:X\to Y$ be a modification. Then $f$ is flat if and only
if it is an isomorphism.
\end{example}

In particular, a nontrivial blowup is not flat.

\subsection{Stable curves} We give a formal definition of the fundamental
notion 
introduced in the introduction:

\begin{definition}
 An $S$ scheme $C\to S$ is called a {\em family of nodal curves}
 over $S$ if it is of finite presentation, proper and flat, and all
 geometric fibers are connected 
 reduced curves with at most ordinary double points (locally $xy=0$) as
 singularities. 
\end{definition}

\begin{remark}
The terminology {\em a nodal curve over $S$} can be used
interchangeably with {\em a family of nodal curves over $S$}. Indeed, if $C\to
S$ comes by way of an extension of a nodal curve $C_\eta$ over the
generic point $\eta_S$ of $S$, it may be natural to call it {\em a
nodal curve over $S$}. 
\end{remark}
 
\begin{definition} The {\em discriminant locus} $\Delta\subset S$ is the
closed subset over which $C\to S$ is not smooth.
 \end{definition}

\begin{definition}[Deligne and  Mumford] A 
family of nodal curves $$f:C\to S,$$
together 
with sections $s_i:S 
\to C , \ \ i=1,\ldots,n$ with image schemes $S_i =s_i(S)$, is called a {\em
family of stable $n$-pointed curves of genus $g$} if 
\begin{enumerate}
\item The schemes $S_i$ are mutually disjoint.
\item The schemes $S_i$ are disjoint from the non-smooth locus $\Sing(f)$.
\item All the geometric fibers have arithmetic genus $g$.
\item The sheaf $\omega_{C/S}(\sum S_i)$ is $f$-ample.
\end{enumerate}

\n
In case $n=0$ we simply call these {\em stable curves} (rather than stable
 $0$-pointed curves).
 \end{definition}

\vn
The definition is made so that a stable pointed curve has a finite automorphism
group (relative over $S$). It agrees with that made (informally) in the
introduction.  

\begin{remark}
In the litterature one sometimes finds the terminology ``$n$-poin\-ted stable
curve'' instead of ``stable $n$-pointed curve''. We try to stick to the latter,
since it effectively conveys the idea that the curve {\em with} the points is
stable. The other terminology might give the impression we are dealing with 
stable curve with some points on them. This would be a different notion in
general! 
\end{remark}

\subsection{Minimal models, existence and uniqueness} An important stepping
stone for understanding moduli of stable curves is the notion of minimal models
of 1-parameter families of curves. 

Let $K$ be a
field, and $C$ a complete, geometrically irreducible
algebraic curve 
smooth over $K$; suppose the  genus of $C$ is at least 2. Let $v$ be
a discrete valuation on $K$, and $R \subset K$ its valuation ring. Pick a
projective model ${\cal C}_0$ of $C$ over $R$. 
Following Abhyankar (1963) we can resolve singularities in dimension 2,
therefore we may assume ${\cal C}_0$ is nonsingular. Following  Shafarevich
(1966) 
and Lipman (1969) we 
have the notion of the {\it minimal model} of $C$ over $S := \Spec(R)$
 (see
Lichtenbaum,  \cite{Lichtenbaum}, Th. 4.4; also 
(see \cite{D.M}, page 87). 
We thus arrive at a family of curves ${\cal C} \to \Spec(R)$ which is a regular
2-dimensional scheme, and which is relatively minimal.
  
\begin{remark} Here we use a special case of resolution of singularities,
namely in the case of schemes of dimension 2.      
\end{remark}

            \section{Results}
\mylabel{Sec:Results}
First recall Hironaka's theorem:
\begin{theorem}[Hironaka] Let $X$ be a variety over a field $k$ of
characteristic 0. Then there exists a sequence of modifications $$X_n \to X_{n-1} \to
\cdots \to X_1 \to X_0 = X,$$ where each $X_i \to X_{i-1}$ is a
blowing up with nonsingular center, and the center lies over the
singular locus $\Sing(X)$. In particular, $X_n \to X$ is a resolution
of singularities in the strong sense.
\end{theorem}
\n
See \cite{hironaka}. Hironaka's theorem and its refinements will be
discussed elsewhere in this volume.

\n
Our main goal is to prove the following result, due to A. J. de Jong:
 
\begin{theorem} Let $X$ be a variety over an algebraically closed
field. There is a separable alteration $Y \to X$ such that $Y$ is quasi
projective and 
regular. 
\end{theorem}

\begin{corollary} Let $X/k$ be a variety. There is a finite extension
$k\subset k_1$ and a separable alteration $Y \to X_{k_1}$ such that $Y$ is
quasi 
projective and regular. 
\end{corollary}

\n 
In order for the induction in the proof to work, de Jong's theorem gives more: 
 
\begin{theorem}[de Jong]\mylabel{Th:deJongMain} Let $X$ be a variety
over an algebraically 
closed field, $Z\subset X$ a proper
closed subset. There is a separable alteration $f:Y \to X$, and an open
immersion 
$j:Y\subset \bar{Y}$, such that $\bar Y$ is projective and regular, and the
subset $j(f^{-1} Z) \cup \bigl(\bar{Y} \setmin Y\bigr)$ is a strict  normal
crossings divisor.
\end{theorem}
\n
See \cite{dj}, 4.1. The proof of this result will be given in Section
\ref{Sec:Proof}.

De Jong's theorem has a few important variants. First,  a theorem of
semistable reduction up to alteration over a one dimensional base: 
\begin{theorem}[de Jong]\mylabel{semistable-theorem} Let 
$R$ 
 be a discrete valuation ring, with 
fraction field $K$ and residue field $k$. Let $X\to \Spec R$ be an integral
scheme 
of finite type such that $X_K$ is a variety. There exists a finite extension
 $R\subset 
R_1$, where $R_1$ is a discrete valuation ring with residue field $k_1$, and an
alteration $Y \to X_{R_1}$, such that $Y$ is nonsingular, and the special fiber
 $Y_{k_1}$ is a reduced,  
strict normal crossings divisor.
\end{theorem}
\n
See  \cite{dj}, 6.5. The proof is detailed in Section
\ref{modifications}. 

This theorem is belongs to a class of theorems about
``desingularization of morphisms''. A ``dual'' case, which can
actually serve
as a building block in proving the alteration type theorems, is the
case where the base is arbitrary dimensional, and the fibers are
curves. A proof can be found in \cite{dj2}.

\begin{theorem}[de Jong] Let $\pi:X \to B$ be a proper surjective morphism
of integral schemes, with $\dim X = \dim B+1$. Let $Z\subset X$ be a proper
closed subset. There exists an alteration
$B_1 \to B$, a modification $X_1 \to \tilde{X}_{B_1}$ of the strict
transform $\tilde{X}_{B_1}$ (see Definition \ref{Def:Strict-Transform}),
sections $s_i:B_1 \to X_1$, and a proper 
closed subset $\Sigma \subset B_1$  such that
\begin{enumerate}
\item
$\pi_1:X_1 \to B_1$ is a family of pointed nodal curves,
\item
$s_i$ are disjoint sections, landing in the smooth locus of $\pi_1$,
and 
\item the inverse image $Z_1$ of $Z$ in $X_1$ is contained in the
union of  $\pi_1^{-1}\Sigma$ (the ``vertical part'') and $s_i(B_1)$
(the ``horizontal part'').
\end{enumerate}
\end{theorem}

\n
The reader who has solved the exercises in Section
\ref{modifications} will be able to complete the proof of this
theorem.  From this de Jong deduced the following refinement of Theorem
\ref{semistable-theorem}: 

\begin{theorem}[de Jong] Let $\pi:X \to B$ be a proper surjective morphism
of integral schemes, $\dim X = \dim B+r$. Assume that $B$ admits a
proper morphism to an excellent two-dimensional  scheme $S$. Then there
are alterations $B_1 \to B$ and $X_1\to \tilde{X}_{B_1}$, a
factorization $X_1 \to X_2 \to \cdots \to X_r \to X_{r+1} = B$, and
subschemes $\Sigma_i = \Sigma_i^{hor} \cup \Sigma_i^{ver}$, 
such that 
\begin{enumerate}
\item $X_i$ are nonsingular and $\Sigma_i$ are normal crossings
divisors, $i=1,\ldots,r+1$;
\item $\pi_i:X_i \to X_{i+1}$ are families of nodal curves, smooth away from
$\Sigma_{i+1}$, and
\item $\Sigma_i^{hor}$ is the union of disjoint sections of $\pi_i$,
lying in the smooth locus of $\pi_i$. 
\end{enumerate}
\end{theorem}

\n
See \cite{dj2}. Alternative proofs of different versions of this
theorem were provided in 
\cite{Abramovich} and \cite{Mochizuki}.

Next, we consider a finite group action:
\begin{theorem}[de Jong]\mylabel{equivariant-theorem}  Let 
$X$ 
be a variety over an algebraically closed field, $Z\subset X$ a proper
closed subset, $G\subset \Aut(Z\subset X)$. There is an alteration $f:Y\to X$,
and a 
finite subgroup $G_1\subset \Aut\ Y$, satisfying:
\begin{enumerate} 
\item there is a surjection $G_1\to G$ such that $f$ is $G_1 $ equivariant, and
the field extension $K(X)^G\subset K(Y)^{G_1}$ is purely inseparable;
\item $Y$ is quasi projective and nonsingular; and
\item $f^{-1}Z$ is a $G$-strict normal crossings divisor.
\end{enumerate}
\end{theorem}

n
See \cite{dj}, 7.3. The proof is detailed in exercises in Section
\ref{modifications}.  

Note that, taking $G = \{id\}$,  this implies:
 
\begin{corollary} Let $X$ be a variety over an algebraically
closed field. There is a purely inseparable
alteration $Y \to X$ where $Y$ is a quotient of a nonsingular variety by
the action of a 
finite group.
\end{corollary}

\begin{remark} For generalizations which combine both
Theorem \ref{semistable-theorem} 
and Theorem \ref{equivariant-theorem}, see \cite{dj2}.
\end{remark} 

\vn
In characteristic 0, any purely inseparable alteration is birational, and the
quotient singularities can be improved:
 
\begin{theorem}[See \cite{aj} and \cite{bp}]  Let $X$ be a variety over an
algebraically closed field of 
characteristic 0. Then there is a projective resolution of singularities in
the weak 
sense $Y \to X$.
\end{theorem}
 
\begin{remark} This is a rather weak version of Hironaka's theorem. The
point is, that new proofs, by Abramovich and de Jong \cite{aj}, and by
Bogomolov and  Pantev \cite{bp}, were given based on de Jong's ideas. The proof
by  Bogomolov and  
Pantev is extremely simple, drawing only on toric geometry. Its proof
is detailed  in Section  \ref{Sec:Bog-Pan}.
\end{remark}

\begin{question}  Can we improve the methods and obtain a weak resolution
of singularities in all characteristics? Or, at least weak resolution up to
purely inseparable alterations?
 \end{question}

\vn
The proof by
Abramovich and de Jong, detailed  in Section  \ref{Sec:A-J}, lends
itself to  generalizations in the flavor of 
de Jong's semistable reduction theorem, such as the following two results:

\begin{theorem}[Abramovich - Karu] Let $X \to B$ be a  dominant
morphism of complex projective varieties. There exists a commutative
diagram 
$$\begin{array}{ccccc}
   U_{X'} & \subset & X' &\to &X \\
  \dar & &\dar & & \dar \\
   U_{B'} & \subset & B' &\to &B
\end{array} $$
such that 
\begin{enumerate} 
\item
$X'\to X$ and $B'\to B$ are modifications,
\item 
$X'$ and $B'$ are nonsingular, 
\item 
$   U_{X'}  \subset  X'$ and
$ U_{B'}  \subset  B'$ are toroidal embeddings, and
 the morphism $X' \to B'$ is a toroidal morphism (see definition in
\ref{Sec:Toroidal}).
 
\end{enumerate}
\end{theorem}

\begin{theorem}[Abramovich - Karu] Let $X\to B$ be a dominant
morphism of complex projective varieties. There exists a commutative
diagram 
$$\begin{array}{ccccc}
   U_{X} & \subset & X_1 &\to &X \\
  \dar & &\ \ \dar \pi_1 & & \ \ \dar\pi \\
   U_{B} & \subset & B_1 &\to &B
\end{array} $$
where $B_1 \to B$ is an alteration, $X_1\to \tilde{X}_{B_1}$ is a {\em
modification} of the strict transform, $ U_{X} \subset  X_1$ and $
U_{B}  \subset  B_1$ are toroidal, the morphism $\pi_1:X_1 \to B_1$ is
toroidal with $\pi_1^{-1}U_B = U_X$, the variety $B_1$ is nonsingular
and
\begin{enumerate} 
\item the morphism $\pi_1$ is equidimensional and
\item all fibers of  $\pi_1$ are reduced.
\end{enumerate}
\end{theorem}
 
\n
See \cite{ak} for details. A refinement is given in \cite{Karu1},
and an application in \cite{Karu2}.

            \section{Some tools}
\mylabel{Sec:Tools}
In this section we gather  some basic tools which we are going to
use. Some of these
tools seem to be of vital importance in algebraic geometry, 
and it is instructive to see them functioning in the context of de Jong's
theorem. We have included some indications of proofs for the
interested reader. For the proof of the alteration theorem only the
following will be necessary: Section \ref{Sec:Strict-Transform},
Lemmas \ref{Lem:Chow},  \ref{Lem:flattening} and
\ref{Lem:Deform}, and Theorem \ref{Th:Gabber's}.

\subsection{The strict transform} \mylabel{Sec:Strict-Transform}
See \cite{dj}, 2.18.

As mentioned in the introduction, we need an operation called the ``strict
transform".  Let us recall the definition.

\begin{definition} \mylabel{Def:Strict-Transform}
Consider a morphism $X \to S$, and a base change
$T \to S$. Assume $T$ to be integral, and let $\eta \in T$ be its generic
point. Then define the {\em strict transform} $\tilde{X}_T \subset T \times _S
X$ as the Zariski closure of the generic fiber $T 
\times _S \eta$: 
$$\begin{array}{ccccccc}
\tilde{X}_T   &\stackrel{\operatorname{def}}{=} &\overline{T\times_S
\eta}^{Zar} 
&\subset & T \times_S X&  \longrightarrow & X \\  
&&& \searrow &\Big\downarrow &  &\Big\downarrow  \\
&&&&T&\longrightarrow &S.
\end{array}$$
Note that if the image of $\eta$ is not in the image of $X \to S$ (i.e. if
$T\times_S X  \to  T$ is not dominant), then the strict transform in the sense
explained here is empty. 
\end{definition} 

\begin{remark} The notion given here is different from the usual notion of the
``strict transform" of a subvariety under a modification (compare with  [HAG],
II.7, the definition after 7.15). For example consider a blowing up $T \to S$
of a surface $S$ in a point $P \in S$, and let $C \subset S$ be a curve in $S$
passing through $P$. The ``strict alteration"  (or ``strict transform" in the
terminology above) of $C$ under $T \to S$ is empty; the ``strict transform" of
$C$ under $T \to S$ in the classical sense, as explained in [HAG], II.7, is a
curve in $T$.

Some people have suggested the use of terminology ``essential pullback of $X$
along $T\to S$'', which may have some merits. After all, $\tilde{X}_T$ contains
only the ``part'' of  $T \times_S X$ which dominates $T$, which is in some
sense its essential part.
\end{remark}


\subsection{Chow's lemma} An algebraic curve and a regular algebraic surface
are 
quasi-projective. However in higher dimension an ``abstract variety" need not
be quasi-projective. A beautiful example by Hironaka (of a variety of dimension
three) is described in \cite{HAG}, Appendix B, Example (3.4.1). However
in certain 
situations (such as the alteration method described below) we like to work with
projective varieties.  

\begin{lemma}\mylabel{Lem:Chow} Given a variety $X$, there is a
modification $Y \to 
X$ such that $Y$ is quasi-projective.
\end{lemma}
\n
See \cite{RedBook} ,  I.10, p. 85, or \cite{HAG}, Exc. II.4.10 p. 107. 
\qed

\subsection{The flattening lemma} In some situations we want to replace a
morphism by a flat morphism. One
can show this is
possible after a {\it modification}  of the base. The general situation is
studied in \cite{R.G}. We only need this in an 
easier, special situation, as follows:

\begin{lemma}[The Flattening Lemma]\mylabel{Lem:flattening} Let $X$
and $Z$ be varieties over  
a perfect field $K$ (more generally,  integral schemes of finite presentation)
 and $X \to Z$ a projective, dominant morphism. There exists a modification
$f: Y \to Z$ such that the strict transform $f': \tilde{X}_Y \to Y$ is flat.
\end{lemma}


The main ingredient in the proof is the existence and projectivity of
the Hilbert scheme. Hilbert schemes were introduced and constructed by
Grothen\-dieck in \cite{FGA}, Exp. 221 (see \cite{Curves-on-Surfaces}
for simplified proofs, \cite{E.H} for discussion). We will come back
to them in  Section
\ref{Sec:Intro-to-moduli}. Their purpose is to parametrize all subschemes
of a fixed projective space. Of course, the set of all subschemes of a
projective space is rather large, so we cut it down into bounded
pieces by fixing the Hilbert polynomial $P(T)$. Grothendieck's result
may be summarized as follows:

\begin{theorem} There is a projective scheme $\cH_{\PP^N, P(T)}$ over $
\Spec \ZZ$ and a closed subscheme $\cX_{\PP^N, P(T)}\subset
\PP^N\times \cH_{\PP^N, P(T)}$ which is flat over $\cH_{\PP^N, P(T)}$,
such that $\cH_{\PP^N, P(T)}$ parametrizes subschemes of $\PP^N$ with
Hilbert polynomial $P(T)$, and $\cX_{\PP^N, P(T)} \to \cH_{\PP^N,
P(T)}$ is a universal family, in the following sense: 

Given a scheme $T$,
let $X\subset \PP^N\times T$ be a closed subscheme which is flat over $T$ and
such that the Hilbert polynomial of the fibers is $P(T)$. Then there
exists a unique morphism $h:T \to \cH_{\PP^N, P(T)}$, such that 
$$X = T \mathbin{\mathop{\times}\limits_{\cH_{\PP^N, P(T)}}}
\cX_{\PP^N, P(T)}$$ 
\end{theorem}
 
\vn
Back to the proof of Lemma \ref{Lem:flattening}.
Note that the generic fiber of $f$ is reduced. By generic flatness,
there exists a dense, open 
subset $i: U \hookrightarrow Z$ such that 
$$f_U: X_U := X\vert_U \longrightarrow U$$
is flat.  Let $P$ be  Hilbert polynomial of the fibers of $f_U$ (all fibers in
a flat family over an irreducible base have the same Hilbert polynomial), and
let ${\cal X} \to {\cal H}$  be the universal family over the Hilbert scheme
associated to this polynomial.
We have a cartesian commutative diagram: 
$$
\begin{array}{ccc}
X_U  & \longrightarrow  &  \cX \\
\downarrow  &  &  \downarrow \\
U &\stackrel{g}{\longrightarrow} & \cH.
\end{array}
$$
Note that $\cX \to \cH$ is a flat morphism. 
We have a morphism $i\times g: U \to Z \times \cH$. Define 
$$Z' := \overline{i\times g(U)}^{Zar} \quad\subset\quad Z \times \cH,$$
and let $X' \to Z'$ be the pull back: $$X' = Z'
\mathop{\times}\limits_\cH \cX.$$  Note that the 
base change of a 
flat morphism is flat, hence $X' \to Z'$ is flat. It follows from
\cite{HAG} III.9.8 that $X'$ is the
strict transform under $Z' \to Z$ of $X \to Z$. 
\qed\ref{Lem:flattening}

\begin{remark} We cited \cite{HAG} III.9.8, which is in fact an
important building block in the construction of Hilbert schemes.
\end{remark}

\begin{remark} We can delete the word ``dominant" in the flattening lemma, and
still prove the conclusion, but we do not gain much: if $X \to Z$ is not
dominant, the identity on $Z$ gives a strict transform (in the sense explained
above)  of $X$ such that $X' = \emptyset$, and flatness trivially follows. 
\end{remark}

\begin{remark} In the proof above we note a general method, which will also be
used in the question of extending curves below: suppose we study a certain
property (e.g. flatness of a map). 
Suppose there is a ``universal family'' having this property (e.g., the Hilbert
scheme). Suppose also that 
in a given family the property
holds over a dense open subset $U$ in the base. Then, after a
modification, or  an alteration of the base, depending on the situation, we can
achieve that property  by  mapping $U$
to the base of the universal family, taking the closure of the graph, and
pulling back the universal family. 
\end{remark}

\subsection{Deforming a node} An important fact underlying the role
of stable curves, which is implicitly invoked in several places in
this paper, is that a node $uv=0$ can only deform in a certain way. To
be precise:
\begin{lemma}\mylabel{Lem:Deform}
 Let $R$ be a complete local ring with maximal
ideal $m$ and algebraically closed residue field. Let $S = \Spec R$
and denote the special point by $s$. Let $X \to S$ be 
the completion of a nodal curve at a closed point $x$ on
the fiber $X_s$ over $s$, so $X_s = \Spec
({R/m[\bar{u},\bar{v}]/(\bar{u}\bar{v})})^\wedge$. Then there is an element
$f\in R$, and liftings 
$u$ of $\bar{u}$ and $v$ of  $\bar{v}$, such that $X \simeq \Spec({
R[u,v]/(uv-f)})^\wedge$. 
\end{lemma}

\vn
One can prove this using the deformation theory of a node: the versal
deformation space (see \cite{Artin}) of the completion $X_s$ of
a nodal curve has dimension 
$\dim \Ext^1(\Omega^1_{X_s},\cO_{X_s}) = 1$, and it is easy to see
that the equation $uv = t$ is versal. An elementary proof by lifting
the equation 
is sketched in \cite{dj}, 2.23.

\subsection{Serre's lemma} A critical result in the theory of moduli of
curves is, that a 1-parameter family of curves admits stable reduction after a
base change (see Theorem \ref{SR}). A crucial point in
the proof is the relationship between the 
automorphisms of a curve and the automorphisms of its jacobian, as in the
following lemma.

\begin{lemma} \mylabel{SL} Let
$C$ be a stable curve defined over an algebraically closed field $k$, let $m
\in \ZZ_{\geq 3}$, not divisible by the characteristic of $k$, and let $\varphi
\in \Aut(C)$ such that  the induced map on locally free sheaves of order $m$
$$\varphi_{\ast}: \Pic^0_C[m] \longrightarrow  \Pic^0_C[m]$$ 
is the identity map. Then $\varphi = 1_C$, the identity morphism on $C$.
\end{lemma}
\n
{\bf Proof (see \cite{S}, or \cite{Deligne}, 3.5.1).} Let
$\tilde{C} \to C$ be the normalization of $C$ (namely the disjoint
union 
of the normalizations of all irreducible components). Denote
$J:={\Pic}^0_C$ and  $X := {\Pic}^0_{\tilde{C}}$. Consider the
``Chevalley  decomposition" (as in 
\cite{C}):  
$$0 \to T \longrightarrow J  \longrightarrow X \to 0,$$ 
i.e. $T \subset  J$ is the maximal connected linear subgroup
in $J$,  the quotient is an abelian variety, and $J/T \cong
\Pic^0_{\tilde{C}}$. Note that $T \cong (\GG _m)^s$ is a split
torus. Define  
$f:= \varphi_{\ast} - 1_J \in \End(J)$.
 Using
$$\Hom(T,X) = 0$$
we obtain a commutative diagram
$$
\begin{array}{ccrcrcrcc}
0 & \to &T &\longrightarrow &J &\longrightarrow& X& \to &0 \\
  &&  g \downarrow &  &  f \downarrow && h \downarrow && \\
0 & \to &T &\longrightarrow &J &\longrightarrow& X& \to &0. \\
\end{array}
$$
By the original lemma of Serre we deduce that $h = 0$; let us sketch the
argument. The automorphism $\varphi$ is of finite order (because $C$ is
stable), hence the induced $\psi \in \Aut(\tilde{C})$ is of finite order, hence
$\psi_{\ast} = 1+h$ is of finite order.  
Note that the ring $\End(X)$ is torsion free, and since $\psi$ is of finite
order the subring $\ZZ[\psi_\ast] \subset \End(X)$ is cyclotomic. By
assumption  
the element $\psi_{\ast} 
\otimes _{\ZZ} \  \ZZ/m = 1$ in $\End(X)\otimes _{\ZZ} \  \ZZ/m$. Since
$\psi_*$ is a root of
unity, and $m \geq 3$,  this 
implies  $\psi_{\ast} = 1_X$, hence $h=0$. 

Moreover an analogous reasoning implies that $g=0$: use that $\End(T) =
\Mat(s,\ZZ)$ is torsion-free. From  $h = 0$ and $g=0$ we deduce that
$f:J \to J$
factors as 
$$J \to X \stackrel{f'}{\longrightarrow} T \to J.$$
Using  
 $$\Hom(X,T) = 0$$
we conclude $f' = 0$, hence $f=0$. Hence $\varphi_{\ast} = 1_J$, and
this implies that  $\varphi = 1_C$. 
\qed\ref{SL}

\begin{remark} We have used the fact that for $p\geq 3$, even modulo
$p$, the root of unity 
$\zeta_p$ is not equal to 1: indeed, the ring $\cO_{\QQ(\zeta_p)}/p$ is
 artinian, with generator $\zeta_p$.
\end{remark}

\subsection{Extending stable curves} Suppose we are given a stable curve
$\cC_U \to U$ over an open set $U\subset S$ of a base scheme $S$. 
Can it be extended to a stable curve $\cC \to S$? In general the answer is
negative. This question is discussed in  \cite{dJ.FO} and
\cite{Mochizuki}, where we find  criteria which ensure that in
certain cases this is possible. The general situation has the following answer%
: an extension to a stable curves is possible after an {\it
alteration} on the base. Note the difference from the Flattening
Lemma, which has to do with extending families of flat subschemes of a
fixed scheme.

\begin{theorem}[Stable Extension Theorem]
\mylabel{Th:Gabber's}  Let 
$S$ 
be a
locally 
noetherian integral  
scheme, let $U \subset S$ be an dense, open subset and let $\cC \to U$ with
sections $s^U_i:U \to \cC$ be a
stable pointed  curve. There exists an alteration $\varphi: T \to S$, and a
stable pointed 
curve $\cD \to T$ with sections $\tau_i:T \to \cD$, such that, if we write
$\varphi^{-1}(U) =: U' \subset T$, we 
have an isomorphism  
$$\cD|_{U'} \quad \stackrel{\phi}{\longrightarrow} \quad U'\quad \times_U \quad
\cC,$$ such that $\phi^*s_i^U = \tau_i$.   
\end{theorem}

\begin{remark} A proof for unpointed curves can be found in \cite{Deligne},
Lemma 1.6. We present  here a  
somewhat different proof. For simplicity of notation the proof is stated in the
case of unpointed curves. 
\end{remark}

\n
The first step is to extend {\em isomorphisms} of stable curves. The
first lemma is the following:

\begin{lemma}\mylabel{Lem:unique-stable-reduction} Suppose $T$ is the
spectrum of a discrete valuation 
ring, and $\cD \to T$ and $\cD' \to T$ are stable (pointed) curves, such
that the 
generic 
fibers  are isomorphic: $\cD_{\eta} \cong \cD'_{\eta}$. Then this
extends a unique
isomorphism: $\cD \cong_T \cD'$.
\end{lemma}

\n
For the proof see  \cite{D.M}, Lemma 1.12. The main point is
that the {\em minimal models} of $\cD$ and $\cD'$ coincide, and
$\cD$ or $\cD'$ 
are obtained from the minimal model in a unique way by blowing down
$(-2)$-curves. This lemma implies the following (see \cite{D.M},
1.11): 

\begin{lemma} Suppose $T$ is a
scheme, and $\cD \to T$ and $\cD' \to T$ are stable (pointed)
curves. Then $\Isom_T(\cD,\cD')\to T$ is finite and unramified. 
\end{lemma}

\n
Indeed, the previous lemma implies that $\Isom_T(\cD,\cD')\to T$ is
proper. Since stable curves have a finite automorphism groups, the morphism
is finite. And since stable curves have
have no nonzero vector fields, the morphism is unramified. 

As a consequence we get the following general result about extending
isomorphisms: 

\begin{lemma} Suppose $T$ is an integral {\em normal}
scheme, and $\cD \to T$ and $\cD' \to T$ are stable (pointed) curves,
such that the 
generic 
fibers  are isomorphic: $\cD_{\eta} \cong \cD'_{\eta}$. Then this induces an
isomorphism: $\cD \cong_T \cD'$.
\end{lemma}

\n
{\bf Proof:} 
The given isomorphism over the generic
point $\eta$ gives a lifting $\eta \to \Isom_T(\cD,\cD')$. The closure of its 
image in $\Isom_T(\cD,\cD')$ maps finitely and birationally to $T$. By
Zariski's Main theorem it is isomorphic to $T$, and therefore gives a section
of $\Isom_T(\cD,\cD')\to T$.\qed  

\begin{exercise} We show that the condition ``normal'' in the
previous lemma is needed. To this end, choose a regular curve $T_0$, and a
smooth curve $\cD_0 \to T_0$. Choose it in such a way that the geometric
generic 
fiber has only the identity as automorphism, and such that there exist closed
points $x, y \in T_0$ and two {\it different} isomorphisms 
$$\alpha, \beta : (\cD_0)_x \stackrel{\sim}{\longrightarrow} (\cD_0)_y.$$
Let $T_0 \to T$ be the nodal curve obtained by identifying $x$ and $y$ as a
nodal point $P \in T$ (and the curves isomorphic outside these
points). Construct $\cD_{\alpha} \to T$ by ``identifying  $(\cD_0)_x$ and
$(\cD_0)_y$ via $\alpha$''. Analogously $\cD_{\beta} \to T$. Show that  
$$\cD_{\alpha} \not\cong _T \cD_{\beta}, \quad\mbox{and}\quad
(\cD_{\alpha})_{\eta_{T}} = (\cD_0)_{\eta_{T'}} = (\cD_{\beta})_{\eta_{T}}.$$
It is instructive to describe $\Isom_T(\cD_{\alpha},\cD_{\beta})$.  
\end{exercise}

\begin{remark}
The phenomenon described in the exercise is characteristic of
situations where one has a {\em coarse} moduli space rather than a
{\em fine} one. See Definition \ref{Def:Coarse-moduli} for details. 
\end{remark}
 
\n
The following is an analogous lemma about isomorphisms of the
geometric generic fibers:

\begin{lemma}\mylabel{Lem:Extend-isom} Suppose $T$ is an integral
scheme, $\cD \to T$ and $\cD' \to T$ 
stable (pointed) curves, such that the {\em geometric} generic fibers  are
isomorphic:   
$$\cD_{\overline{\eta}} \quad\cong\quad \cD'_{\overline{\eta}}.$$
Then there exists a finite surjective morphism $T' \to T$ and an isomorphism
$$\cD \mathop{\times}\limits_T  T'  \quad\mathop{\cong}\limits_{T'}\quad
\cD'  
\mathop{\times}\limits_T  T'.$$  
\end{lemma}

\begin{remark} It is easy to give examples where the isomorphism requested does
not exist even over the generic point of $T$. 
\end{remark}

\n
{\bf Proof:} As in \cite{D.M}, 1.10 we consider $\Isom_T(\cD,\cD')$. The
condition in the lemma assures that this is not empty,  it is finite and
dominant over $T$, and the lemma follows. \qed

\n
{\bf Proof of Theorem \ref{Th:Gabber's}:}  Here we use the fact that there
exists a ``tautological family'' of curves over the compactified
moduli space of 
curves with a level structure. For stable curves without points, this
is given in \cite{vG.FO}.  For the case of 
stable pointed curves, use the results in Section
\ref{omg}.
 Another proof, more in
the line of \cite{D.M}, is sketched in
\ref{kollar's-method}

 Let us suppose that there exists $m \in \ZZ_{\geq 3}$ such that $S \to 
 \Spec(\ZZ[1/m])$. Hence the family $\cC \to U$  defines a moduli morphism 
$$f: U \to \oM[1/m] := \oM \mathop{\times}\limits_{\Spec\ZZ} \
\Spec\ZZ[1/m].$$ 
We write $M := \omM$ (after having fixed $g$ and $m$).   
We have a  curve $\cZ \to M$ such that the associated  morphism to
$\oM[1/m]$ is
 the natural morphism
$\pi: M \to \oM[1/m]$ (we say that $\cZ \to M$ is a {\em tautological
family}; see \ref{taut}). Let $U'':= U \times _{\oM} \ M$, let $U' \subset U''$ be 
a reduced, irreducible component of $U''$ dominant over $U$,  and let $\cC'$ be
the pull back $\cC' =\cC\times_U U'$. Let $\cZ' \to U'$ be the pull back of the
tautological family
$\cZ' = \cZ \times_M U'$. By 
the previous lemma we can replace $U'$ by a finite cover (call it again $U'$)
for which there is a $U'$ isomorphism $\cC' \simeq \cZ'$.  Let $S'$ be the
normalization of $S$ in the function field of $U'$. We define  
$V \subset S' \times M$ to be the image  of $U'$ by the two morphisms into $S'$
and $M$, and let $T = \overline{V}^{Zar} \subset S' \times M$.

By construction there is a stable curve over $T$, obtained by pulling back $\cZ
\to M$, which moreover by construction extends the pull back of $\cC' \to
V$. This proves the theorem in case  $S \to S_m := \Spec(\ZZ[1/m])$.   

In case $S \to  \Spec(\ZZ)$ is surjective, one does the construction for two
different values of $m$, and then one pastes the result  using Lemma
\ref{Lem:Extend-isom}.  
\qed\ref{Th:Gabber's}

\subsection{Contraction and stabilization}\mylabel{Sec:Knudsen} In
\cite{Knudsen},  II, Section 
3, 
pp. 173-179, we find a description of the following two
 constructions.  

\begin{enumerate} \item
Consider a stable $(n+1)$-pointed curve $(\cX,\cP) \to S$ with $2g-2+n >
0$. Deleting one section gives a nodal 
$n$-pointed curve (with an extra section),  which need not be a  stable
$n$-pointed curve. However, if necessary one can contract ``non-stable
components" of fibers (regular rational curves containing  not enough
singularities and marked points). After this blowing down one obtains a stable
$n$-pointed curve $(X',\cQ) \to S$, and an $S$-morphism $\cX \to \cX'$ mapping
the first $n$ sections of $\cP$ to $\cQ$. This process, which arrives at a
unique solution to this problem,  is called {\it "contraction"}. 

\item
Consider a stable $n$-pointed curve $(\cY,\cQ) \to S$ plus an extra
section $\sigma:S \to \cY$ not in $\cQ$. This extra section may meet sections
in $\cQ$, or meet the nodes of $\cY \to S$. One 
can blow up $\cY$ in such a way that the strict transforms (in the old sense)
of elements of $\cQ$ and of the extra section give a stable $(n+1)$-pointed
curve $(\cX,\cP) 
\to S$,  and an $S$-morphism $\cX \to \cY$ mapping the first $n$ sections of
$\cP$ to $\cQ$. This process, which arrives at a unique solution to this
problem, is  called {\it "stabilization"}. 
\end{enumerate}

            \section{Proof of de Jong's main theorem}
\mylabel{Sec:Proof}
 One striking feature of
the proofs of 
de Jong's theorem and its derivatives is, that all the ingredients, with  the
exception of one  subtle, but still natural, result (the Three Point
Lemma), were known and understood nearly two decades before. The way
they are put together is quite 
ingenious. 

\subsection{Preparatory steps and observations}\mylabel{reduction-steps}
 The proof of de Jong's theorem
starts with a series of simple reduction steps. 

\noindent{\bf The situation.} We want to prove de Jong's Theorem
\ref{Th:deJongMain}. Thus we are given a variety $X$ defined over an
algebraically closed field
$k$,  and a Zariski - closed subset $Z\subset X$. We   perform some
elementary reductions:

\noindent{\bf Replacing $X$ by an alteration.}  In order to prove the theorem
for a 
 variety $X$ and a closed subset $Z$, 
 it is enough to prove it for an alteration $X'$ of $X$ while
replacing $Z$ by its inverse image $Z'$ in $X'$. Thus in several stages of the
proof, once we find an alteration $X' \to X$ which we like better than $X$, we
simply replace the pair $(X,Z)$ by $(X', Z')$. 

\noindent{\bf Making $Z$ into a divisor.}  By blowing up $Z$ in $X$, and using
the 
observation above, we may assume 
that $Z$ is the support of an effective Cartier divisor. We will slightly 
abuse terminology, and say that ``$Z$ is a divisor''  when we mean that $Z$ is
a closed  subset supporting an effective Cartier divisor.

\noindent{\bf Enlarging $Z$.} Suppose $Z_i\subset X$ are  
divisors and 
$Z_1\subset Z_2$, then to prove the theorem for $(X,Z_1)$ it suffices to
prove it 
for $(X,Z_2)$. Indeed, if $f:Y\to X$ is an alteration such that $Y$ is
nonsingular and $f^{-1}(Z_2)$ is a strict normal crossings divisor, then
$f^{-1}(Z_1)$ is a Cartier divisor contained in $f^{-1}(Z_2)$, and it is clear
from the definition that it i a strict normal crossings divisor as
well. Thus we may always enlarge the divisor $Z$.

\noindent{\bf Making $X$ quasi-projective.} Using Chow's Lemma
\ref{Lem:Chow}, 
we may 
assume $X$ is quasi-projective. 
Indeed, by Chow's lemma there is a modification $X'\to X$ such that
$X'$ is quasi-projective. We may replace $X$ by $X'$. 

\noindent{\bf Enlarging $X$.} Suppose $X\subset X_1$ is an open embedding of
varieties, 
$Z_1\subset 
X_1$  a divisor which containing $X_1\setmin
X$, and $Z = X \cap Z_1$. Then 
evidently to prove the theorem for $(X,Z)$ it suffices to prove it for  $(X_1,
Z_1)$. 

\noindent{\bf Making $X$ projective.}  Since $X$ is quasi-projective,
there is an open  
embedding $X \subset \overline{X}$ where $\overline{X}$ is projective. Denote
$Z_1 = \overline{Z} \cup (\overline{X}\setmin X)$.  We may replace
$\overline{X}$ 	by the blowup of $Z_1$, thus  we may assume that $Z_1$ is the 
support of a Cartier divisor. By 
the previous observation it is enough to prove the result for $(\overline{X},
Z_1)$. 

\noindent{\bf  We may assume $X$ is normal.} Indeed, we can simply replace $X$
by its  
normalization.

\vn
To summarize, one may assume that the 
variety $X$ is projective and normal, and the subset $Z$ is the support of an
effective Cartier divisor. Moreover, one may always replace $Z$ by a larger
divisor. 

\subsection{Producing a projection}\mylabel{good-projection}
The next step is to produce a projection with some nice properties.  
We first
start with some general facts about projections in projective spaces.

Let $Y\subset \PP^N$ be a projective variety over an algebraically closed
field (in fact, separably closed would suffice). For any closed point $p\in
\PP^N \setmin Y$ we
have a projection $pr_p:Y \to \PP^{N-1}$.

\begin{lemma} Suppose $\dim Y < N-1$. Then there is a nonempty open set $U
\subset \PP^N$, such that if $p\in U$ then $pr_p$ sends $Y$ birationally to
its image. 
\end{lemma}
\n
{\bf Proof.} Let $q\in Y$ be  a regular point. Define the cone $C_{Y,q}$ over
$Y$ with vertex $q$ to be the Zariski closure of the union of all secant lines
lines containing $q$ and another $q'$, for all $q'\in Y$. It is easy to see
that  $C_{Y,q}$ 
has dimension $\leq \dim Y + 1 < N$. Note that $C_{Y,q}$ contains (as ``limit
points'') the projective tangent   space   ${\mathbb T}_{Y,q}$ at $q$. Therefore
if $p\in  \PP^N \setmin C_{Y,q}$ then
 the line through $p$ and $q$ meets $Y$ transversally, at $q$ only. This
property holds as well for the line through $p$ and $q'$, for any $q'\in Y$ in
a neighborhood of $q$. Hence the lemma. \qed

\begin{lemma} Suppose $\dim Y = N-1$. Then there is a nonempty open set $U
\subset \PP^N$, such that if $p\in U$ then $pr_p$ maps $Y$ generically \'etale
to $\PP^{N-1}$.
\end{lemma}
\n
{\bf Proof.} Same as before, using  ${\mathbb T}_{Y,q}$ instead of
$C_{Y,q}$. \qed


\n
We  go back to our $X$ and
$Z$.  

\begin{lemma}\mylabel{Lem:projection} There exists a modification
$\phi:X'\to X$ and a morphism 
$f:X 
\to 
\PP^{d-1}$ such that 
\begin{enumerate} \item There exists a finite set of nonsingular closed points
$S \subset X_{ns}$ disjoint from $Z$, such that $X'$ is the blowup of $X$ at
the points of $S$.
\item $f$ is equidimensional of relative dimension 1
\item The smooth locus of $f$ is dense in all fibers
\item Let $Z' = \phi^{-1}Z$. Then $f|_{Z'}$ is finite and generically \'etale
\item At least one fiber of $f$ is smooth.
\end{enumerate}
\end{lemma}

\n
{\bf Proof.} First project $\pi:X\to \PP^d$ using the previous lemmas
$N-d-1$ times. 

Let $B\subset \PP^d$ be the locus over which $\pi$ is not \'etale.

If we choose a general $p\in \PP^d$, then $pr_p:\pi(Z)\to \PP^{d-1}$ is
generically \'etale - simply use the lemma above for all irreducible components
of $\pi(Z)$. 

We choose such a $p$ away from  $B$. By the local description of blowing up, we
can identify the variety
$$X' = \{(x,\ell) \in X\times \PP^{d-1} | \pi(x) \in \ell\}$$
with the blowing up of $X$ at the points in $\pi^{-1}(p)$. 

We define $f:X' \to \PP^{d-1}$ to be the second projection.

We can identify the fibers: the fiber over a point $\ell$ is the scheme
theoretic inverse image $\pi^{-1}(L)$ where $L$ is the line corresponding to
$\ell$. 

It follows immediately that $f$ is equidimensional: all fibers have dimension
at most 1, and are defined by $d-1$ equations (the equations of $L$).

Since no line through $p$ is contained in $B$, every fiber has a dense smooth
locus.

The last assertion follows by Bertini's theorem, since the fibers are obtained
by intersecting $X$ with linear subspaces.

\begin{lemma} The morphism $f$ has connected fibers. \end{lemma}

\n
{\bf Proof.} Since the smooth locus is dense in every fiber, the Stein
factorization is \'etale. Since projective space has no nontrivial finite
\'etale covers, the Stein factorization is trivial.

\begin{remark} 
\begin{enumerate} 
\item The last assertion is not really necessary: if $f$ did not have
connected fibers, we could replace $f:X' \to \PP^{d-1}$ by its Stein
factorization. 
\item The projection above is the only point where it is crucial that
$X$ should be 
normal, to guarantee that the generic fiber is smooth. From here on we will
allow ourselves to make reductions after which $X$ might not be normal.
\end{enumerate}
\end{remark}

To summarize, one may assume that we have a morphism of varieties $X
\to P$, for some variety $P$, which makes $X$ into a generically smooth family
of curves, satisfying 
some nice properties, in particular  $Z \to P$ is finite and
generically \'etale.

\subsection{Enlarging the divisor $Z$}\mylabel{good-divisor}
In order to ``rigidify'' the situation, it will be useful to enlarge $Z$ so it
meets every fiber ``sufficiently''. This is done as follows:

\begin{lemma} Let $X \to P$ be as above. There exists a divisor $H\subset
X$ such that 
\begin{enumerate}
\item $f|_H : H \to P$ is finite and generically \'etale, and
\item  for  any irreducible component $C$ of a geometric fiber of $f$, we have
$$\# sm(X/P)\cap C \cap H\geq 3.$$
Here we count the points {\em without} multiplicities.
\end{enumerate}
\end{lemma}

\n
{\bf Proof.} Let $n\geq 1 $ be an integer. Given a very ample line bundle
$\cL$ on $X$, consider the embedding 
$$i: X \hookrightarrow \PP = \PP(\Gamma(X, \cL^{\otimes n}))$$ associated to $
\cL^{\otimes n}$. 

\n
{\bf Claim.} Given any irreducible curve $C\subset X$, the image
$i(C)\subset \PP$ is not contained in any linear subspace of dimension $n-1$.

\n
{\bf Proof of claim.}
Since $\cL$ is very ample, the image of $\Gamma(X, \cL)\to \Gamma(C, \cL_{|C})$
containes a base point free pencil $V \subset \Gamma(C, \cL_{|C})$. The map
$\Sym^mV\to\Gamma(C,\cL_{|C}^{\otimes n})$ has rank $\geq n+1$, therefore $
\Gamma(X, \cL^{\otimes n}) \to \Gamma(C, \cL_{|C}^{\otimes n})$ has rank $\geq
n+1$, which is what we claimed.\qed(Claim)

  We consider the collection of ``bad'' divisors and show that
there are ``good'' ones left. So consider
$$T = \{(H,y)\in \PP^\vee \times P| \dim f^{-1}y\cap H = 1\} \subset
\PP^\vee \times P. $$

It is clear that $T$ is a Zariski closed subset. We can describe the fibers of
$pr_2:T \to P$ using irreducible components of the fibers:
$$pr_2^{-1}(y) = \mathop{\bigcup}\limits_{C\subset f^{-1}y}\{H|i(C)\subset
H\}.$$ 

But by the fact that $i(C)$ is not contained in any linear subspace of
dimension $n-1$, we have
$$\codim (pr_2^{-1}(y), \PP^\vee) \geq n.$$

\n
Therefore $\dim T \leq \dim P + \dim\PP^\vee - n$.

Thus if $n$ is large enough, $pr_1(T)\subset \PP^\vee$ is of large codimension
(at least $n-\dim P$). In particular  $pr_1(T)\neq \PP^\vee$.

We  fix such large $n$. Thus there are plenty of $H$ which map finitely to
$P$.

\n
For a fixed closed point  $y\in P(k)$ consider the set
$$U(y)= \left\{ H\in \PP^\vee(k) \left| \begin{array}{l} H \protect{\not\in}
pr_1(T) 
\\ H\cap f^{-1}y \subset sm(X/P) \\ H\cap f^{-1}y \mbox{ is reduced
}\end{array}\right. \right\}.$$  
This is clearly a nonempty open set of $\PP^\vee$. Moreover, if $H \in U(y)$
then $H\in U(y')$ for all $y'$ in a neighborhood of $y$. 

If moreover $n\geq 3$, then we have that $\#H\cap f^{-1}y \geq 3$. so we are
done for all points in a neighborhood $V$ of $y$.

We deal with points in $P\setmin V$ in the same way. Using Noetherian induction
we are done.\qed(Lemma)

\subsection{The idea of simplifying the fibers} De Jong's idea is to
simplify the fibers of 
the morphism $X \to P$. Then by induction on dimension one can simplify the
base $P$, and finally put these simplifications together.

The method of simplifying the fibers uses the deepest ingredient in the
program: the theory of moduli of curves (see
Section \ref{Sec:Intro-to-moduli} for discussion).

Here is the general plan. First, as we will see below, it is easy to make an
alteration of $P$, and 
replace $X$ and $Z$ by their pullbacks, such that $Z$ becomes the union of
sections of $X \to P$. 

We can think of the generic fiber of $X \to P$ as a smooth curve with a
number of points marked on it. Say the genus of this curve is $g$, and the
number of points is $n$. By the Stable Extension Theorem
\ref{Th:Gabber's}, the generic fiber can be extended, {\em after an
alteration $P_1 \to P$}, to a family of stable curves $X_1 \to P_1$:
$$
\begin{array}{ccc} X_1 & \das & X \\
         \downarrow & & \downarrow \\
                   P_1 &\to & P \end{array}
$$

\n
The new morphism $f_1 : X_1 \to P_1$ is much nicer than $f$, since at least the
fibers are as nice as one can expect: they are nodal curves. Moreover, $Z$ was
made much nicer: it is replaced by $n$ sections which are mutually disjoint,
and pass through the smooth locus of $f_1$. 

If we can resolve $P_1$ (say using induction on dimension), then it is easy to
resolve $X_1$ as well. 

There is a problem though: 
 if we want to repeat this inductively, we cannot
allow a rational map $X_1 \das X$ which is not a morphism - since we
cannot pull back nicely along 
rational maps. So we want to find a way to make sure that $X_1 \das X$
is actually 
a regular map.

\begin{remark} If one is satisfied with proving a weaker result,
namely that every variety admits a ``rational alteration'' by a
nonsingular variety, then there is an alternative way to avoid the
issue. This is carried out by S. Mochizuki in \cite{Mochizuki}.
\end{remark}
\begin{remark} 
Another way to circumvent the issue of extending $\beta$ to a
morphism, is to ensure that it extends automatically, by using a
moduli space into which a morphism is built in: the space of stable
maps. This was carried out, in characteristic 0, in \cite{Abramovich},
Lemma 
4.2. Unfortunately the details of constructing moduli spaces of stable
maps have not yet been written out in positive or mixed
characteristics, although this would not be difficult: the results of
\cite{Behrend-Manin} imply that the moduli of stable maps forms a proper
Artin stack, and the results of \cite{Keel-Mori} imply that this
stack admits a proper algebraic space as a coarse moduli space. One
should even be able to modify the argument of \cite{Kollar},
Proposition 4.5 and show that this space is projective, but this is
not essential for the argument.
\end{remark}

\n
Let us go into details.

\subsection{Straightening out $Z$} 
\begin{lemma} There exists a normal variety $P_1$ and a separable finite
morphism $P_1 \to P$ satisfying 
the 
following property:

Let $X_1  =\tilde{X}_{P_1}$ be the strict transform (see Definition
\ref{Def:Strict-Transform}), and let $Z_1$ be the inverse image of $Z$ in
$X_1$. Then there is an integer $n\ge 3$, and $n$ distinct sections $s_i:P_1
\to 
X_1,\, i=1,\ldots,n$ such that $$Z_1 = \mathop{\bigcup}\limits_{i=1}^n
s_i(P_1).$$ 
\end{lemma}

\n
{\bf Proof.} This can be proven by induction on the degree $n$ of $Z \to P$ as
 follows: 
Let $Z_1$ be an irreducible component of $Z$ and let $P'
 := Z_1^\nu$  be its normalization. We have a generically \'etale morphism $P'
 \to 
 P$. Denote  $X' = \tilde{X}_{P'}$ and let $Z'$ be the inverse image of
 $Z$. The morphism $P'\to Z$ gives rise to a section $s_{k+1}:P'\to Z'$,
 and 
 therefore we can write  $Z' = s_{k+1}(P')\cup Z''$. We have
 $\deg(Z''\to P') = \deg(Z\to P) - 1$, and therefore the inductive
 assumption holds for $Z''$. \qed 

\subsection{Producing a family of stable pointed curves}
Let $X\to P, \quad s_i:P \to X$ be the new family. Let $U\subset P$ be an open
set satisfying the following assumptions: 

\begin{enumerate} \item $X_U\to U$ is smooth;
\item the sections $s_i|_U: U \to X_U$ are disjoint.
\end{enumerate}

Such an open set clearly exists.

 Since $n\ge 3$ this gives the morphism $X_U\to U$ the structure of a family of
stable $n$-pointed curves. 

And here comes the point where moduli theory is used: by Theorem
\ref{Th:Gabber's}, there exists an alteration $P_1 \to P$, a family 
of stable pointed curves $\cC\to P_1$, with sections $\tau_i: P_1 \to \cC$,
such that over the open set $U_1 = P_1 \times_PU \subset P_1$ we have an
isomorphism  $\beta:\cC_{U_1} \to U_1 \times_P X$, satisfying $\beta^*s_i =
\tau_i$. 

\subsection{The three point lemma} As usual, we replace $P$ by $P_1$ and $X$
by its strict transform. Thus we may assume that we have a diagram as follows:

$$\begin{array}{ccc} C & \stackrel{\beta}{\das} & X \\ 
 & \searrow & \dar \\
 & & P \end{array} $$

\n
The crucial point, for which we needed to ``enlarge $Z$'' in a previous step,
is the following: 

\begin{lemma}[Three Point Lemma] \mylabel{Lem:three-point}
 Suppose $Z$ meets the smooth locus of
every irreducible component of every fiber in at least three points. Then, at
least 
after a modification of $P$, the rational map $\beta:C \das X$ extends to 
a morphism.
\end{lemma}

\n
The proof of this lemma, which is detailed in the next few paragraphs,
is probably the most subtle point in this chapter.

\subsection{Flattening the graph}
Let $T\subset X \times_P C$ be the closure of the graph of the rational map
$\beta$. We have two projection maps $pr_1:T\to C$ and  $pr_2:T\to X$.

\begin{claim} There exists a modification $P'$ of
$P$ such that the strict transform of $X$, and the closure of the graph of
$C\das X$ are both flat. Thus we might as well assume $X\to P$ and $T\to P$ are
flat. 
\end{claim}

\n
{\bf Proof.} By the Flattening Lemma \ref{Lem:flattening} there exists a
modification $P'\to P$ such 
that $\tilde{X}_{P'}$ and $\tilde{T}_{P'}$ are both flat. Evidently the closure
of the graph 
of the rational map $C \times_PP' \to \tilde{X}_{P'}$ is contained in
$\tilde{T}_{P'}$, and since $\tilde{T}_{P'}$ is flat they 
coincide by \cite{HAG} III.9.8. \qed(Claim)

 Let $p$ be a geometric point on $P$, and denote by  $X_p, T_p, C_p$
the fibers over 
$p$. There exists a finite set $W\subset X_p$ such that $T_p\to X_p$ is finite
away from $W$. Indeed, the flatness implies that  $\dim T_p = \dim X_p = 1$.

 Thus, for any  $x\in X_p\setmin W$,  there is an open neighborhood $x \in
V\subset X$  such that $pr_2^{-1}V \to V$ is finite and birational.

 In case  $x\in \Sm(X_p)\setmin W$, 
we may choose $V\subset \Sm(X \to P)$.  Using the assumption
that $P$ is normal, it follows that $V$ is normal. In this case, by
Zariski's main theorem,   
$pr_2^{-1}V \to V$ is an isomorphism.

Note that the assumption that $x\in \Sm(X_p)$ excludes only finitely many
points, since our projection $X \to P$ is smooth at the generic point of each
component of the geometric fiber $X_p$. Therefore we conclude that the
following lemma  holds:

\begin{lemma}\mylabel{Lem:unique-component} If $X'\subset X_p$ is an irreducible
component, then there is 
a {\em unique} irreducible component $T'$ of $T_p$ mapping finitely
onto $X'$ via $pr_2:T\to X$. Moreover, $T' \to X'$ is birational.
\end{lemma}

\n
Repeating the argument for $pr_1:T\to C$, we also have:

\begin{lemma}\mylabel{Lem:unique-component-C} If $C''\subset C_p$
is an irreducible component, then there is 
a  unique irreducible component $T''$ of $T_p$ mapping finitely onto
$C''$ via $pr_1:T\to C$. Moreover, $T'' \to C''$ is birational.
\end{lemma}

\n
Let $x\in \Sm(X_p)$ be a closed point. Considering the Stein
factorization $T\to \tilde{X}\to X$, we have that the fiber $pr_2^{-1}(x)$ is
connected. Indeed, since $X$ is normal at $x$,  we have that $\tilde{X}\to X$
is an isomorphism at $x$.

\subsection{Using the three point assumption} Let $X'\subset X_p$ be an
irreducible 
component, and $T' \subset T_p$ the unique component mapping finitely
(and birationally) onto it,
as in  Lemma \ref{Lem:unique-component} above. We will prove that
$pr_1:T'\to C$ 
is non-constant. Assume by contradiction that $pr_1(T') = c$ is a point. 

We will  use the three point assumption. Let $s_i:P\to X, i=1,\ldots,3$
be three of the given sections such that $s_i(p) = x_i$ are three distinct
points on $\Sm(X')$. Let $T_i = pr_2^{-1}x_i$. Let $\tau_i(p) = c_i\in C_p$.

 Note that the point $t_i=(c_i,x_i) \in C\times X$ is in $T$.

Assume $c\not\in \{c_i,i=1,\ldots,3\}$. Then each of
$T_i,i=1,\ldots,3$ contains an irreducible component $T_i'$ whose image in $C$
is again a curve passing through $c$. These image components are {\em
distinct}. Indeed, $T_i$ are disjoint subschemes of $T_p$, whose images in $C$
connect $c$ with $c_i$, and therefore each has an irreducible component whose
image contains $c$. These components are distinct, and by
Lemma \ref{Lem:unique-component-C} their images are distinct.

This 
contradicts the fact that $C_p$ is nodal. Thus $c$ is among the $c_i$.

 Assume, without loss of generality, $c=c_1$. Repeating the argument of the
previous paragraph we conclude that there are two distinct components of $C_p$
passing through $c$.
This contradicts the fact that $C_p$ has a marked point at $c=c_1$. 

Thus we  conclude that $pr_1:T\to C$ is finite and birational.

By  Serre's criterion   $C$ is normal: it is clearly regular in codimension 1,
and condition $S_2$
follows since $C\to P$ has reduced one-dimensional fibers and $P$ is normal.

We conclude that $T\to C$ is an isomorphism, hence $\beta$ extends as a
morphism! \qed\ref{Lem:three-point}

\subsection{Induction}

We arrived at the following situation:
$$\begin{array}{ccc} C & \stackrel{\beta}{\to} & X \\ 
 & \searrow & \dar \\
 & & P \end{array} $$
\n
We may replace $X$ by $C$, and $Z$ by its inverse image in $C$. Note that 
$Z$ is no longer finite over $P$: it has a``finite part'', the union of the
sections $\tau_i: P \to C=X$, but there is a ``vertical'' part $Z_{vert}$,
which is the union of irreducible components of {\em singular} fibers of $X  
\to P$.  

Let $\Sigma\subset P$ be the closed subset over which $f:X \to P$ is not
smooth. By the induction assumption there is a projective alteration $P_1 \to
P$ 
such that $P_1$ is nonsingular and the inverse image of $\Sigma$ is a strict
normal crossings divisor. We may replace $X$ by its pullback to $P_1$, and
replace $P$ by $P_1$. It is 
 convenient to replace $Z$ by its union with $f^{-1}(\Sigma)$.

We  arrived at a situation where both $P$, and the morphism
$f:X \to P$, are simplified. The resulting variety has very simple
singularities, 
and its desingularization results from the following exercises.

\subsection{Exercises on simple blowup of nodal families}\label{resolve-nodal}
Let $f:X \to S$ be a flat morphism of varieties over an algebraically closed
field $k$, with $n=\dim X = \dim S + 1$. Let $D\subset S$ be a reduced
divisor. We make the following assumptions.
\begin{itemize} 
\item[N1] The base $S$ is nonsingular.
\item[N2] The divisor $D$ has strict normal crossings.
\item[N3] The morphism $f$ is smooth over $S\setmin D$.
\item[N4] The morphism  $f: X \to S$ is a nodal curve.
\end{itemize}
Let $x\in X$ be a closed point and let $s = f(x) \in S$. By assumption we may
choose a regular system of parameters $t_1,...,t_{n-1}$ at $s$ such that $D$
coincides on a neighborhood with the zero locus of $t_1\cdots t_r$ for some
$r\leq n-1$. It can be seen that if $x$ is a singular point of $X$, then the
completed local ring of $X$ at $x$ can be described as
$$(*) \quad k[[u,v]] / (uv - t_1^{n_1}\cdot\cdots \cdot t_r^{n_r})$$.
 
\n
{\bf Step 1:} Assume $\codim_X\Sing(X) = 2$.
\begin{enumerate}
   \item Show that there is an irreducible component $D_1\subset D$ and
   $\Sigma_1\subset \Sing(X)$ such that $f(\Sigma_1) = D_1$.
   \item Fix a point $x\in \Sigma_1$, and use formal coordinates as in $(*)$,
   such that $D_1= V(t_1)$. Show that the power $n_1$ of $t_1$ is $>1$.
   \item Show that the ideal of $\Sigma_1$ in the formal completion is
   $(u,v,t_1)$. 
   \item Conclude that $\Sigma_1 \to D_1$ is \'etale, in particular  $\Sigma_1$
   is nonsingular. 
   \item Let $X_1 = Bl_{\Sigma_1}X.$ Show that $X_1 \to S$ satisfies conditions
   N1-N4, there is at most one component of $\Sing(X_1)$ over $\Sigma_1$, with
   the exponent $n_1$ replaced by $n_1 - 2$    
   \item Conclude by induction that there is a blowup $X' \to X$ centered above
   $\Sing(X)$, such that $X'$ satisfies N1-N4, and $\codim_X\Sing(X) > 2$.
   \item Show that  each component of $\Sing(X')$ is defined by
   $u=v=t_i=t_j$ in equation $(*)$, in particular it is nonsingular. 
\end{enumerate}

 \n
{\bf Step 2:} Assume $\codim_X\Sing(X) > 2$.
Define $Z = f^{-1} D$. Unfortunately here we need to abandon the structure
$X\to B$ of a family of nodal curves. Instead we look at $X$ itself. The
situation is as follows:
\begin{itemize} 
\item[T1] whenever $x$ is a nonsingular point of $X$, $Z$ has normal
crossings at $x$.
\item[T2] whenever $x\in \Sing(X)$, we have formal description
$$(**) \quad k[[u,v]] / (uv - t_1\cdots  t_s),\quad 2\leq s \leq r\leq n-1$$
and $Z = V(t_1\cdots  t_r)$.
\item[T3] All components of $\Sing(X)$ are nonsingular.
\end{itemize}
\begin{enumerate}
   \item Let $E\subset \Sing(X)$ be an irreducible component. Show that the
   blowup $Bl_EX$ satisfies T1-T3, and its singular locus has one fewer
   irreducible component.
   \item Conclude by induction that there is a resolution of singularities
   $X'\to X$.
\end{enumerate}

This concludes the proof of Theorem \ref{Th:deJongMain}\qed

%
	\section{Modifications of the proof for Theorems
	   \ref{semistable-theorem} and
	   \ref{equivariant-theorem}}
\mylabel{modifications}  
%
\subsection{Exercises on removing the conditions on the
projection}\mylabel{Ex:remove-conditions} 
An important step in the proof of de Jong's theorem was, that given the
projection $X \to P$, one could construct an alteration $P_1 \to P$ and a 
diagram
$$\begin{array}{ccc} C & \stackrel{\beta}{\to} & X \\ 
 \dar &  & \dar \\
P_1 & \to & P \end{array} $$
where $C\to P_1$ was a family of nodal curves. In order for the proof to go
through, we made several assumptions on the projection $X \to P$. Here we will
show that even if these conditions fail, we can still reduce to the case where
they do hold.

\begin{exercise}
Using an alteration, show that the condition  that $Z \to P$ be finite
in Lemma \ref{Lem:projection} (4) is
unnecessary for the rest of the proof.
\end{exercise}

\begin{exercise}
 Show that, if one is willing to accept inseparable alterations in the
theorem, the condition that $Z\to P$ be generically \'etale in Lemma
\ref{Lem:projection} (4) is unnecessary for the rest of the proof.
\end{exercise}

\begin{exercise}{\LARGE${\mathbf *}$} By reviewing the arguments, show that the
condition 
that every component of every fiber of 
$X\to P$ be generically smooth is unnecessary. 

\n
Here a modification of the three point lemma is be necessary! In \cite{dj2}, de
Jong uses a trick of ``raising the genus of the curves'' with finite
covers. Another way goes as follows: in the proof of the Three Point Lemma
\ref{Lem:three-point}, after flattening $X$ and $T$, one works with fibers of
{\em the normalizations} $ X^{\nu}$ and $T^{\nu}$. This way one avoids the need
for $\Sm(X_p)$ to be dense. One notes that the sections
$s_i$ lift to  $ X^{\nu}$, and at least three meet every component of every
fiber, since $Z$ is the support of a Cartier divisor! The details are left to
the reader.
\end{exercise}

\begin{exercise} Using the flattening lemma and the previous exercise, show
that the condition that $X \to P$ be equidimensional is unnecessary.
\end{exercise}

\begin{exercise} Show that, if one is willing to accept inseparable alterations
in the theorem, the condition that the generic fiber of $X\to P$ be smooth is
unnecessary. 
\end{exercise}

\subsection{Exercises on Theorem \ref{semistable-theorem}}
Let us address Theorem \ref{semistable-theorem} on semistable reduction up to
alteration. Suppose
$S=\Spec R$ where $R$ is a discrete valuation ring,  $X \to S$ is a morphism
as in the theorem, and $Z$ a proper closed subset.

\begin{exercise} Show that we may assume $X$ projective over $S$, that $Z$ is
the support of a Cartier divisor.
\end{exercise}

\begin{exercise} Show that we may assume that the generic fiber is a normal
variety, and that $X$ is a normal scheme. (You may need an inseparable base
change!) 
\end{exercise}

\begin{exercise} Let $d$ be the dimension of $X_\eta$.  Produce a projection
$X\to \PP^{d-1}_S$ with connected fibers.
\end{exercise}

\begin{exercise} Use the semistable reduction argument, with the Three Point
Lemma, and the results of Section \ref{Ex:remove-conditions} to
replace $X$ by a nice 
family of curves $X \to P \to S$. 
\end{exercise}

\begin{exercise} Use induction on the dimension to conclude the proof of the
theorem.  
\end{exercise}

\begin{exercise} Can you think of other situations where a similar theorem can
be proven, 
where $S$ is not necessarily the spectrum of a discrete valuation ring? (This
is interesting even in characteristic 0!)
\end{exercise}

\subsection{Exercises on Theorem \ref{equivariant-theorem}}
We  address the equivariant version of the theorem. Suppose $X$ is a
variety, $Z$ a proper closed subset and $G$ acts on $X$ stabilizing $Z$. We
wish to prove Theorem \ref{equivariant-theorem}.

\begin{exercise}  Produce an equivariant version of Chow's lemma, so that we
may assume $X$ 
is projective.
\end{exercise}

\begin{exercise}  Show that, to prove the theorem, it suffices to
consider the case where $Z$ is a divisor.
\end{exercise}

\begin{exercise}  Show that we may replace $Z$ by a bigger {\em equivariant}
divisor; in 
particular we may assume $Z$ contains the fixed point loci of elements in $G$.
\end{exercise}

\begin{exercise}  Using a projection of $X/G$, show that we may assume we have
an 
equivariant projection $X \to P$ making $X$ into a nice family of curves.
\end{exercise}

\begin{exercise}  Consider the case $X = \bfA^2_k$ where $\chara k = p$, and
$G=\ZZ/p$ 
acting via $(x,y) \mapsto (x, x+y)$. Show that the fixed point set maps
inseparably to the image. In particular, the map $Z\to P$ in the previous
exercise might be inseparable!
\end{exercise}

\begin{exercise}  Making an alteration ``Galois'': Given a variety $W$, a
finite group 
action $H \subset \Aut\ W$ and an alteration $V_0 \to W$, show that there exists
an alteration $V \to V_0$, and a finite group $H'$ with a surjection $H'\to H$,
and a lifting of the $H$  action $H'\subset \Aut (V\to W)$ such that 
the extension of fixed fields $K(W)^H \subset K(V)^{H'}$ is purely
inseparable. 
\end{exercise}

\begin{exercise}{\LARGE$\mathbf *$} Use the uniqueness in the stable reduction
theorem to show that there 
is an alteration $P'\to P$, a family of stable pointed curves $C \to P'$ and a
finite group $G'$ with a quotient $G'\to G$ and a diagram
$$\begin{array}{ccc} C &\to& X \\ \dar & & \dar \\ P' &\to& P\end{array}$$
on which $G'$ acts equivariantly,  such that $C$ is birational to
$\widetilde{X_{P'}}$, and the extensions $K(X)^G\subset  K(C)^{G'}$ and 
  $K(P)^G\subset  K(P')^{G'}$ are purely inseparable.
\end{exercise}

\begin{exercise}  Use induction on the dimension and a suitable modification of
the 
elementary blowups argument to conclude the theorem.
\end{exercise}
 
            \section{Toroidal geometry}
\mylabel{Sec:Toroidal}
\subsection{Basic definitions} For simplicity we work over an algebraically
closed field.
We recall the notion of a toric variety (a more thorough discussion is
available in \cite{cox}):

\begin{definition} A variety $X$ together with an open dense embedding
$T\subset 
X$ is called a {\em toric variety} if $X$ is normal, $T$ is a torus
(geometrically 
isomorphic to ${\mathbb G}_m^k$), and the action
of $T$ on itself by translations extends to an action on $X$.
\end{definition}

\n
To get an intuitive idea about the singularities of a toric variety, it is
worth noting that  a normal,  affine variety, defined by equations between
monomials is toric, and every toric variety is locally of this type. 

For many purposes toric varieties are too restrictive. A more general notion
was introduced by Mumford in \cite{kkms}:

\begin{definition} A variety $X$ together with an open embedding $U\subset
X$ is called a {\em toroidal embedding} if any point $x\in X$ has an \'etale
neighborhood $X'$ such that $X'\times_XU\subset U$ is isomorphic to an \'etale
neighborhood of a point on a toric variety. 
\end{definition} 

\n
Thus a toroidal embedding looks locally like a toric variety, and the big open
set $U$ is a device which ties together these ``local pictures''. In a sense,
this notion is suitable for studying varieties whose singularities are like
those 
of toric varieties. 

In this section we recall facts about toric varieties and briefly
indicate how one can obtain analogous facts about toroidal
embeddings. The details are available in \cite{kkms}.

\begin{remark} A more sheaf - theoretic approach was introduced by K. kato, see
\cite{kato}, \cite{kato1}.
\end{remark}

\begin{definition} 
A toroidal embedding is said to be {\em strict} (or a {\em
toroidal embedding without self intersections}) if every irreducible component
of $X\setmin U$ is normal. 
\end{definition}

\n
For instance, if $X$ is a nonsingular variety, $D\subset X$ is a strict
normal crossings divisor, and $U = X \setmin D$, then $U\subset X$ is a strict
toroidal embedding. 

{\em We will only work with strict toroidal embeddings}.

\begin{definition} 
If $G\subset \Aut(U\subset X)$ is a finite group, we say that $G$ {\em acts
toroidally} if for any point $x\in X$, the stabilizer $G_x$ of $x$  can be
identified with a subgroup of the torus in an appropriate \'etale
neighborhood of $x$.
\end{definition}

\begin{definition}
A morphism between toric varieties is called a {\em toric morphism}, if it is
surjective and torus-equivariant. A morphism of toroidal embeddings
$(U_X\subset X) \to (U_Y \subset Y)$ is called a {\em toroidal morphism} if
locally on $X$ it looks like a toric morphism.
\end{definition}

\subsection{The cone} 

 First recall  some notation (see \cite{cox}):

 $M= \Hom(T, \GG_m)$ - this is the group of algebraic characters of $T$; 

$M_\RR = M\otimes\RR$

$N = \Hom(M, \ZZ) = \Hom( \GG_m,T)$  - this  is the group of 
  1-parameter subgroups on $T$;

$N_\RR = N\otimes\RR$

It is common to call the functions  defined by elements of $M$ the {\em
monomials}. One uses the notation $x^m$ for the monomial associated with the
element $m\in M$.

Recall the basic correspondence 
between 
\begin{center} $\{$affine toric varieties $T\subset X\}$  \\
 and  \\
  $\{$strictly convex rational polyhedral cones $\sigma \subset N_\RR\}$
\end{center} 
which can be  defined in one direction via
$$V_\sigma = \Spec k[\sigma^\vee\cap M],$$
 and in the other direction by 
 $$\sigma\ = \ \hbox{\parbox{9.5cm}{ the cone spanned by the 1-parameter
subgroups $\phi:\GG_m\to T$  such that the limit
$\lim_{z\to 0} \phi(z)$ exists in $V$.}}$$

There is another, less well known characterization of  $\sigma$, which is less
dependent on the torus action, and is therefore useful for  toroidal
embeddings: 

Any  monomial $m\in
M$ defines a Cartier divisor $\Div(x^m)$ supported on $V\setmin T$. If
$\sigma^\vee$ contains a line through the origin, then for any $m$ on this line
the divisor is easily seen to be trivial (both $m$ and $-m$ give regular
functions).

We use the following notation:

$\sigma^\perp = \{m\in M_\RR| \langle m,\sigma\rangle = 0\}$

 $M^\sigma=$  Cartier divisors  supported on $V\setmin T$. 

One can easily see that $M^\sigma = M / \sigma^\perp\cap M$.

 $N_\sigma=span(\sigma).$ Clearly $N_\sigma = \Hom(M^\sigma , \ZZ)$.

Let $M^\sigma_+\subset M^\sigma$:  the {\em effective} Cartier divisors. 

We have that $M^\sigma_+ = \sigma^\vee\cap M / \sigma^\perp\cap M$. 

 It is not hard to see that $\sigma = (M^\sigma_+)_\RR^\vee$, the dual
cone of the cone spanned by $M^\sigma_+$.

In short: $\sigma$ is the dual cone to the cone of effective Cartier divisors
supported on $V\setmin T$.


\subsection{The toroidal picture}  We wish to mimic the construction
of the cone
in the toroidal   case. 
We follow \cite{kkms}, Chapter II.

Let $U \subset X$ be a
strict 
toroidal embedding. $X \setmin U  = \cup D_i, $ where $D_i $ normal.

We decompose $\cap_{i\in I} D_i = \cup X_\alpha$; the locally closed subsets
$X_\alpha$ are called {\em strata}.  

Each stratum has its star: $\Star(X_\alpha)= \cup_{X_\alpha \subset
\overline{X_\beta}} X_\beta$.  

Note: $X_\alpha$ is the unique closed stratum in $\Star(X_\alpha)$. In a sense
it is analogous to the unique closed orbit in an affine toric variety.

  Define:

 $M^\alpha=$ group of Cartier divisors supported on $\Star(X_\alpha)\setmin U$;

$M^\alpha_+=$  subset of effective Cartier divisors;

 $N_\alpha=\Hom(M_\alpha,\ZZ)$;

 $\sigma_\alpha = (M^\alpha_+)_\RR^\vee$.

Thus, to each stratum we associated a strictly convex rational polyhedral cone.

\begin{remark} The cone $\sigma_\alpha$ has a description analogous to the
toric one using 
1-parameter subgroups, in terms of valuations. Let $RS(X)$ be the discrete
valuations on $X$. Let $v$ be a valuation centered in $\Star(X_\alpha)$. Let
$f_j$ be rational function defining generators of $M_\alpha$ on a small affine
open. Then $v(f_j)$ is a vector in $\sigma_\alpha$, and in fact $\sigma_\alpha$
can be described as a set of equivalence classes of discrete valuations
centered in $\Star(X_\alpha)$, the equivalence being defined by equality of the
valuations of these functions $f_j$.
\end{remark}

\subsection{Birational affine morphisms} 
\mylabel{Sec:tor-birational}
Recall: if $\tau\subset \sigma$ are two 
strictly convex rational polyhedral cones, then $\tau^\vee\supset \sigma^\vee$
gives rise to a morphism  
$V_\tau\to V_\sigma$, which is birational and affine. 

Note that $V_\tau\to V_\sigma$ can be described in the following invariant
manner:

$$ V_\tau=  \Spec_{V_\sigma} \sum_{E\in M^\tau_+} 
\cO_{V_\sigma}(-E), $$
where the sum is taken inside the field of rational functions of
$V_\sigma$.

This clearly works over $\Star(X_\alpha)$ in the toroidal case as well.

\subsection{Principal affine opens} If $m\in \sigma^\vee$ then $\tau = \{n\in
\sigma | <n,m>=0\}$ is a face of $\sigma$. We have $\tau^\vee = \sigma^\vee+
\RR \cdot m$, and therefore $V_\tau$ is the principal open set on $V_\sigma$
obtained by inverting the monomial $x^m$.

Again, this can be described divisorially in terms of $\Div(m)$. Thus the same
is true for $\Star(X_\alpha)$: given a face $\tau$ of $\sigma_\alpha$, we get an
open set $$\Star(X_\beta)\subset \Star(X_\alpha)$$ such that $\tau =
\sigma_\beta$. 

The most important face of a cone is the vertex. It corresponds to the open set
$T\subset V$. In the toroidal case you get $U$.

\subsection{Fans and polyhedral complexes}
Recall: if $\sigma_1$ and $\sigma_2$ intersect along a common face $\tau$, then
$V_{\sigma_1}$ and $V_{\sigma_2}$ can be glued together along the common open
set 
$V_\tau$, forming a new toric variety.

In general, whenever you have a {\em fan} $\Sigma$ in $N$, namely a collection
of cones $\sigma_i$ intersecting along faces, you can glue together the
$V_{\sigma_i}$  and get a toric variety $V_\Sigma$.

It is not hard to see that {\em every} toric variety is obtained in this way in
a unique manner. The   point is that every toric variety is covered by affine
open toric varieties.

In the toroidal case, 
$X$ is covered by the open sets $\{\Star(X_\alpha)\}_\alpha$. 

In general  $$\Star(X_{\alpha})  \cap \Star(X_{\beta}) = \cup
\Star(X_{\gamma_i}),$$ so $\sigma_{\gamma_i}$ are possibly several faces of both
$\sigma_\alpha$ and  $\sigma_\beta$. 

Still these can be glued together, as a rational conical polyhedral complex.
The main difference from the toric case, is that it is abstractly defined, and
in general it is not linearly contained in some vector space $N_\RR$.

\subsection{Modifications and subdivisions} Let $\Sigma$ be a fan, and
$\Sigma'\to \Sigma$ a (complete) subdivision. This corresponds to a toric
modification $V_{\Sigma'} \to V_\Sigma$. 

Since the construction is local (the $\Spec$ construction, as in
Section \ref{Sec:tor-birational}, and gluing) it works  
word for word in the toroidal case. There is a small issue in checking that the
resulting modification  is still a strict toroidal embedding; this is discussed
in detail in \cite{kkms}.

In \cite{kkms} (see also \cite{cox}) it is shown that a modification
is {\em projective} if and only if the subdivision is induced by a
{\em support function} - one associates to a support function an
ideal, whose blowup gives the modification. This works in the toroidal
case as well.

\subsection{Nonsingularity} Recall: an affine toric variety $V_\sigma$ is
nonsingular if and only if $\sigma$ is simplicial, generated by a basis of
$N_\sigma$ (namely, part of a basis of $N$). Such a cone is called 
nonsingular.   

In general: a toric variety $V_\Sigma$ is nonsingular if and only if every
cone $\sigma\in \Sigma$ is nonsingular. 

This is a local fact, so it is true in the toroidal case as well.

\subsection{Desingularization} Recall that it is easy to
resolve singularities of a toric variety: one finds a simplicial subdivision
such that every cone is nonsingular.

Obviously, the same works in toroidal case! We obtained:

\begin{theorem} For any toroidal embedding $U\subset X$ there is a
projective toroidal modification $U\subset X' \to X$ such that $X'$ is
nonsingular.
\end{theorem}
\n
See \cite{kkms},  Theorem 11*, page 94.

\subsection{Exercises on toric varieties and toroidal embeddings}
\begin{enumerate}
\item Show that $\GG_m^n \subset \AAA^n$ is a toric variety. Describe its
cone.  
\item Show that $\GG_m^n \subset \PP^n$ is a toric variety. Describe its
fan.  
\item Let $X \subset \AAA^n$ be a normal variety defined by monic monomial
equations of type $$\prod x_j^{n_j} =\prod x_j^{m_j}.$$ 
Show that $X$ is toric. (Identify the torus!)
\item  Do the same if the monomial equations are not necessarily  with
coefficients $=1$. 
\item Describe the cone associated to the affine toric variety defined by $$xy
= t_1^{k_1}\cdots t_r^{k_r}.$$
\item Look at the affine 3-fold  $xy = zw$. Let $X'\to X$ be the blowup of $X$
at the ideal $(x,z)$. Describe this blowup, show that it is toric, and describe
the cone subdivision associated to it.
\item Let $X = \AAA^2$, $D=\{xy(x+y-1)=0\}$, $U = X\setmin D$. Show that
$U\subset X$ is a toroidal embedding. Describe its conical polyhedral complex.
(Compare with the fan of $\PP^2$!)
\item  Do the same for $D = \{y(x^2+y^2-1)=0\}$. Show that the resulting
complex can  not be  linearly embedded in a vector space.
\item Consider the surface $X = \{z^2 = xy\}$, $U = \{z\neq 0\}$. Show that
$U\subset X$ is toric and describe its cone.
\item Consider the surface $X$ above, let $D_1 = \{x=0\}$, $D_2 = \{y=x(x-1)^2
\mbox{ and  } z = x(x-1)\}$. Let $U = X\setmin (D_1\cup D_2)$. Show that
$U\subset X$ is toroidal. Describe its conical polyhedral complex. Make sure to
describe the integral structure!
\end{enumerate}

\subsection{Abhyankar's lemma in toroidal terms}
Abhyankar's lemma about fundamental groups (see \cite{g},
\cite{gr-mu}) describes the local 
tame fundamental group of a variety around a normal crossings divisor. Let
$X = \Spec k[[t_1,\ldots,t_n]]$ and let $D = V(t_1\cdots t_n)$. Let $Y \to X$
be a finite alteration which is tamely branched along $D$, and \'etale away
from $D$. For $m$ prime to $\chara\ k$, denote $X_m = \Spec
k[[t_1^{1/m},\ldots,t_n^{1/m}]]$. Abhyankar's lemma says that the normalization
of $Y \times_XX_m$ is \'etale over $X_m$.

In the following exercises we interpret this in toroidal terms.
\begin{exercise}
Let $U \subset X$ be a {\em nonsingular} strict toroidal embedding. Let $f:Y\to
X$ be a 
finite cover, which is tame, and \'etale over $U$. Then $f^{-1}U \subset Y$ is
a strict toroidal embedding.
\end{exercise}

\begin{exercise}
Suppose further that $Y \to X$ is Galois, with Galois group $G$. Show that $G$
acts toroidally on $Y$.
\end{exercise}

            \section{Weak resolution of singularities I}
\label{Sec:Bog-Pan}
Given the existence of toroidal resolution, the proof of weak
resolution of singularities in 
characteristic 0 
by Bogomolov and Pantev is arguably the simplest  available. It does not even
require surface resolution, or moduli spaces.

We will  go through this proof. The steps of proof here include some
simplifications on the arguments in \cite{bp}, which came up in discussions
with T. Pantev. These and additional 
simplifications were discovered independently by K. Paranjape
\cite{Par}, and we have   
used his exposition in some of the following exercises.

\subsection{Projection} Let $X$ be a variety over an algebraically closed
field 
of 
characteristic 0, 
and $Z\subset X$ a proper closed subset. Let $n=\dim X$, and again assume
we know the weak resolution theorem for varieties of dimension $n-1$.

First a few reduction steps:
\begin{enumerate} 
\item Show, as in \ref{reduction-steps} that we may assume $X$ projective and
normal, and $Z$ the support of 
a Cartier divisor.
\item Show that there is a finite projection $X \to
\PP^n$.
\item Let $P\to  \PP^n$ be the blowup at a closed point. Show that 
$$P \simeq
\PP_{\PP^{n-1}} (\cO_{\PP^{n-1}} \oplus \cO_{\PP^{n-1}}(1)).$$
Denote by $E$ the exceptional divisor of $P\to  \PP^n$.
\item By blowing up  a general point on $\PP^n$, and blowing up $X$ at the
points above, show that we may assume we have a finite morphism $f:X \to P$,
which is \'etale along $E$, such that the image of $Z$ is disjoint from $E$,
and maps finitely to $\PP^{n-1}$.
\item By the Nagata - Zariski purity theorem, note that the branch locus of $X
\to P$ is a divisor $B$, disjoint from $E$, mapping finitely to $\PP^{n-1}$. 
\end{enumerate}
We  replace $Z$ by $Z\cup f^{-1}B$.

\subsection{Vector bundles}
The next  steps are aimed at replacing $P$ by another $\PP^1$-bundle 
$Q \to \PP^{n-1}$, such that the branch locus in $Q$  of $X\to Q$ becomes
simpler. 
Let $Y$ be any variety, $F$ a rank-2 vector bundle on $Y$, $P= \PP_Y(E)$.
Let $E\subset P$ be a divisor which is a section of $\pi:P\to Y$ and let
$D\subset 
P$ be another effective divisor disjoint from $E$. Denote by $\cO_P(1)$ the
tautological bundle, and by $d$ the relative degree of $D$ over $Y$.

\begin{enumerate} 
\item Consider the exact sequence $$0 \to \cI_D(d) \to \cO_P(d) \to \cO_P(d)|_D
\to 0. $$ Use this to show that there is an invertible sheaf $\cL_D$ on $Y$
such that 
$\cI_D(d)\simeq \pi^*\cL_D$.
\item If $D_1, D_2\subset P$ are any two {\em disjoint} divisors finite of
degree $d$ over $Y$, show that there is an embedding of vector bundles
$\cL_{D_1} \oplus \cL_{D_2} \subset sym^d F$ inducing a {\em surjection}
$\pi^*(\cL_{D_1} \oplus \cL_{D_2})\to \cO_P(d).$ 
\item Assume the characteristic is 0. Consider the case $D_1 = dE, D_2 =
D$. show that the resulting morphism  $P\to \PP_Y(  \cL_{D_1} \oplus
\cL_{D_2})=P'$ maps $E$ to a section $E'$ and $D$ 
to a disjoint section; and its branch locus is of the form $(d-1)E + D'$ where
$D'$ has degree $d-1$ over $Y$ and is disjoint from $E'$.
\item Continue by induction to show that there is a $\PP^1$ bundle $Q\to Y$
and a morphism $g:P\to Q$ over $Y$, such that the image of $D$ and the branch 
locus of $g$ form a union of sections of $\pi_Q:Q \to Y$.
\end{enumerate}

\subsection{Conclusion of the proof}
Back to our theorem, where $Y = \PP^{n-1}$. Composing with the morphism $f:X
\to P$, we obtain that the image $g(f(Z))\subset Q$ is the image $D_1$ of a
section 
$s_1:\PP^{n-1} \to Q$  of $Q \to
\PP^{n-1}$ and the branch locus of $g\circ f$ is the union of images $D_i$
sections $s_i:\PP^{n-1} \to Q$  as well. Denote  $\Delta = \pi_Q(\cup_{i\neq
j}D_i\cap D_j)$. 

The following steps use moduli theory; however it has been shown (in the
preprint version of \cite{bp}, and in Paranjape's exposition
\cite{Par}) that this can 
be easily
circumvented. 

\begin{enumerate}
\item Use the stable reduction argument to show that there is a {\em
modification} 
$Y' \to Y$, and a modification $Q' \to Q \times_YY'$ such that $Q'\to Y$ is a
family of nodal curves of genus 0, and the sections lift to {\em disjoint}
sections $s_i':Y'\to Q'$. We  replace $Y$ by $Y'$, $\Delta$ by its inverse
image, etc.
\item  Use induction on the dimension to replace $Y$ by a nonsingular variety
such that $\Delta$ becomes a strict divisor of normal crossings.
\item Use either toroidal geometry, or Section \ref{resolve-nodal}, to
replace $Q$ by a nonsingular variety, such that the inverse image of $D$ is a
strict normal crossings divisor.
\item Let $\tilde{X}$ be the normalization of $Q$ in the function field of
$X$. Use Abhyankar's lemma to show that  $\tilde{X}$ has a toroidal structure,
such that the inverse image $\tilde{Z}$ of $Z$ is a toroidal divisor.
\item Conclude that there is a weak resolution of singularities $r:X'\to X$
such that $r^{-1}Z$ is a strict divisor of normal crossings.
\end{enumerate}

            \section{Weak resolution of singularities II}
\label{Sec:A-J}

The weak resolution argument according to Abramovich - de Jong starts very much
like de Jong's theorem: a projection $X \to P$ is produced, and a Galois
alteration $P_1 \to P$ over which one has stable reduction $X_1 \to P_1$,
equivariant under the Galois group $G$, is produced. Induction on the dimension
for $P$ allows one to assume that $X_1$ and $P_1$ are toroidal, and the Galois
action on $P_1$ is toroidal as well. The only point left is to make the group
action on $X_1$ toroidal, so that the quotient should be toroidal, and
therefore admit toroidal resolution. 

Let us  go through the steps. Let $X$ be a variety over an algebraically
closed field of characteristic 0, and let $Z \subset X$ be a Zariski-closed
subset. We want to find a nonsingular, quasi-projective variety $X'$ and a
modification $f: X' \to X$ such that $f^{-1}Z$ is a divisor with simple normal
crossings. 
\subsection{Reduction steps}

\begin{exercise} 
Show that it is enough to prove the result when $X$ is projective and normal,
and $Z$  a Cartier divisor. 
\end{exercise}

\begin{exercise} 
Reduce to the case where there is a projection $X \to P$, such that the generic
fiber is a smooth curve.
\end{exercise}

\begin{exercise}\hardexercise  
Using the trick of enlarging $Z$ and stable reduction, show that there is a
diagram as follows: 
$$\begin{array}{ccc} X_1  &\to& X \\
                     \dar & & \dar \\
		     P_1  &\to& P 
\end{array}
$$
such that $P_1 \to P$ is an alteration, $X_1 \to \tilde{X}_{P_1}$ is
birational, and $X_1  \to P_1 $ has section $s_i: P_1 \to X_1$ making it a
family of stable pointed curves, and the image of these sections in $X$
contains $Z$.
\end{exercise}

\begin{exercise}
Show that you can make $P_1 \to P$ a Galois alteration. Call the Galois group
$G$. Show, using the uniqueness of stable
reduction\ref{Lem:unique-stable-reduction}, that the action of $G$ on 
$P_1$ lifts to an action on $X_1$, which permutes the sections $s_i$. 
\end{exercise}
   
You can  replace $X$ by $X_1/G$ and $P$ by $P_1/G$

\begin{exercise}
Use induction on the dimension to reduce to the case where:
\begin{enumerate} 
\item $P$ is nonsingular, with a normal crossings divisor $\Delta$;
\item The branch locus of $P_1 \to P$ is contained in $\delta$;
\item The locus where $X_1 \to P$ is not smooth is contained in $\delta$.
\end{enumerate}
\end{exercise}

\begin{exercise}
Show that in this case $P_1 \to P$ is a toroidal morphism, $G$ acts
toroidally on $X$, 
and $X_1 \to P_1$ is a toroidal morphism as well.
\end{exercise}

\n
The only point left is to make the action of $G$ on $X_1$ toroidal - if it
were, then $X$ would be toroidal and we could easily resolve its
singularities. 

Looking locally, the question boils down to the following
situation:

Let $T_0\subset X_0$  be an affine torus embedding, $X_0 = \Spec\
R$. 
Let $G\subset T_0$ be a finite subgroup of $T_0$, let $p_0\in X_0$ be
a fixed point of the action of $G$, and let $\psi_u$ be a character
of $G$. Consider the torus embedding of $T=T_0\times \Spec\ k[u,u^{-1}]$
into $X= X_0\times \Spec\ k[u]$, and let $G$ act
on $u$ via the character $\psi_u$.
 Write $p=(p_0,0)\in X$ and
write $D=(X_0\setmin T_0) \times \Spec\ k[u]$. We wish to find a canonical
blowup  $ X_1 \to X$, such that if $U\subset  X_1$ is the inverse image
of $T_0$, then it is a toroidal embedding, and the group $G$ acts toroidally.

\subsection{The ideal} Let $M\subset R[u]$ be the set of monomials. For each
$t\in M$ let $\chi_t$ be the associated character of $T$, and let
let $\psi_t:G\to k^*$ be the restriction of $\chi$ to $G$. Define  $M_u=\{t |
\psi_t = \psi_u\}$, the set of monomials on which  $G$ acts as it acts on
$u$. Define $I_G =  \left< M_u\right >$, the ideal generated by
$M_u$. 
\begin{exercise}[canonicity]  Show that if $X_0',\ T_0',\ G',\ p_0'$ and
$\psi'_u$ is a 
second set of such data, 
and if we have an isomorphism of formal completions
$$ \varphi : \hat{X}_p\stackrel{\sim}{\to} \hat{X'}_{p'}, $$
which induces isomorphisms $G\cong G'$ and 
$\hat{D}_{p}\cong \hat{D'}_{p'}$, then
$\varphi$ pulls back $I_G$ to the ideal $I_{G'}$.
\end{exercise}

\begin{exercise}[gluing property] If $q_0$ is any point of $X_0$ and if
$G_q\subset G$ is the stabilizer of $q=(q_0,0)$ in $G$, then the stalk of $I_G$
at $q$ is the same as the stalk of $I_{G_q}$ at $q$.
\end{exercise}

\begin{exercise} Show that $I_G$ is generated by $u$ and a finite number
of monomials $t_1,\ldots,t_m$ in $M_u\cap R$. 
\end{exercise}

\begin{exercise} Let $X' = B_{I_G}(X)$ be the blowup. Let $X'_u$ be the
chart with 
coordinates $u, t_j/u$. Show that the action of $G$ on $X'_u$ is toroidal.
\end{exercise}

\begin{exercise} Let $X'_i$ be the chart on $X'$ with coordinates $t_i, v=
u/t_i, 
s_j=t_j/t_i$. Show that $G$ acts trivially on $v$, and that $X'_i = \Spec\
R'_i[v]$ where $R'_i$ is generated over $R$ by $s_j$. 
\end{exercise}

\begin{exercise}  Let $X_1$ be the normalization of $X'$. Show that if
$U\subset  X_1$ is 
the inverse image of $T_0$, then it is a toroidal embedding, and the group $G$
acts toroidally. 
\end{exercise}

            \section{Intersection multiplicities}
%

Intersection theory has a long history, and certainly we are not going
to say much about it here. One aspect is, that it is not so easy to have a
good definition for intersection multiplicities.

\noindent{\bf Remark, exercise}: let $C \subset \PP^2_k$ be a
plane algebraic curve, $P \in C$ a closed point at which $C$ is
regular, and $D = {\cal Z}(F) \subset \PP^2$ a plane curve given by a
homogeneous polynomial $F$; suppose $F$ is not identical zero on a
neighborhood of $P$ in $C$ (i.e. no component of $D$ contains the
component of $C$ containing $P$). Show that the following two
definitions of {\it the intersection multiplicity $i(C,D;P)$ of $C$
and $D$ at} $P$ are equivalent: 
\begin{itemize}
\item the dimension of the $k$-vector space 
$${\cal O}_{C,P} \ \ \otimes_{{\cal O}_{\PP^2,P}} \ \ {\cal O}_{D,P},$$
\item the value of the valuation $v = v_{C,P}$ defined by the discrete
valuation ring  ${\cal O}_{C,P}$ computed on the function on $C$ given
by $F$, 
\end{itemize}
see [HAG], Exercise (5.4) on page 36, and Remark (7.8.1) on page 54.

 Consider two varieties $V, W \subset \PP^n$ which have an isolated
point of intersection at $P \in V \cap W$. 
One could try to define the intersection of $V$ and $W$ at $P$ as the
length of  
$${\cal O}_{V,P} \ \ \otimes_{{\cal O}_{\PP^n,P}} \ \ {\cal O}_{W,P}.$$
Analogous situations of intersections of arbitrary schemes in some
regular ambient scheme can be considered. 

\begin{exercise} (See Gr\"obner \cite{Groebner}, 144.10/11, also see
\cite{Serre}, \cite{Roberts}, \cite{Berthelot}, see [HAG],
I.7):
{\bf a)}\ \  Let $C \subset \PP^3$ be the space curve with parameterization
$$(x_1 \ : \ x_2 \  : \  x_3\  :\  x_4) = (t^4 \ :\  t^3{\cdot}s \ :\
t{\cdot}s^3\  : \ s^4)$$ 
(we work over some field $K$). Show that the prime ideal given by this
curve equals 
$$j := (T_1^2T_3-T_2^3 \ \ , \ \  T_1T_4-T_2T_3 \ \ , \ \
T_1T_3^2-T_2^2T_4 \ \ ,\ \  T_2T_4^2-T_3^3) \subset
K[T_1,T_2,T_3,T_4].$$ 
{\bf b)}\ \  Consider $C$ as a curve embedded in $\PP^4$: choose the
hyper plane  
$\PP^3  \cong {\cal Z}(T_0) = H$, and we get $C \subset H \subset
\PP^4$. Let $P:=(x_0=1:0:0:0:0) \in \PP^4$. Define $V \subset \PP^4$
as the cone with vertex $P$ over the curve $C \subset \PP^4$, i.e. $V$
is defined by the ideal 
$$J :=  K[T_0,T_1,T_2,T_3,T_4] \cdot j, \quad V = {\cal Z}(J).$$
Note that the dimension of $V$ equals two, that the degree of
$V\subset \PP^4$ equals four.\\ 
{\bf c)}\ \   Let $W$ be the 2-plane given by
$$I:=(T_1,T_4) \subset  K[T_0,T_1,T_2,T_3,T_4], \quad  W := {\cal Z}(I).$$
Note that $P \in W$. Remark that (set-theoretically):
$$W \cap V = \{P\}$$
(use the geometric situation, or give an algebraic computation).\\
We like to have a B\'ezout type of theorem for the situation $W \cap V
\subset \PP^4$, however:\\ 
{\bf d)}\ \  Define
$$M := {\cal O}_{W,P}, \quad  A:={\cal O}_{{\PP^4},P}
, \quad     N:= {\cal O}_{V,P},$$
and compute
$$\dim_K \left(M \otimes_A N\right)$$
(surprise: this is not equal to four).\\
{\bf e)}\ \   Compute
$$\dim_K \left(\Tor^A_i(M,N)\right), \ \  \forall i$$
(either using, or reproving $\chi_A(M,N)=4$, for notation see below).
\end{exercise}

\n
Hence we see that just the length of the appropriate tensor product
does not define necessarily the correct concept. Serre proposed in
1957/58 to define the intersection multiplicity as the alternating sum
of the lengths of the $\Tor_i$ (note that $\Tor_0 = \otimes$), i.e. by

$$\chi_A(M,N) := \sum\limits_{i \geq 0} \ \ \length_A \Tor^A_i(M,N)$$ 
(we follow notation of \cite{Serre}, also see \cite{Berthelot},
6.1, see \cite{Roberts}), here $A$ is a regular local ring, and $M$
and $N$ are $A$-modules such that $M \otimes_A N$ has finite
length. In equal characteristic this is the right geometric concept
(i.e. satisfies B\'ezout's theorem, coincides with previously defined
intersection multiplicities etc.).

The following theorem was conjectured by Serre, proved by Gabber
(using de Jong's alteration result), 
and written up by Berthelot (in
\cite{Berthelot}, 6.1):

\begin{theorem}   Let the characteristic of $A$ be equal to
zero. Suppose $p \in \mathfrak{m}^2$, hence  its residue field
$A/\mathfrak{m}$ 
has characteristic $p >0$.  Then: 
$$ \chi_A(M,N) \geq 0.$$ 
\end{theorem}

\newpage\part{Moduli of curves}
            \section{Introduction to moduli of curves}
\mylabel{Sec:Intro-to-moduli}

It is an important feature of algebraic geometry, that the set of all objects
(e.g. smooth projective curves)  of the same
a fixed geometric nature (e.g. genus) often has the structure of an algebraic
variety itself.   Such a space is a ``moduli space'', which gives a good
algebraic meaning to the problem of ``classification''. It is fair to
say that this ``self referential'' 
nature of algebraic geometry is one of the main reasons for the depth
of the subject - it is impossible to overestimate its importance.

The first instances of this phenomenon to be discovered were those of
{\em embedded} variety: the projective space as a parameter space for
lines in a vector space; Grassmannians parametrizing vector subspaces
of arbitrary dimension; the projective space (of dimension $(d^2 +
3d)/2$) parametrizing all plane curves  of degree $d$, and so
on. The case of abstract varieties, such as smooth curves, had to
await for some technical advances, although already Riemann knew that
algebraic curves of genus $g$ ``vary in $3g-3$
parameters''; see \cite{R}, page 124:  ``Die $3p-3$ \"ubrigen Verzweigungswerthe in jenen Systemen gleichverzweigter $\mu$-werthiger Functionen k\"onnen daher beliebige Werthe annehmen; und es h\"angt also eine Klasse von Systemen gleiverzweigter $(2p+1)$-fach zusammenhangender Functionen und die zu ihr geh\"orende Klasse algebraischer Gleichungen von $3p-3$ stetig ver\"anderlichen Gr\"ossen ab, welche die Moduln dieser Klasse genannt werden sollen."
 although already Riemann knew that
algebraic curves of genus $g$ ``vary in $3g-3$
parameters''; see \cite{R}, page 124:  ``Die $3p-3$ \"ubrigen Verzweigungswerthe in jenen Systemen gleichverzweigter $\mu$-werthiger Functionen k\"onnen daher beliebige Werthe annehmen; und es h\"angt also eine Klasse von Systemen gleiverzweigter $(2p+1)$-fach zusammenhangender Functionen und die zu ihr geh\"orende Klasse algebraischer Gleichungen von $3p-3$ stetig ver\"anderlichen Gr\"ossen ab, welche die Moduln dieser Klasse genannt werden sollen."

Historically moduli spaces of curves, or of curves with points on
them, were constructed with more or less ad hoc
methods. Moduli spaces for 3 or 4 points on rational curves have been known for
ages, using the so called ``cross ratio'' (see exercise
\ref{rational-case} below). For genus 1, the modular function
$j$ was used (see exercise
\ref{j-function}). The case of genus 2 was already quite
difficult to achieve by 
algebraic methods \cite{Ig}. For years, moduli spaces for higher genus
were  only known   
to exist using Teichm\"uller theory. 

One problem which took years to solve was, that no good understanding
of ``what moduli spaces really are'' was available. Then
Grothendieck introduced the notion of ``representable functor", describing the
best possible meaning for moduli spaces. This had a great success with
the development of Hilbert schemes. For a while one 
hoped that nature would be as ideal as expected (see Grothendieck hopeful
Conjecture 8.1 in \cite{FGA} 212-18, and its retraction in the Additif
of \cite{FGA} 221-28). But it was soon seen that in general 
moduli functors are not representable, or as we say now, some  moduli functors
do not give rise to 
``fine moduli spaces'' due to existence of automorphisms. Finally, Mumford
pinned down 
 the compromise  notion of a ``coarse moduli scheme'',  which enables us to
have a good insight in various aspects of moduli theory. This is what we shall
try to describe here. It should be said that since then, other good
approaches were developed, by way of ``enlarging the category of
schemes'' to include some ``moduli objects'', called stacks. For
details see \cite{D.M}, Section 4. We will not pursue this direction here.

In this section we gather  some basic definitions on functors of moduli for
curves. In the next sections we discuss existence theorems for moduli
spaces of curves, and for complete moduli 
spaces with extra structure carrying a ``tautological family". 
 Section \ref{Sec:Moduli-other} is devoted to
some further questions, examples and facts, not needed for the methods of
alterations, but in order to give a more complete picture of this topic. 

\subsection{The functor of points and representability}
To any scheme $M$ one naturally associates a contravariant functor 
$$\cF_M: \{\mbox{Schemes}\} \to  \{\mbox{Sets}\} $$
via 
$$ X \mapsto \Mor(X,M).$$ This is known as {\em the functor of points of $M$},
see \cite{RedBook}, II \S 6.

We say that a  contravariant functor $\cF:\{\mbox{Schemes}\} \to
\{\mbox{Sets}\}$  is {\em representable} by a scheme $M$, if it is isomorphic
to 
$\cF_M$, 
 i.e. there is a functorial isomorphism 
$$\xi: \cF(-) \stackrel{\sim}{\longrightarrow} \Mor_S(-,M).$$

\begin{remark}
Strictly speaking, it is  the {\em pair} $(M,\xi)$, consisting of the
object $M$ {\em and}
the isomorphism $\xi$, which represent $\cF$. But it has become customary to
say ``$M$ represents $\cF$'', suppressing $\xi$. 
\end{remark}

\n
Already in the early Bourbaki literature one finds this notion in  the disguise
of a ``universal property''. The question of representability of
functors can also be posed in 
categories other than the category of schemes.

\begin{definition}\label{taut} Suppose $T$ is a scheme, and let $f: T \to \oM$ be a morphism. A family $\cC \to T$ of stable curves is called a {\it tautological family} if it defines $f$. In particular this implies that for every geometric point $t \in T$ the fiber $\cC_t$ is a stable curve which defines the image uder $f$, i.e.: $[\cC_t] = f(t)$.
\end{definition}

\n
Note: If $\cC \to T$ represents a functor of stable curves (plus some extra structure), i.e. if this is a fine moduli scheme, then this universal family is tautological.\\
Note: There exist cases (and we shall give examples), where a moduli functor is not representable, where there is no (unique) universal family, but where a tautological family does exist. In such cases the use of the word ``tautological", and the distinction between ``universal"  and  ``tautological" is necessary, and it pins down the differences.

The terminology ``tautological" will also be used in cases such as pointed curves, curves with a level structure, and so on.

\begin{exercise} \begin{enumerate} 
\item Fix an integer $N$, and let $V$ be a vector space of dimension
$N+1$ over $\CC$. Consider the functor $\cF_\CC$ that associates to a
scheme $T$ over $\CC$, the set $\{\cL \subset T\times V\}$ of all line
sub-bundles of the trivial vector bundle
$T\times V$. Show that $\cF_\CC$ is represented by $\PP^N_\CC$.
 \item  Let $\cF$ be the functor that associates to any scheme $T$
(over $\ZZ$) the set $\{\cL \subset \cO_T^{N+1} \}$ of all locally free
subsheaves of rank 1 of the trivial sheaf $\cO_T^{N+1}$ having locally
free quotient. Show that
$\cF$ is represented by the projective scheme $\PP^N$.
\item In general, show that the Grassmannian scheme ${\mathbf{
Grass}}(n,r)$ represents the functor of locally free subsheaves of rank
$r$ of the trivial free sheaf $\cO_T^n$ of rank $n$ having locally
free quotients.     
\end{enumerate}
\end{exercise}

\begin{exercise} Fix  integers $N$ and $d$, and let $\cG$ be the
functor that   associates to a
scheme $T$ the set $\{\cX \subset \PP^N_T\}$ of all flat families of
hypersurfaces of degree $d$ in projective $N$-space over $T$. Show
that $\cG$ is represented (over $\ZZ$) by a projective space $\PP^{M-1}$,
where $M = {\binom{N+d}{d}}$ is the dimension of the space of
homogeneous polynomials of degree $d$ in $N+1$ variables.
\end{exercise}

\begin{exercise} Show that the Hilbert scheme $\cH_{\PP^N,P(T)}$
represents the ``Hil\-bert'' functor, that associates
to a scheme 
$T$, the set of all 
subschemes $X\subset \PP^N_T$, which are flat over $T$ and such that
the geometric fibers have Hilbert polynomial $P(T)$.
\end{exercise}

\subsection{Moduli functors and fine moduli schemes}
Suppose a contravariant functor $\cF$ has the nature of a {\em moduli functor},
namely,  it 
assigns to a scheme $S$ the set $\{C \to S\}/\cong$ of isomorphism
classes of certain families of 
objects over 
$S$. As a guiding example, let us fix an integer $g$, with $g \in \ZZ_{\geq
0}$, and define the moduli functor for smooth curves: $${\cal M}_g(S)=
\{\mbox{isom. classes of families of curves 
of genus $g$ over   } S\}.$$ 
A morphism $T \to S$ defines (by pulling back families) a map of sets
in the opposite direction: $\cM(T) \leftarrow
\cM(S)$, and we have obtained a contravariant functor.

Assume the functor $\cF$ were represented by a scheme $M$.
Then we would call $M$ a {\em fine moduli scheme} for this functor $\cF$, and
the object $${\cal C} \rightarrow M \quad\mbox{corresponding to the
identity}\quad id \in \cF(M)  = \Mor_S(M,M)$$ 
would be called a {\em universal family}.

\begin{remark} Note that in the exercise above on the Hilbert scheme,
we can view it as a fine moduli scheme, if we agree that ``families up to
isomorphism'' means  ``up to isomorphisms as subfamilies of the fixed
$\PP^N_T$, 
namely up to equality.
\end{remark}

\n
It is a  fact of life that {\em for every $g \geq 0$ the functor
$\cM$ is not 
representable.} We will explain later why this is true in general, but for the
moment let us consider the easiest case:

\begin{exercise} Let us say that $C$ is a ``curve of genus 0'', if it is an
algebraic 
curve defined over a field $K$, and over some extension of $K \subset L$ it is
isomorphic with $\PP^1_L$. In other words: $C$ is geometrically irreducible,
reduced, it is complete and of genus equal to zero.
\end{exercise}

\begin{enumerate}
\item Let $K$ be a field. Show there exist  an extension $K \subset K'$, and
two  curves of genus $0$ over $K'$ which are not isomorphic. 
\item  For every algebraically closed field $k$, the set ${\cal M}_0(k)$
consists of one element, ${\cal M}_0(k)= \{\PP^1_k\}$.
\item Show that the moduli functor ${\cal M}_0$ is not representable. 
\end{enumerate}

\subsection{Historical interlude}  
The first case of a highly nontrivial algebraic construction of a
 moduli space of curves in all charactersitics, appeared in Igusa's work
 \cite{Ig}. This is a construction of
 a  ``moduli scheme for non-singular curves of
genus two in all characteristics'', which would now be denoted by $M_2 \to
\Spec(\ZZ)$. This happened almost concurrently with
 Grothendieck's 
 study of representability of functors. But notice that, when Samuel
 discussed these beautiful results by Igusa in 
S\'eminarie Bourbaki (see \cite{Sam}), his very first comment was: 

\begin{quote}
{\em ``Signalons ausit\^ot que le travail d' {\sc Igusa} ne r\'esoud pas
pour les courbes de genre 2, le ``probl\`eme des modules'' tel 
qu'il a \'et\'e pos\'e par {\sc Grothendieck}  \`a diverses reprises 
dans  ce S\'eminaire."}
\end{quote}

\n
It really seemed that Nature was working against algebraic geometers, refusing
to provide us with these fine moduli schemes...

The truth is, Nature does provide us with a replacement. Indeed,  not much
later, Mumford (see \cite{GIT}, 5.2) discovered how to follow nature's
dictations 
and come to  a good working definition, requiring that the scheme
gives geometrically what you want, and does it in the best possible way.

\subsection{Coarse moduli schemes}
Here is the definition:
\begin{definition}\mylabel{Def:Coarse-moduli}
 A scheme $M$ and a  morphism  of functors
$$\varphi: F \rightarrow \Mor_S(-,M)$$
is called a {\em coarse moduli scheme} for $F$ if:
\begin{enumerate}
\item  for every algebraically closed field $k$ the map
$$\varphi(k): F(\Spec(k)) \rightarrow \Mor_S(\Spec(k),M)=M(k)$$
is bijective, and 
\item  for any scheme $N$ and any morphism  $\psi: F(-) \rightarrow
\Mor_S(-,N)$ there is a unique $\chi: N \rightarrow M$ factoring $\psi.$  
\end{enumerate}
\end{definition}

\n
And here is its first triumphant success:

\begin{theorem}[Mumford]  Suppose $g\geq 2$. The functor ${\cal M}_g$ of
smooth cur\-ves of 
genus $g$ admits a quasi-projective coarse moduli scheme.
\end{theorem}

\n
See \cite{GIT}, Th. 5.11 and Section 7.4,
 or \cite{D.M},
Coroll. 7.14.  We will denote the coarse moduli scheme of ${\cal M}_g$ by   $M_g
\rightarrow  \Spec \ZZ$.

\n 
We note some properties:
\begin{itemize}
\item As we mentioned before, for every  $g\geq 2$ the functor $\cM$
is not representable: there does not exist a universal family of
curves over $\M$ which can give an isomorphism between $\cM$ and
$\M$.  
\item  For every  $g\geq 2$ and for any field $K$, the variety $(\M)_K
= \M\times_{\Spec \ZZ}\Spec K$ is not complete. A fortiori, the
morphism  $M_g  
\rightarrow  \Spec (\ZZ)$ is not proper.
\item At least for the sake of de Jong's theorem, we need a moduli
space of curves with points on them. 
\end{itemize}

\n
The first problem is solved by introducing a finite covering $M \to
\M$ admitting a {\em tautological family}, namely a family
realizing the morphism $M \to
\M$ as its moduli morphism. The nicest way of doing
this is by introducing a new moduli functor, of smooth curves ``enriched''
with a
finite amount of ``extra structure'', which does admit a fine moduli
scheme. See Section \ref{Sec:Tautological-pointed}.

In order to ``compactify" these spaces, the notion of {\em stable
curves} was invented. Historically, the influential paper \cite{D.M}
by Deligne and Mumford seems to be one of the first printed versions
in which the concept of stability, especially in the case of algebraic
curves is explained and used. In \cite{Mumford-WH} 
we see that
already in
1964 Mumford was trying to find the appropriate notions assuring good
compactifications. In \cite{GIT}, page 228, Mumford attributes the
notion of a stable curve to unpublished joint work with Alan Mayer. 

As it turns out, the third problem was solved almost concurrently with
the second. 

First, the moduli space of smooth pointed curves:

\begin{theorem}  Let $g \in  \ZZ_{\geq 0}$,
and $n \in \ZZ_{\geq 0}$ such that $2g-2+n > 0$. Consider the functor
$\cnM$ of isomorphism classes of families of stable smooth $n$-pointed curves
of  genus $g$.  
 This
functor admits a quasi-projective coarse moduli scheme. 
\end{theorem}

\n
We will denote this moduli space by $\nM 
\rightarrow  \Spec(\ZZ)$.

\begin{remark} \begin{enumerate}
\item  Note that this includes the previous theorem. 

\item For $g=0$ and $n=3$ this space is proper over
$\Spec(\ZZ)$. However in all other cases in the theorem  $\nM$  is not
proper. In many cases it will not represent the functor (see Section
\ref{Sec:Moduli-other} for a further discussion), in other words, in
general this is not a   fine moduli  scheme. 

\item It is important to note that these spaces exist over
$\Spec(\ZZ)$, which is useful for
arithmetical applications.

\item The litterature poses difficulties in choosing notations. In
\cite{GIT} 
the subscript $n$ denotes a level structure, but in \cite{Knudsen} it
indicates the number of marked points. We have chosen to indicate the
markings as lower index, using $n$, and the level structure as
upper index, using $(m)$. 
\end{enumerate}
\end{remark}

\n
Finally the moduli space of stable pointed curves:

\begin{theorem}[Knudsen and Mumford]
 Let $g \in  \ZZ_{\geq 0}$, and $n \in \ZZ_{\geq 0}$ such that $2g-2+n> 0$. 
Consider the functor $\conM$ of isomorphism classes of families
 of stable $n$-pointed curves of 
 genus $g$.  This functor admits a {\em projective} coarse moduli scheme.
\end{theorem}
\n
See \cite{Knudsen}, part II, Theorem 2.7 and part III, Theorem 6.1, or
\cite{Gies}, Theorem 2.0.2. 
 We
will denote this moduli scheme by $\onM 
\rightarrow  \Spec(\ZZ)$.

\n
The following exercise should give you an idea why the moduli space
$\overline{{\cal M}_{0,n}}$ is complete. This is discussed in further
detail in the next section.

\begin{exercise} Let $K$ be a field, and let $R \subset K$ be a
discrete valuation ring having $K$ 
as field of fractions. Consider the projective line
$\PP^1$ over $K$ and  suppose $n \geq 3$, let $P_1, \cdots , P_n \in
\PP^1(K)$ are distinct points. Write $P = \{P_1, \cdots ,
P_n\}$. Construct a stable $n$-pointed curve $(\cC,\cP)
\to \Spec(R)$ extending $(C,P)$. (You will need to blow up closed
points over the special fiber where the Zariski closures of $P_i$
meet. Then you may need to blow down some components! See \cite{Knudsen}.)
\end{exercise}

\begin{exercise} 
\begin{enumerate}\mylabel{rational-case}
\item Let $K$ be a field. Given three distinct finite points $P_1,
P_2$ and  $P_3$ on
$\PP^1_k$ consider the cross ratio $$\lambda(P_1,P_2,P_3,z) = {\frac{(z-P_1)
(P_2-P_3)}{(z-P_3) (P_2-P_1)}}.$$
Show that, as a function of $z$, the cross ratio is an automorphism of
$\PP^1$ carrying $P_1,P_2,P_3$ to $0, 1,$ and $\infty$,
respectively. Show that this automorphism is the unique one with this
property. Check that this definition can be extended to the case where
one of the points is $\infty$. 
\item Using the cross ratio $\lambda$ defined above, describe
$M_{0,3}$. 
\item Show that $M_{0,3}$ is a fine moduli scheme by exhibiting a
universal family over it!
\item  Show that 
$M_{0,3} = \overline{M_{0,3}}$. 
\item Use the cross ratio to give an explicit description of
$M_{0,4}$. Show that it is a fine moduli scheme by explicitly
constructing a universal family. 
\item Use the above (possibly together with the previous exercise) to
describe $M_{0,4} \subset \overline{M_{0,4}}$.
\item  Show that $\overline{M_{0,4}}$ is a  fine
moduli scheme, and give  explicit descriptions of the universal
family.
\item Show that the universal family over $\overline{M_{0,4}}$ is
canonically isomorphic to $\overline{M_{0,5}}$.
\end{enumerate}
\end{exercise}

\begin{exercise} Give an alternative description of
$\overline{M_{0,4}}$ as follows: consider the projective space $\PP$ of
dimension 5 parametrizing conics in $\PP^2$. Choose four points in
general position in $\PP^2$ (for instance $(1:0:0), (0:1:0), (0:0:1),
(1:1:1)$ will do). Let $M\subset \PP$ be the subscheme parametrizing
conics passing through these four points. Show that $M =
\overline{M_{0,4}}$ and the universal family of conics is a universal
family  for $\overline{M_{0,4}}$.
\end{exercise}

\begin{exercise} \hfill
\begin{enumerate}
\item Show that $M_{0,n}$ exists and is a fine moduli scheme (you may
exhibit it as an open subscheme of $(\PP^1)^n$).
\item Show that, assuming $\overline{M_{0,n}}$ is a fine moduli
scheme, then there is a canonical morphism $\overline{M_{0,n+1}}\to
\overline{M_{0,n}}$ which exhibits $\overline{M_{0,n+1}}$ as the
universal family over $\overline{M_{0,n}}$.
\item\hardexercise\
Show that for every $n \geq 3$, the scheme $\overline{M_{0,n}}$ is a
fine moduli 
scheme. (You may want to use Knudsen's stabilization technique.)
\end{enumerate}
\end{exercise}

\begin{remark} For every $n \geq 3$, let $(C,P)$ be a
stable $n$-pointed rational curve. Then $\Aut(C,P) = \{id\}$. You do
not need to
know this in the previous exercise, but it ``explains" why the
result should be true. 
\end{remark}

\begin{exercise}
\mylabel{j-function}
 Let $k$ be a field of characteristic $\neq 2$ and let
$(E,O)$ be an elliptic curve, namely a projective, smooth and
connected curve $E$ of genus 1 with a $k$-rational point $O$ on it.
\begin{enumerate}
\item Considering the linear series of $\cO_E(2O)$, show that $E$ can
be exhibited as a branch covering of $\PP^1$ of degree 2.
\item Show that the branch divisor $B$ on $\PP^1$ is reduced and has
degree 4.  
\item If $k$ is algebraically closed, show that $E$ is determined up
to isomorphism by the divisor $B$. 
\item Conclude that $M_{1,1}$ is isomorphic to the quotient
of $M_{0,4}$ by the action of the symmetric group $S_4$, permuting the
four points. 
\item Assume further that $\chara k \neq 3$, so that every elliptic
curve can 
be written in affine coordinates as $y^2 = x^3 + ax + b$. Show that 
$$j(E) = 1728{\frac{4a^3}{4a^3 + 27 b^2}}$$ is an invariant
characterizing the $\bar k$-isomorphism class of $E$, exhibiting
$M_{1,1} = \AAA^1$. 
\end{enumerate}
\end{exercise}

\section{Stable reduction and completeness of moduli spaces}
\subsection{General theory}
In order to understand the reason why $\onM$ is projective, let us
recall the following:

\begin{theorem}[The valuative criterion for properness] A morphism
$$f: X \to 
Y$$ of finite type is proper, if and only if the following holds: 

Let $R$ be a  discrete valuation ring,
and let  $S 
:= \Spec(R)$ be the corresponding ``germ of a  non-singular curve'',  with
 generic point $\eta$.  Let $\varphi: S
\to Y$ and let $\psi_\eta: \eta  \to X$ be a lifting:
$$\begin{array}{ccc} \eta & \stackrel{\psi_\eta}{\to} & X \\ 
			\dar & 			& \dar \\
		      S   & \stackrel{\varphi}{\to} & Y.
\end{array}$$
 Then there
is an extension $\psi: \Gamma \to X$, lifting $\varphi$:
$$\begin{array}{ccc} \eta & \stackrel{\psi_\eta}{\to} & X \\ 
		    \dar & \stackrel{\psi}{\nearrow}& \dar \\
		      S   & \stackrel{\varphi}{\to} & Y.
\end{array}$$
\end{theorem}
\n 
See \cite{HAG},
II, Theorem (4.7) for a precise formulation.

\n
Let us translate this to our moduli scheme. Keeping in mind the
relationship between the functor $\conM$ and the space $\onM$, one
might hope that every family of stable pointed curves over $\eta$ as
in the theorem above
might extend to $R$. This is not the case, as we shall see
later. However, a weaker result, sometimes called ``the {\em weak} valuative
criterion for properness'', does hold for the functor $\conM$, and it
does imply the valuative criterion for $\onM$. The first case to
consider is when the generic fiber is {\em smooth} and $n=0$.
This is
the content of the following result, the Stable Reduction Theorem for a one
parameter family of curves:

\begin{theorem} \mylabel{SR}
Let $S = \Spec(R)$ be the spectrum of a discrete valuation ring $R$,
$\eta\in S$ the generic 
point, corresponding with the field of fractions $K$ of $R$. Let
$C_\eta\to \eta$ be a smooth stable curve of genus $g>1$. There exists a 
finite extension of discrete valuation rings $R \hookrightarrow R_1$,  with
$S_1 = \Spec R_1$ and
generic point $\eta_1$, and an extension 
$$\begin{array}{ccl} C_{\eta_1}& \hookrightarrow  & C_1 \\
 \downarrow & &  \downarrow \\
			\{\eta_1\}& \hookrightarrow  & S_1 ,
\end{array}$$ 
such that $C_1 \to S_1$ is a family of stable curves.
\end{theorem}

\n  
Proofs of this theorem, using different methods, may be found in
  various references. One proof which works  in pure characteristic 0 is
  relatively simple. As the reader will notice,
  none of the general proofs  is easy or elementary. 

 Most  proofs of this theorem use resolution of singularities of
 2-dimensional schemes (Abhyankar).

\begin{exercise} Suppose $R$ is of pure characteristic 0. Let $s\in
 \Spec R$ be
 the closed point.
\begin{enumerate}
\item Show that there exists an extension $\pi:C \to S$ such that
 $\pi$  is proper and flat, $C$ is nonsingular, and $C_s\subset C$ is a
 normal crossings divisor.
\item Let $ x\in C_s$ be a singular point. After passing to the
algebraic closure of the field of constant, let  $\bar{x} \in
\hat{C}_{\bar{k}}$ be the completion. Show that one can find local
parameters $u,v$ at 
$\bar{x}$ and $t$ at $\bar{s} \in \hat{S}_{\bar{k}}$, and positive
integers $k_x, 
l_x$, such that $t = u^{k_x} v^{l_x}$.
\item Let $S_1 \to S$ be a finite cover obtained by extracting the
$n$-th root of a uniformizer, where $n$ is divisible by all the
non-zero $k_x, l_x$ given above. Let $C'_1$ be the normalization of $C\times_S
S_1$. Show that the special fiber is reduced and nodal.   
\item Show that the minimal model $C_1$ of $C'_1 \to S_1$ is
stable. 
\end{enumerate}
\end{exercise}

\n
We list some approaches for positive and mixed characteristic:

\n
{\bf Artin-Winters.} This proof can be found in \cite{A.W}. A
precise and nice description and analysis of the proof is given by
Raynaud, see  \cite{Raynaud}.

In this proof one attaches an numerical invariant to a given genus,
and one proves that by choosing a prime number $q$ larger than this
invariant, and not equal to the residue characteristic, and by
extending the field of definition of a curve of that genus such that
all $q$-torsion point on the jacobian are rational over the extension,
then one acquires stable reduction. The proof consists of a careful
numerical analysis of the possible intersection matrices of
components of degenerating curves. The proof does not rely on a lot of
theory, but is quite subtle.

\n
{\bf Grothendieck, Deligne-Mumford.} This proof can be found in
\cite{D.M}, Theorem (2.4) and Corollary (2.7). 

In this proof one shows that a curve has stable reduction if and only if its
jacobian has stable reduction. Then one shows following Grothendieck
that eigenvalues of algebraic $\ell$-adic monodromy are roots of unity
(see \cite{S.T}, Appendix).
Moreover, again following Grothendieck
one shows that  these eigenvalues are all equal to one iff the abelian
variety in question has stable reduction. The advantage of this proof
is that it has a more conceptual basis. The big disadvantage is that
it relies on the theory of N\'eron models, whose foundations are quite
difficult.

\n
{\bf Hilbert schemes and GIT - Gieseker.} See \cite{Gies}, Chapter 2,
Proposition(0.0.2). He says on the first page of the introduction:
``...we use results of Chapter 1 to give an indirect proof that the
$n$-canonical embedding of a stable curve is stable if $n \geq 10$,
and to construct the projective moduli space for stable curves. As
corollaries, we obtain proofs of the stable reduction theorem for
curves, and of the irreducibility for smooth curves."  The proof uses
 Geometric Invariant Theory to prove directly that $\oM$
exists and is projective, and then one can easily derive the
theorem. This proof does not use resolution of singularities for
surfaces in any explicit manner.

\begin{remark} This theorem is an instance of the {\em semistable reduction
problem.} In \cite{Berthelot}, 1.3,  the
definition of  semistable reduction, over a one-dimensional
base, and {\em arbitrary} fiber dimension, is recalled. As we have
seen above, it  
is true that 
if the relative dimension is one, stable reduction,
hence semistable reduction,  exists over a one-dimensional base. 
 For higher
relative dimension an analogous result holds in pure
characteristic zero - see \cite{kkms}.  The general case
is an important open problem, which seems difficult. 
\end{remark}

\n
Once Theorem \ref{SR} is known, it is easy to generalize
it. The
pointed case can be easily proven using Knudsen's stabilization technique:

\begin{exercise} Let $R$ be a discrete valuation ring, with field of
fractions $K$, suppose $(C,p_1,\ldots,p_n)$ is a {\em smooth, stable
$n$-pointed curve} of genus $g>1$  
defined over $K$. There 
exists a finite extension $R \subset R_1$ of valuation rings, with
$K_1$ the field of fractions of $R_1$,  such that $C_1 = C \otimes
K_1$ extends to a stable $n$-pointed curve $\cC_1 \to \Spec(R_1)$.
\end{exercise} 
  
\n
The case of genus zero follows from  Exercise
\ref{rational-case}. We
will discuss the case of genus 1 in Section \ref{elliptic-case} below.

 We can also consider the case when the generic fiber is not
necessarily smooth:

\begin{exercise} Let $R$ be a discrete valuation ring, with field of
fractions $K$, suppose $C$ is a stable curve defined over $K$. There
exists a finite extension $R \subset R_1$ of valuation rings, with
$K_1$ the field of fractions of $R_1$,  such that $C_1 = C \otimes
K_1$ extends to a stable curve $\cC_1 \to \Spec(R_1)$. [Below we
formulate a generalization to stable pointed curves of this.]
\end{exercise}

\n
We  give a full generalization of (\ref{SR}):
\begin{exercise}\mylabel{SRg}  
Let $S$ be the spectrum of a discrete valuation ring, $\eta\in S$ the generic
point. Let $(C,P) \to \{\eta\}$ be a stable $n$-pointed curve of
genus. I.e., $C$ is a complete, nodal curve  defined over a field
$K$, and  $P := \{P_1,
\cdots , P_n\}$ are distinct closed points $P_j \in C(K)$, with such
that $(C,P)$ is  stable $n$-pointed over 
$K$. 

Then there exists 
a 
finite extension of discrete valuation rings $S \hookrightarrow S_1$, with
generic point $\eta_1$, and an extension 
$$\begin{array}{ccc} (C_{\eta_1},P)& \hookrightarrow  & (\cC,\cP) \\
 \downarrow & &  \downarrow \\
			 \Spec(K_1)& \hookrightarrow  & S_1
\end{array}$$
such that $(\cC,\cP) \to S_1$ is a family of stable $n$-pointed
curves. 
\end{exercise} 
\n
This is the ``weak valuative criterion for properness'' of the functor
$\conM$.

\begin{remark} Here is a hint about a technical detail which can be
used in solving the previous exercises, {\it ``The normalization of a
stable $n$-pointed curve"}: 
 Suppose given a stable $n$-pointed curve $(C,P)$ over a field
$K$, with $P = \{P_1, \cdots , P_n\}$. There exists a finite extension
$K \subset L$, a finite disjoint union $(D,Q)$  of stable pointed
curves, and a morphism (``the normalization") $\varphi: (D,Q) \to
(C,P)_L = (C,P) \otimes_K L$ such that: $D = \coprod D^{(t)}$, let the
singular points of $C_L$ be: $R_j \in C(L)$, with $1 \leq j \leq d$,
moreover  $Q = \{Q_1, \cdots , Q_n\} \cup \{S_j, T_j \mid 1 \leq j
\leq d\}$, for every irreducible component of $C_L$ there is a unique
component of $D$ mapping birationally onto it, the morphism $\varphi$
is an isomorphism outside $\Sing(C_L)$, the markings $Q_i$ corresponds
with the markings $P_i$ of the pointed curve $(C,P)$, and the markings
$\{S_j, T_j\}$ are precisely the points mapping to $R_j$.

You need to show this choice can be made, and show it is unique in
case $K=k$ is an algebraically closed field.  
\end{remark}

\begin{corollary} Let $g \in  \ZZ_{\geq 0}$, and $n \in \ZZ_{\geq 0}$ such
that $2g-2+n > 0$. The coarse moduli scheme $\pi: \onM 
\rightarrow  \Spec(\ZZ)$ is proper over $\Spec(\ZZ)$.
\end{corollary}

\n
{\bf Proof:} We use the valuative criterion for properness setting
$X = \onM$ and $  Y = \Spec(\ZZ).$
Suppose $R$ is a discrete valuation ring, with field of fractions $K$,
and  suppose given  
$$\begin{array}{ccc} \Spec K & \stackrel{\psi_K}{\to} & X \\ 
			\dar & 			& \dar \\
		      \Spec R   & \stackrel{\varphi}{\to} & Y.
\end{array}$$

By the definition of a coarse moduli scheme, there is a finite
extension $K \subset K'$ such that  the point  $\psi_K(\Spec(K))  \in
X$ corresponds to a stable pointed curve $(C,P)$ over $K'$. By the
stable reduction theorem there is a finite extension $K' \subset K_1$
such that $(C,P) \times_{\Spec K'} \Spec K_1$ extends to a stable
pointed curve; this 
defines a morphism $\tau: \Spec(R_1) \to X$, ``extending'' $\varphi$ and
$\psi_K$. It factors over $\Spec(R)$, because $R = K \cap R_1$. This
shows that the condition for the valuative criterion holds in our
situation, hence  that  $\pi: \onM 
\rightarrow  \Spec(\ZZ)$  is proper. \qed

\subsection{Stable reduction for elliptic curves.}\mylabel{elliptic-case}

In these exercises we illustrate the concept of stable
reduction by studying the case of elliptic curves. In
this case examples are easy to give because in many cases
we can choose plane models (Weierstrass equations). 
These exercises can be used at motivation for more abstract methods
which apply for higher genus. You can do the exercises by explicit
methods and calculations.  

 For details on
elliptic curves, Weierstrass equations,  the $j$-function, and related
 issues, 
see \cite{Silv}, Chapters III and VII.

 A non-singular one-pointed curve of genus one is called
an elliptic curve. In other words: an elliptic curve is an algebraic
curve $E$ defined over a field, absolutely irreducible, non-singular,
of genus one, with a marked point $P \in E(K)$. Morphisms are supposed
to respect the marked point.
 
The following exercise is an easy exercise using the theorem of Riemann-Roch.

\begin{exercise} Show the following are equivalent:
\begin{enumerate}
\item $(E,P)$ is an elliptic curve over $K$.
\item $E \subset \PP^2_K$ is a plane, nonsingular cubic curve,
with a marked point $P \in E(K)$. 
\item $(E,P)$ is an abelian variety of dimension one over $K$.
\end{enumerate}
\end{exercise}

\begin{definition} Let $R = R_v$ be a discrete valuation ring, with
$K=$ fract$(R)$ its field of fractions, and $k = R_v/m_v$ the residue
class field.
\begin{enumerate}
\item 
An elliptic curve $E$ defined over $K$ is said to have {\it good reduction}
(at the given valuation) if there exists a smooth proper
morphism ${\cal E} \rightarrow \Spec(R)$ with generic fiber
isomorphic to $E  \rightarrow \Spec(K)$. If $E$ does not have good
reduction, we say that it has {\em bad} reduction.

\item We say  $E$ has {\it stable reduction} at $v$ if either it has good
reduction, or  there exists a {\em nodal} ${\cal E} \rightarrow \Spec(R)$
with generic fiber isomorphic to $E  \rightarrow \Spec(K)$.
\end{enumerate}
\end{definition}

\begin{definition}
  We say that $E$ has {\em potentially good reduction}, if there
exists a finite extension $K \subset L$, where $B$ is the integral
closure of $R$ in $L$, and $w$ a valuation over $v$, such that
$E \otimes L$ has good reduction at $w$.

We define {\em potentially stable reduction} analogously.
\end{definition}

\n
Here are some exercises to warm up:

\begin{exercise} Suppose $R=k[T]$, with char$(k) \not= 2, \not= 3$, and
let $E$ over 
$K = k(T)$ be given by the equation $Y^2 = X^3 + T^6.$
Show that $E$ has good reduction at the valuation given by $v(T)=1.$
\end{exercise}

\begin{exercise}
  Suppose $R=k[S]$, with char$(k) \not= 2,\not= 3 $,
and let $E$ over 
$K = k(S)$ be given by the equation $Y^2 = X^3 + S.$
Show that $E$ has bad reduction at the valuation 
$v$ given by $v(S)=1.$
\end{exercise}

\begin{exercise}
 Suppose $R=k[S]$, with char$(k) \not= 2,\not= 3 $, and
let $E$ over 
$K = k(S)$ be given by the equation $Y^2 = X^3 + S.$
Show that $E$ has potentially good reduction at the valuation 
$v$ given by $v(S)=1.$ 
\end{exercise}

\n
Suppose that $E$  is given over
$K$ by a Weierstrass equation with 
coefficients in $R$ (see \cite{Silv}, III). Such an equation defines
an {\em affine}
plane curve  
${\cal E} \subset  \AAA^2_B$
over $\Spec(R) = B$, and it is easy to see that the curve
$E_0 :=  {\cal E} \otimes_R k$ is irreducible and has at most
one singular point. The curve $E$ is obtained by adding the point at
infinity to $E_0$.  Suppose the Weierstrass equation is {\em  minimal} at
$v$. If this singular point is a cusp, we say that this  
reduction is of {\em additive type}, if it is a node we say that this
reduction is of {\em multiplicative type}, or we say in this case
the reduction is stable. 
ind a model which has stable reduction over a finite extension.

\begin{exercise}
Show that the notion of ``good reduction" as defined earlier is equivalent
by saying there is Weierstrass equation defining good reduction. Show
that a reduction of multiplicative type is a stable reduction. 
  \end{exercise}

\n
A reduction given by a minimal Weierstrass equation of additive type is
bad reduction which is non-stable; non-stable bad reduction is
sometimes called 
cuspidal reduction.

\begin{exercise}
  Suppose $R=k[T]$, with $\chara k \neq 2$, and let $E$
over 
$K = k(T)$ be given by the equation $Y^2 = X^2 \cdot (X-1).$
Show that any model of this curve given by a Weierstrass 
equation has stable reduction at the valuation given by 
$v(T) = 1$. Show that this curve does not have potentially
good reduction.
\end{exercise}

\begin{exercise}
 Let $R$ be a discrete valuation ring, with residue characteristic
$\neq 2$, and fraction field $K$. Let $E$ be an elliptic curve over $K$.
\begin{enumerate}
\item Show that after a suitable extension of $R$, the curve $E$
admits a minimal Weierstrass equation of the form 
$$ y^2 = x (x-1) (x-\lambda)$$ for some $\lambda\in R$.
\item   Conclude that this curve has
potentially stable reduction.
\end{enumerate}
\end{exercise}

\begin{exercise}\hardexercise\
 Let $R_v$ be a DVR, with residue characteristic
$\neq 3.$ Suppose $E$ is an elliptic curve over $K$ given
by a Weierstrass equation $E = V(F) \subset \PP^2_K$
such that all flex points of $E$ have coordinates in $K$.
\begin{enumerate}
\item Show that this curve admits a plane equation (not a Weierstrass
equation!) over  $K$  of the form $$\lambda (X^3+Y^3+Z^3) = 3\mu XYZ,$$
for some $(\lambda:\mu)\in \PP^1_{R_v}$.  
\item  Show that
$E \otimes _K L$ has stable reduction at $v$.
\item Show that $E$ has potentially stable reduction. 
\end{enumerate}
\end{exercise}

\n
{\bf Conclusion.} Every elliptic curve over a field $K$ with a
discrete valuation has potentially stable reduction at that
valuation.  

This is a special case of \ref{SRg}, the stable reduction
theorem for stable  curves of 
arbitrary genus.

\begin{exercise}\hardexercise
\begin{enumerate}
\item  Let $R_v$ be a DVR, and $E$ an elliptic curve of
the fraction field $K$. Show that $E$ has potentially good reduction at $R_v$
if and only if $j(E) \in R_v$.
\item Can you formulate (and prove?) the same result for
curves of arbitrary genus? 
\item If $E$ is an elliptic curve over a field $K$, and
$\End_K(E) \neq \ZZ$, then $E$ has potentially good reduction at every place
of $K$. 
\end{enumerate}
\end{exercise}

\subsection{Remarks about monodromy}

 Let $C$ be a non-singular curve over a
field $K$, and let $v$ be a discrete valuation of $K$. Consider
properties of good reduction, bad reduction at $v$, and so on. We have
quoted that $C$ has stable reduction at $v$ iff $J := \Jac(C)$ has
stable reduction, see \cite{D.M}, Proposition (2.3).  

However note that it may happen that $J$ has good reduction, and $C$
has bad reduction; this is the case if the special fiber $C_0$ of the
minimal model of $C$ at $v$ has a generalized jacobian $J_0 =
\Jac(C_0)$ which is an abelian variety. Such a curve $C_0$ is called a
curve of ``compact type", or a ``nice curve" (and sometimes called a
``good curve", but we do not like that terminology, because a curve
reducing to a ``good curve" may not have good reduction...). In this
case the special fiber $C_0$ is a tree of non-singular curves,
i.e. every irreducible component is non-singular, and in the dual
graph of $C_0$ there are no cycles. The easiest example is: a join of
two non-singular curves, each of genus at least one, meeting
transversally at one singular point. For example a curve of genus two
degenerating to a transversal crossing of two curves of genus one is
the easiest example. Here is another example: take $\PP^1$ with three
marked points, and attach three elliptic tails via normal crossings at
the markings, arriving at a nice curve of genus three 

Monodromy (action of the local fundamental group of the base on
cohomology) decides about the reduction of an abelian variety being
bad or good, see \cite{S.T}, Theorem 1 on page 493. In the analytic
context one can take the local fundamental group of a punctured disc
acting on cohomology; in all cases one considers the inertia-Galois
group of $v$ acting on $\ell$-adic cohomology, where $\ell$ is a prime
number not equal to the residue characteristic of $v$. 

Note that algebraic monodromy has {\it eigenvalues which are roots of
unity.} This was proved by Landman, Steenbrink, Brieskorn in various
settings, and we find a proof by Grothendieck in the appendix of
\cite{S.T}. For a sketch of that proof, see \cite{FO}, for further
references, see \cite{FO4}. 

{\it Algebraic monodromy is trivial iff $X=J$ has good reduction, iff
$C$ has compact type reduction} (which may be either good reduction or
bad but ``nice" as explained above). 

The algebraic monodromy is 
 unipotent (all eigenvalues are equal to
one) if and only if  $X$ has stable reduction, if and only if $C$ has
 stable reduction.  

But, how can we distinguish  for curves the difference between good
reduction and bad compact type reduction? As we have seen, this is not
possible via algebraic monodromy on cohomology. But, in a beautiful
paper, \cite{A.M.O} we find a method which for curves unravels these
subtle differences for curves: the local fundamental group of the base
acts via outer automorphisms on the fundamental group of the generic
fiber (again, here one can work in the analytic-topological context,
or in the $\ell$-adic algebraic context). {\it This action is trivial
iff $C$ has good reduction.}

%
        \section{Construction of moduli spaces}
\mylabel{Sec:Moduli-construction}
Early constructions of the moduli spaces of smooth curves $M_g$ included a
complex - analytic constructions via Teichm\"uller theory and via the
construction of moduli of abelian variety using locally symmetric spaces. These
constructions are not algebraic in nature and therefore cannot be generalized
to positive or mixed characteristics. 

A first algebraic approach, which  is still commonly used today, was given
by Mumford using his Geometric Invariant Theory \cite{GIT}. We will sketch
one version of this approach, due to Gieseker, which automatically gives also
the moduli spaces of stable curves $\oM$.
There is 
another approach,
due to Artin and Koll\'ar \cite{Kollar}, which circumvents the use
of Geometric Invariant 
Theory. Nowadays both approaches work over $\ZZ$.

How does one start? It is evident that if we want to parametrize {\em
all} stable curves of a certain genus, we had better have {\em some}
family of curves in which all these curves appear. We know of two
general approaches for that. One method uses parameter spaces for
curves embedded in projective space, such as Hilbert schemes (or Chow
varieties). We will 
follow this approach. The other approach, due to Artin \cite{Artin},
uses versal 
deformation spaces. It works in greater generality but involves a
number of technicality which we would rather avoid here. 

It is easy to see
that for any stable curve $C$ of genus 
$g>1$, and any $\nu\geq 3$, the $\nu$-canonical series $H^0(C,\omega_C^\nu)$
gives an embedding of $C$ as a curve of degree $d:=\nu(2g-2)$ in
a projective space of dimension $N:=\nu(2g-2)-g$. Thus the Hilbert scheme
$ \cH_{\PP^N,P(T)}$ (over
$\ZZ$!) 
parametrizing subschemes of $\PP^N$ with Hilbert polynomial $P(T) := dT +1-g$
carries a universal family $\cC_{\PP^N,P(T)} \to \cH_{\PP^N,P(T)}$ in which
each stable curve of genus $g$ appears at least once.

There are two problems with this family:
\begin{enumerate}
\item Each curve appears more than once in the family. Indeed, the embedding of
the curve $C$ in $\PP^N$ involves two choices: a choice of a line bundle of
degree $d$, and a choice of a basis for the linear series. And of course the
curves could also be embedded in a projective subspace using a subseries.
\item There are many curves in $\PP^N$ with Hilbert polynomial $P(T)$ which are
far from stable.
\end{enumerate}

\n
Since a nodal curve can only deform into nodal curves, it is easy to see that
there is an open subset $\cH_{st} \subset\cH_{\PP^N,P(T)}$ which parametrizes
{\em stable} curves, embedded by a {\em complete} linear system in
$\PP^N$. Denote the restriction of the universal family to $\cH_{st} $ by
$\pi:\cC_{st} \to\cH_{st}$.  Considering the locus in $\cH_{st}$ where
$R^1\pi_*(\cO(1) \otimes \omega_{\cC_{st} \to\cH_{st}}^\nu)$ jumps in
dimension, we immediately see that there is a {\em closed} subscheme
$\cH_g\subset \cH_{st}$ parametrizing stable curves embedded by a complete
$\nu$-canonical series. The restriction of the universal family will be denoted
$\cC_g\to \cH_g$. 

There is a natural action of the projective linear group $PGL(N+1)$ on $\cH_g$
via changing coordinates on $\PP^N$. It is easy to see that the ``ambiguity''
for choosing the embedding of a curve $C$ in the latter universal family is
fully accounted for by the action of this group. In other words, stable curves
correspond in a one-to one manner with $PGL(N+1)$ orbits in $\cH_g$. Thus, at
least set theoretically, $\oM = \cH_g / PGL(N+1)$. 

\subsection{Geometric Invariant Theory and Gieseker's approach}

We arrived at the following questions:
\begin{enumerate} 
\item Does the quotient $ \cH_g / PGL(N+1)$ exist as a scheme?
\item Can we show that it is projective?
\item Does it satisfy the requirements of a coarse moduli scheme?
\end{enumerate}

Geometric Invariant Theory is a method which allows one to approach the first
two questions simultaneously. The third question then becomes an easy gluing
exercise.

The general situation is as follows: Let $X\subset \PP^n$ be a quasi-projective
scheme and suppose $G$ is an algebraic group acting on  $\PP^n$ and stabilizing
$X$. One wants to know whether or not a quotient $X/G$ exists as a scheme and
whether or not it is projective. 

A natural approach is to look for a space of invariant sections of some line
bundle. Thus assume that the action of $G$ on $\PP^n$ lifts to
$\cO_{\PP^n}(l)$. Then it also lifts to any power $\cO_{\PP^n}(l\cdot m)$, and
we can look at the ring of invariants $R := \oplus (\cO_{\PP^n}(l\cdot
m)^G$. We have a natural rational map $q:\PP^n \das R$. We would like to know
whether or not this map is well defined along $X$, and what the image is like. 

First, an easy observation. For any point $x\in X$, the map $q$ is well defined
at $x$ if and only if there exists a nonconstant invariant $f\in R$ such that
$f(x) \neq 0$.  

We want to check whether $q$ is a quotient map at $x$.
To go any
further, we need to assume that the group $G$ is {\em reductive}. Assuming that
$G$ is reductive, then the question whether map $q$ is a quotient map at a
neighborhood of $q$ 
can be translated to a question about the closure $\overline{Gx}$ of the orbit
of $x$: one needs to check that for any point $y\in \overline{Gx}\setmin Gx$
there is an invariant $f\in R$ which vanishes at $y$ but not at $x$. A point
$x$ is called {\em GIT-stable} if it satisfies this condition.

Mumford's numerical criterion for stability (see \cite{GIT}) gives a
 way to check GIT-stability in some situations. 

Let us consider our situation. The scheme $\cH_g$ is quasi projective
- from its 
construction one sees that it
naturally sits inside a Grassmannian, which has a Pl\"cker embedding in
some $\PP^n$. It is easy to see that the action of $PGL(N+1)$ extends to
$\PP^n$, and lifts to some line bundle $\cO_{\PP^n}(k)$. Applying this
criterion systematically,
Gieseker verified in \cite{Gies} that
\begin{enumerate}
\item If a point $x\in \cH_{\PP^N,P(T)}$ corresponds to a scheme which is not
a stable curve, or to a curve which is not embedded by a complete linear
series, then {\em every} nonconstant invariant vanishes at $x$. 
\item If a point $x\in \cH_{\PP^N,P(T)}$ corresponds to a stable curve embedded
by the complete $\nu$-canonical linear series, then $x$ is GIT-stable.
\end{enumerate}

Using the two statements, and the fact that $G$ is reductive, it is not
difficult to realize that  
\begin{enumerate}
\item the map $\cH_g \to \PP^n$ is a quotient map, and
\item the image of $\cH_g$ is projective.
\end{enumerate}
This proves the existence and projectivity of $\oM$.

\subsection{Existence of $\onM$} There is no known analogue of Gieseker's
result for stable pointed curves. It is not difficult to construct a
Hilbert-type scheme for stable pointed curves, with a reductive group
action, and such 
that the quotient is set-theoretically $\onM$. But in order to tell that the
quotient is isomorphic to  $\onM$ as a scheme, we first need to
construct $\onM$ 
in some other way.

But there is a very useful trick, which reduces the
construction of $\onM$ to the existence of $\oM$ for some larger value of $g$. 
 We  give the reduction over a field, but it works similarly over $\ZZ.$

Fix $n$ irreducible stable curves $C_i$ of genus $g_i>g$, all nonisomorphic to
each other, and fix a rational point $x_i \in C_i$. For any stable $n$-pointed
curve $(C,p_1,\ldots,p_n)$ of genus g, we can construct a stable curve $C'$ of
genus $g'=g+\sum 
g_i$  as follows: $C' = (\cup C_i ) \cup C $, where we glue together $C$ and
$C_i$ by identifying $p_i$ with $x_i$.

Clearly this construction gives a set theoretic embedding $\onM \to
\overline{M_{g'}}$. The image set is easily seen to be a scheme, and
by working 
the construction in a family it is easy to see that it is a coarse moduli
scheme.


	\section{Existence of tautological families}
\mylabel{Sec:Tautological-pointed}

For almost any application of moduli spaces of curves including the
alteration theorem, it is necessary to know that there exists a
family $C \to M$ over a scheme $M$ such that the associated morphism
to the moduli space is finite and surjective. Such a
family is called a {\em tautological family}; see \ref{taut}. 
Various authors
have devised general methods of showing this, but for the moduli spaces of
curves there is a ``very nice'' way to find such a cover, using level
structures. The case of the moduli space $\oM$ of stable (unpointed) curves is
discussed in detail in \cite{vG.FO}. In this section we describe how this can
be generalized for sable {\em pointed} curves as well. We rely throughout on
the treatment in \cite{vG.FO}. 
In Section \ref{kollar's-method}
we outline another way to construct tautological
families, which works in greater generality.

\subsection{Hilbert schemes and level structures}
Fix: 
\begin{itemize}
\item  an integer $g \in \ZZ_{\geq 0}$ (the genus),
\item   an integer $n \in \ZZ_{\geq 0}$ (the number of marked points), \item
such that $2g-2+n > 0$,  
\item   and an integer $m \in \ZZ_{\geq 1}$ (the level).
\item Fix an integer $\nu \in \ZZ_{\geq 5}$, which will be used to study
$\nu$-canonical embeddings of curves into a projective space. 
\end{itemize}
\begin{remark} If $n=0$ or $m=1$ these data will be omitted from the notation,
e.g. $M_{g,0} = M_g$. If $g=0$, the level structure is irrelevant,
$M_{0,n}^{(m)} = M_{0,n}$.
\end{remark}

\n
Let $C$ be  a curve whose jacobian is an abelian variety. By a level-$m$
structure on $C$ we mean a 
symplectic level structure as explained in \cite{vG.FO}. If  a
level-$m$ structure is considered we assume that all schemes, varieties are
over a base on which $m$ is invertible, i.e. are schemes over
$\Spec(\ZZ[1/m])$. 

Recall that there is a Hilbert scheme $H_P$ parametrizing curves $C\subset
\PP^N$, where $N = \nu{\cdot}(2g-2+n)-g$, with Hilbert polynomial $P(t) =
 \nu{\cdot}(2g-2+n)\cdot t - g + 1$. We want  to find a scheme parametrizing
{\em pointed} curves - this is done in a standard way as follows. Observe that
there is a closed subscheme 
$H_{P,n}\subset H_P \times (\PP^N)^n$ parametrising pairs $(C,
(p_1,\ldots,p_n))$ where $p_i \in C$.  There is an open subscheme
$H_{st} \subset H_{P,n}$ where the curves are nodal, the points are distinct
and regular points on the curves, and the pairs  $(C,
(p_1,\ldots,p_n))$ are stable. Last, there is a closed subscheme
$H_{g,n}\subset H_{st}$ where the embedding line bundle of $C\subset \PP^N$ is
isomorphic to $(\omega_C(p_1+\ldots +p_n))^\nu$. 

Over $H_{g,n}$ there is a universal family $C_{g,n} \to H_{g,n}$ with sections
$s_i:  H_{g,n} \to C_{g,n}$ of stable pointed curves, embedded in $\PP^N$ by
the chosen line bundle.  The linear group $\PGL = PGL(N)$ acts on $C_{g,n} \to
H_{g,n}$ equivariantly, and 
$\onM  =   H_{g,n}\mbox{\large $/$} \PGL$ 
is the quotient. 

 Note that there is an open subset $H^0_{st} \subset
H_{g,n}$ parametrizing {\em smooth} stable pointed curves.

\subsection{Moduli with level structure}

\begin{theorem}\mylabel{mg}
 For $m \geq 3$, and $2g-2+n > 0$, there exists a {\em fine} moduli scheme
 $M_{g,n}^{(m)}$ for smooth stable $n$-pointed curves with level-$m$ structure.
 In particular there
 exists a  universal curve with level structure over $M_{g,n}^{(m)}$. 
 This moduli scheme is  smooth over $\Spec(\ZZ[1/m])$. 
\end{theorem}

\n
Note in particular that $\mnM$ is a normal scheme, and that $\mnM \to \nM$ is a
Galois cover with Galois group $\Syp(2g,\ZZ/m)$.   

We use the notation $S_m := \Spec(\ZZ[1/m])$.

\begin{definition} Let $g \in \ZZ_{\geq 1}$. Fix  $n \in \ZZ_{\geq 0}$, with
$2g-2+n>0$. For any $m \in \ZZ_{\geq 3}$, the scheme  
$$\overline{M_{g,n}^{(m)}} \longrightarrow S_m$$
is defined as the normalization of $\onM[1/m] = \onM \times_{\ZZ} S_m$ in
$\mnM$.  
\end{definition}

\n
For simplicity of notation in this section, we write $M =
\overline{M_{g,n}^{(m)}}$ and $M^0 = \mnM\subset M$. 

\begin{theorem}\mylabel{omg} Fix $g, n$, and $m$ as above. Suppose $m
\geq 3$.   There exist   
a stable $n$-pointed curve $(\cC,\cP) \to  M$, and a level-$m$-structure
$\alpha$ on $\cC ^0 := \cC|_{M^0}$ such that
$$ (\cC,\cP) \to  M \quad\mbox{\it is tautological for}\quad M \to \onM,$$
 and such that
$$ (\cC ^0, \cP ^0,\alpha) \to M^0$$
represents the functor $\cmnM$.
\end{theorem}

We give an argument for \ref{mg} and  \ref{omg} following the line of
\cite{vG.FO}. This is a kind of ``boot-strap'' argument, which uses the idea
that {\em once one quotient space exists, many others follow}. We also sketch
another argument which  reduces the problem to the case of \cite{vG.FO}.

 There is a relative jacobian scheme $J(C^0_{g,n})
\to H^0_{g,n}$. This is an abelian scheme, so we can look at its
group-subscheme of $m$-torsion points. Taking a symplectic rigidification of
this group scheme we arrive at 
$H_{g,n}^{(m),0}$ - the Hilbert scheme of smooth stable $n$-pointed 
curves with symplectic level-$m$ structure - embedded in projective space as
above.  

The action of  $\PGL$ on $ H^0_{g,n}$ clearly lifts to $H_{g,n}^{(m),0}$.
This immediately implies that 
$$ 
\mnM\ =\ H_{g,n}\ \hbox{\large $/$} \ \PGL
$$  
exists, since it
is 
finite over $\nM$.  By Serre's lemma this action has no fixed points, and it
also lifts 
 to $C_{g,n}^{(m),0} = C_{g,n}\times_{H_{g,n}}H_{g,n}^{(m),0}$. 
This means that the quotient $\PGL \mbox{\large $\backslash$} C_{g,n}^{(m),0}
\to  \mnM$ is a universal family of smooth stable pointed curves with level
structure. 

 This proves Theorem \ref{mg}. 

The normalization of $H_{g,n}$ in $H_{g,n}^{(m),0}$ will be denoted by
$H_{g,n}^{(m)}$. The argument of \cite{vG.FO}, (2.6)  works word for
word, and shows  that $\PGL$ still acts without fixed
points on $H_{g,n}^{(m)}$. This gives the existence of  the quotient 
$$ 
\omnM \ =\ H_{g,n}^{(m)}\ \hbox{\large $/$}\ \PGL.$$ 
Again
the universal family over $H_{g,n}^{(m)}$ descends to a family over $\omnM$,
this extends the universal family over $\mnM$, and clearly it is
tautological. This proves  Theorem \ref{omg}. \qed

\subsection{Proof by reduction to the unpointed case} Starting from 
$\omM$ and its tautological family we can construct 
$\omnM$ and its tautological family by induction on the number of points $n$ in
the manner described below.  

Denote by $D \to \omnM$ the tautological family. It is easy to see that in fact
$D = \overline{M_{g,n+1}^{(m)}}$. So $D\times_{\omnM} D \to
\overline{M_{g,n+1}^{(m)}}$ is a family of stable $n$-pointed curves with level
structure, but with an additional section given by  the diagonal. 
 Using the stabilization
process as described in \cite{Knudsen} (see Section
\ref{Sec:Knudsen} above) one blows this scheme up, to obtain 
the  tautological family over $\overline{M_{g,n+1}^{(m)}}$  as desired.

\begin{remark}
The moduli space $\mnM$ is smooth over $S_m$ for $m \geq 3$; this follows from
Serre's lemma and deformation theory. However, the moduli space
$\omnM$ is singular if $g > 2$;  Serre's lemma holds also 
in this situation, but the space is not the coarse (or fine) moduli
space of a moduli functor whose deformation spaces coincide with the
deformations of stable curves. For more explanation, see \cite{Mostafa} or
\cite{vG.FO}. 
\end{remark}

\n
The argument above works for $g>1$ when $\omM$ exists.
For rational curves these theorem are relatively easy, and known, since the
moduli spaces are fine moduli spaces in genus 0. For 
elliptic curves these theorems are known by the theory of modular curves.

\subsection{Artin's approach via slicing} \mylabel{kollar's-method}
A general approach for constructing tautological families over finite
covers of coarse
moduli space was developed by Artin (see description in
\cite{Kollar}). Here we present a version of this 
approach adapted to stable pointed curves.

\vn
{\bf Step 1: slicing}
Consider the locally closed subset of the Hilbert scheme $H_{g,n}$
discussed 
above. It carries a universal family of stable pointed curves
$C_{g,n}\to H_{g,n}$ suitably
embedded
in a projective space. This family induces a natural morphism
$H_{g,n}\to \onM$. The fibers coincide with the  $G$-orbits
associated to the
embedded curves, where $G=\PGL$.

Fix a point $x\in H_{g,n}$. By repeatedly taking hyperplane sections, we
can find a locally closed subscheme $V_x\subset H_{g,n}$ such that

\begin{enumerate}
\item $Gx\cap  V_x\neq \emptyset$;
\item If $x'\in H_{g,n}$ and $Gx'\cap V_x\neq \emptyset$, then there exists a
neighborhood $x'\in U$ such that for any $y\in U$ we have 
that $Gy\cap V \neq \emptyset$; and
\item for any $y\in H_{g,n}$ we have that $V_x\cap Gy$ consists of
finitely many 
closed points.
\end{enumerate}

These
$V_x$ are  ``multi-sections'' of the map $H_{g,n}\to \onM$ in a
neighborhood of $Gx$. 
The essential point is that all orbits in $H_{g,n}$ are of
the same dimension.

Using the Noetherian property, we can choose finitely many of these, say
$V_1,\ldots,V_l$, such that every orbit meets at least one of them.

\vn
{\bf Step 2: normalization.}
Let $K$ be the join of the function fields of $V_i$ over $\onM$. Let
$V$ be the normalization of $\onM$ in the Galois closure of $K$. The
scheme $V$ admits many rational maps to the $V_i$. It is not hard to
see that for every point $v\in V$ at least one of these maps is well
defined at $v$! Pulling back the  
families on $V_i$, we see that $V$ is covered by open sets, each of
which carries a  family of stable pointed curves, compatible with the
given morphism $V \to \onM$.

\vn
{\bf Step 3: Gluing.}
Now we can use Lemma \ref{Lem:Extend-isom} inductively. We obtain a
finite surjective $M \to V$ over which the families glue
together to a family $C \to M$ such that the associated moduli
morphism is the composition $M \to V \to \onM$. Since $V$ is finite
over  $\onM$, this forms a tautological family.  
 
\begin{remark}
It is not hard to construct a tautological as above {\em without using
the existence of $\onM$!} One can use this to construct the moduli
space ``from scratch'' as a proper algebraic space, which is roughly
speaking a quotient of a 
scheme by a finite equivalence relation. Koll\'ar in \cite{Kollar}
has shown how to use this to prove, without GIT, that $\onM$ is
projective. 
\end{remark}
%
%


%
		\section{Moduli, automorphisms, and families} 
\mylabel{Sec:Moduli-other}

This section will not be needed in the proofs above. The central theme
here is the relationship between automorphisms, coarseness of moduli,
and the existence of  families. The main principle which will emerge
is:
\begin{center} a moduli space $M$ is a fine moduli space
$$ \Updownarrow $$
objects parametrized by $M$ have no nontrivial automorphisms
 $$ \Updownarrow $$
$M$ carries a unique tautological family.
\end{center}

We also touch on the issue of
singularities of moduli spaces.


 For rational curves, and $n \geq 3$,  the moduli schemes $M_{0,n}$
 and $\overline{M_{0,n}}$ exist, these are smooth over $\Spec(\ZZ)$,
 these are fine moduli schemes, i.e. they carry a universal family.

However, the moduli space $M_{1,1}$ and the moduli spaces $\M$ for  $g
>1$ are not  fine for the related moduli functor. 

\begin{exercise}[Deuring]  Let $K$ be a field, let $x \in K$. Then
there exists an elliptic curve $E$ defined over $k$ with $j(E) = x$.
[ Suppose char$(K) \not= 2, \not= 3$, suppose $E$ is is given over $K$

by the equation $Y^2 = X^3 + AX + B,$ with $4A^3 + 27B^2 \not= 0$. 
Then define
$$j(E) := 1728{\cdot}4{\cdot A^3} / (4A^3 + 27B^2).$$
For the definition of the $j$-invariant, see \cite{Silv}.]
\end{exercise}

\n
This can partly be made more precise as follows:

\begin{exercise}  Consider $M_{0,1} \cong \bfA^1_{\ZZ}$, and remove
the sections $j=0$ and $j = 1728$: 
$$U := \bfA^1_{\ZZ} \setmin \{0,1728\}_{\ZZ}.$$ 
There exists a tautological curve
$${\cal E} \rightarrow U.$$
\begin{enumerate} \item 
 This cannot be extended over any of the deleted points. 
\item  This family is not at all unique. 
\end{enumerate} 
\end{exercise}

\begin{exercise} Consider
$U:=\CC \setmin \{0,1728\}$. Show: up to isomorphisms there exist exactly 4
tautological curves (stable, one pointed smooth curves of genus 1 with
$j$ invariant different from 0 and 1728)
over this moduli space. 
Show that for the ground field $K = \QQ$ there are {\em infinitely many}
tautological curves over the moduli space $\bfA^1_{\QQ} \setmin
\{0,1728\}$. Characterize them all. 
\end{exercise}

\n 
We have seen the difference between a universal curve and a
tautological curve: the moduli problem for elliptic curves with
geometrically no non-trivial automorphisms admits a coarse moduli
scheme; over that scheme there is a tautological curve, but the scheme
is not a fine moduli scheme (not every family is a pull-back  from one
chosen tautological curve). Here is another example: 

\begin{definition} A curve $\cC \to S$ is called a {\it
hyperelliptic curve} if it is smooth, of relative genus $g$ with $g
\geq 2$, and if there exists an involution $\iota \in \Aut(\cC/S)$
such that the quotient $\cC/<\iota> \to S$ is a smooth family of
rational curves.
\end{definition}

\begin{remark} Elliptic curves and rational curves are not called
``hyperelliptic", but sometimes the terminology ``quasi-hyperelliptic"
is used for curves having an involution with rational quotient. 
\end{remark}

\begin{theorem} Consider the moduli space  $\operatorname{Hip}_g$ of
hyperelliptic 
curves of genus $g \geq 2$ (even over $\CC$).  
  If $g$ is even there {\em does not exist}  a curve
defined over the function field $\CC(\operatorname{Hip}_g)$ having as
moduli point the 
  generic point of 
$\operatorname{Hip}_g$. 
\end{theorem}

(See Shimura \cite{Shimura}, Theorem 3.)

\vn
In different terminology: For no open dense subset $U \subset H_g$
does there exist a tautological curve when $g$ is {\it even}.  

There does exist a open dense subset $U \subset H_g$ and  a
tautological curve ${\cal C}_U$ when $g$ is {\it odd}.  

\begin{corollary} No dense open subset in
$M_2$ or in  $M_2 \otimes K$ carries a tautological curve.
\end{corollary}

\begin{exercise} Choose $g \in \ZZ_{> 2}$, and consider nonsingular
curves of genus $g$. 
\begin{enumerate}
\item Show that there exists such a curve which has no nontrivial 
automorphisms.
\item (variant:) Show that a general curve of genus $>2$ has no nontrivial 
automorphisms.
\end{enumerate}
\end{exercise}

\begin{remark} There is a morphism $M_{g,n+1} \to
M_{g,n}$ ("forgetting the last marking"). Sometimes this is called the
``universal curve over $M_{g,n}$",  but we think in general this
terminology is not justified in all cases possible. 
\end{remark}

\begin{theorem} Let $U \subset M_g$ with $g \geq 3$ fixed, be the set
of points corresponding with curves which have geometrically no
non-trivial automorphisms. This set is dense and open. Let ${\cal
M}_U$ be the corresponding moduli functor. This functor is representable.
\end{theorem}

\n
In other terminology: there does exist a (unique) universal curve
${\cal C}_U \to U$  for the moduli problem of curves of genus $g \geq
3$ with geometrically no non-trivial automorphisms.

In particular: Let $K$ be a field, $g \in \ZZ_{\geq 3}$, and $\eta$ be
the generic point of $\M \otimes K$. There exists an algebraic curve
defined over $K(\eta)$ having $\eta$ as moduli point. However the
universal family as indicated above over $U \subset M_g$ does not
extend to any smooth family of curves over $\M$.     

\begin{exercise} Formulate and prove a generalization of
previous theorems to the case of stable pointed curves. 
\end{exercise}

\begin{exercise} Let $n>2g+2$ and let $(C,P_1, \cdots ,
P_n)$ be any stable $n$-pointed curve 
of genus $g$. Suppose that $C$ is {\it regular} (and hence irreducible).
Show that $$\Aut((C,P_1, \cdots , P_n)) = \{1\}$$
(if you want, assume that char$(k) = 0$). 
\end{exercise}

\begin{exercise}
 Let $g\in \ZZ_{\geq 1}$ and $2-2g < n
\leq 2g+2$ and $0 \leq n$. Show that $\nM$ is a coarse, but not a fine
moduli space. 
\end{exercise}

\begin{exercise}
 Choose $g\in \ZZ_{\geq 0}$, and let
$n>2g+2$. Show that $\nM$ is a fine moduli space. Show that the
universal curve over $\nM$ is not smooth if $n \geq 2$. 
\end{exercise}

\begin{exercise}
 Consider all stable $n$-pointed curves of
genus $g$. 
Suppose that
$$2g-2+n \geq 3.$$
\begin{enumerate}
\item Show that there exists such a curve which has no nontrivial 
automorphisms.
\item (variant:) Show that a general  curve as above has 
no nontrivial 
automorphisms.
\end{enumerate}
\end{exercise}

\begin{exercise}
 Choose some $g$ (e.g. $g = 3$), choose a
very large 
integer 
$n$ (e.g. $n = 1997$), and construct a stable $n$-pointed curve of 
genus $g$ which has a non-trivial group of automorphisms.

Variant: Let $2g-2+n \geq 2$; show that there exist stable 
$n$-pointed curves of  genus $g \geq 3$ in codimension two in the
moduli space with  non-trivial groups of automorphisms.
\end{exercise}

\begin{exercise}
 Let  $g\in \ZZ_{\geq 1}$, and $n > 2-2g$
and $n \geq 0$. Show that $\onM$ is not a fine moduli space.  
\end{exercise}

 \n
Choose $2g-2+n > 0$, choose $m \geq 1$ and let $M$ be
one of the following spaces: $\mnM$, or $\onM$ (all these spaces are
defined by a moduli functor). Let $x \in M(k)$ be a geometric point,
and let $X_0 := (C,P,\alpha)$ be the corresponding object over $k$ (if
$C$ is non-smooth there is no level structure, the genus of $C$ is
$g$, we have $P = \emptyset$ if $n=0$, we have $\alpha = id$ if
$m=1$).  Let $D = \Def(X_0)$ be the universal deformation space;
i.e. consider $\Lambda = k$ if char$(k)=0$, and $\Lambda =
W_{\infty}(k)$ in case of positive characteristic, consider all local
artin $\Lambda$-algebras, and consider the object prorepresenting all
deformations of $X_0$ over such algebras (see
\cite{Schlessinger}). This universal deformation 
object exists, and it is formally smooth over $\Lambda$ on $3g-3+n$
variables; in case $n=0$ this can be found in \cite{D.M}, page 81,
the case of pointed curves follows along the same lines; in case
$m>1$, we have required that $m$ is invertible in $k$, finite, flat
group schemes of $m$-power order on such bases are \'etale, and
deformations of level structures are unique by EGA IV$^4$, 18.1. Let
$G := \Aut(X_0)$. Note that $G$ {\it is a finite group} (because we
work with stable curves). Note that $G$ acts in a natural way on  $D =
\Def(X_0)$ by ``transport of structure".

\begin{theorem} In the cases described, the formal completion of
$M$ at $x$ is canonically isomorphic with the quotient

$$ \Def(X_0) \  \hbox{\large $/$} \ G \ 
\stackrel{\sim}{\longrightarrow} \
M^{\wedge}_x.$$ 
\end{theorem}
This is well-known, e.g. see \cite{H.M}, \S 1.  

\begin{exercise}
 (Rauch, Popp): Let $g \in \ZZ_{\geq 4}$,
and let $A \subset \M$ be an irreducible component of the set of all
points corresponding with curves with non-trivial automorphisms. Show
that the codimension of  $A \subset \M$ is $\geq 2$. (In
positive characteristic this is also correct, but you might need some
extra insight to prove also those cases.) 
\end{exercise}

\begin{remark} Stable {\it rational} pointed curves have
no non-trivial automorphisms. For {\it elliptic} curves there are
curves with more than 2 automorphisms in codimension one. For curves
{\it of genus two} we find a description of all curves with ``many
automorphisms" in \cite{Ig}. Note that hyperelliptic curves of genus
three are in codimension one.\\ 
\end{remark}

\begin{exercise}
  Show that non-hyperelliptic curves of genus three
with non-trivial automorphisms are in codimension at least two. 
\end{exercise}

\begin{exercise}
 (Rauch \cite{Rauch}, Popp
\cite{Popp3}):  Let $g \in \ZZ_{\geq 4}$, and let $[C] = x \in M_g$
be a geometric point. Show that $x$ is a singular point on $M_g$ iff
$\Aut \not= \{id\}$. 
[You might like to use: \cite{A.K}, Coroll. 3.6 on page 95: A
quasi-finite local homomorphism of regular local rings having the same
dimension is flat. Also you might like to use purity of branch locus:
a ramified {\it flat} covering is ramified in codimension one.] 
\end{exercise}

\begin{remark} For singularities of $M_2$ see
\cite{Ig}. Show that for genus three non-hyperelliptic points are
singular iff there are non-trivial automorphisms, e.g. see
\cite{FO2}. For singularities of moduli schemes of abelian
varieties, see \cite{FO3}. 
\end{remark}

\begin{remark} As we have seen in \cite{vG.FO}, the
moduli schemes $\omM$ have singularities for all $g \geq 3$ and $m
\geq 3$ (these spaces cannot be handled with the methods just
discussed, these spaces are not given by ``an obvious" moduli functor
!). As Looijenga, see \cite{Loo}, in characteristic zero,  and
Pikaart and De Jong, see \cite{P.dJ} showed, there exist a finite
map $M \to \M$ with $M$ regular (using non-abelian level structures)
(it is even true that $M$  is smooth over $\QQ$, or smooth over
$\ZZ[1/r]$ for some natural number $r>1$). 
\end{remark}

\noindent{\bf Summary} about
$$\mnM \hookrightarrow \omnM \longrightarrow \Spec(\ZZ[1/m]) =: S_m$$
for 
$$g \in \ZZ_{\geq 0}, \ \ n \in \ZZ_{\geq 0}, \ \  m \in \ZZ_{\geq 1},
\quad\mbox{with}\quad  2g-2+n > 0,$$ 
$\mnM$ and $\onM$ exist as coarse moduli schemes, we have constructed $\omnM$.
We have seen:
\begin{itemize}
 \item  For $g \geq 2$ the coarse moduli scheme $\M \to S =
 \Spec(\ZZ)$ exists. These are not fine moduli spaces. They do not
 carry a tautological family. For every $g$ this is singular.  
\item For $g \geq 2$ the coarse moduli scheme $\oM \to S = \Spec(\ZZ)$
exists. These are not fine moduli spaces. They do not carry a
tautological family. They are singular. 
\item A dense open set in $M_{1,1}$ carries a tautological family, and
it is not universal. 
\item No dense open set in $M_2$ carries a tautological family.
\item For  $g \geq 3$ a dense open set in $\nM$ carries a universal family.
\item For $n \geq 3$ the moduli spaces $M_{0,n} \subset
\overline{M_{0,n}}$ exist, they are fine moduli spaces, they are
smooth over $S = \Spec(\ZZ)$. 
\item For $2g-2+n > 0$, and $m \geq 0$ the moduli spaces $\nM \to S$,
and $\onM \to S$ and $\mM \to S_m$ exist, they coarsely represent a
moduli functor. For $n > 2g+2$ the moduli space $\nM$ is fine, and
smooth over $\Spec(\ZZ)$ (but the universal family is not smooth for
$n>1$). For $m \geq 3$ the space $\mM$ is fine and smooth over $S_m$. 
\item For $2g-2+n > 0$, and $m \geq 0$ there is a moduli space, and a
tautological family, with properties as in \ref{omg}. For $g \geq 3$
the morphism $\omnM \to S_m$ is not smooth. 
\end{itemize}

\newpage
            
\bigskip\bigskip\bigskip
\noindent
\parbox{8cm}{
Dan Abramovich	 		\newline
Department of Mathematics	\newline
Boston University		\newline
111 Cummington Street		\newline
Boston, MA 02215		\newline 
USA 				\newline
{\tt abrmovic@math.bu.edu	\newline
http://math.bu.edu/people/abrmovic}
}
\parbox{8cm}{
Frans Oort 			\newline 
Mathematisch Instituut 		\newline
University of Utrecht 		\newline
Budapestlaan 6 			\newline
NL - 3508 TA UTRECHT 		\newline
The Netherlands 		\newline
{\tt oort@math.uu.nl		\newline
http://www.math.uu.nl/staff/oort.html}
}

\begin{thebibliography}{XXXXXX}
%


\bibitem[Alteration]{dj}   
A. J. de Jong, {\em Smoothness, semistability, and
alterations}, Publications Math\'ematiques I.H.E.S. {\bf 83}, 1996, pp. 51-93.

\bibitem[GIT]{GIT}  
D. Mumford, J. Fogarty and F. Kirwan {\em Geometric invariant
theory.} 
Springer, Berlin, 1994. 



\bibitem[HAG]{HAG}   
R. Hartshorne, {\em  Algebraic geometry.} Springer, New York, 1977.

\bibitem[Red Book]{RedBook}  
D. Mumford, {\em  The red book of varieties and schemes.} Lecture Notes in Math., 1358, Springer, Berlin, 1988.

\bibitem[SGA]{SGA}  
A. Grothendieck (with M. Raynaud and D. S. Rim), {\em   Groupes de monodromie
en g\'eom\'etrie alg\'ebrique  I.} (S\'eminaire de g\'eom\'etrie
alg\'ebrique du Bois-Marie.) Lecture Notes in Math., 288, Springer,
Berlin, 1972.  

\bibitem{Abramovich}  
D. Abramovich,  {\em A high fibered power of a family of varieties of
general type dominates a variety of general type,} 
Invent. Math. {\bf 128} (1997), no.~3, 481--494.

\bibitem{aj}  
 D. Abramovich and A. J. de Jong, {\em
Smoothness, semistability and toroidal geometry},  J. Algebraic
 Geom. {\bf 6} (1997), no.~4, 789--801. 

\bibitem{ak}  
D. Abramovich and K. Karu,  {\em Weak semistable reduction in
characteristic 0}, preprint. {\tt alg-geom/9707012}

\bibitem{aw}   
 D. Abramovich and J. Wang, {\em Equivariant
resolution of singularities in characteristic 0},
 Math. Res. Lett. {\bf 4} (1997), no.~2-3, 427--433. 


\bibitem{A.K}  
A. Altman and S. Kleiman, {\em Introduction to Grothendieck duality theory.}
Lecture Notes in Math., 146, Springer, Berlin, 1970.

\bibitem{Artin}  
M. Artin, {\em Versal deformations and algebraic stacks,} Invent. Math. {\bf
27} (1974), 165--189.

\bibitem{A.W}  
M. Artin and G. Winters, {\em Degenerate fibres and reduction of
curves.} Topology {\bf 10} (1971), 373-383.
 
\bibitem{A.M.O}  
M. Asada, M. Matsumoto and T. Oda, {\em Local monodromy on the fundamental
groups of algebraic curves along a degenerated stable curve.} Journ. Pure
Appl. Algebra {\bf 103} (1995), no.~3, 235--283.


\bibitem{Behrend-Manin}  
K. Behrend and Yu. Manin, {\em Stacks of stable maps and Gromov-Witten
invariants,} Duke Math. J. {\bf 85} (1996), no.~1, 1--60 

\bibitem{Berthelot}  
P. Berthelot, {\em Alt\'erations de vari\'et\'es alg\'ebriques
(d'apr\`es A. J. de Jong),}  Ast\'erisque No. 241 (1997), Exp.\ No.\
815, 5, 273--311. 

\bibitem{Bie-Mil}  
E. Bierstone and P. D. Milman, {\em Canonical desingularization
in characteristic zero by blowing up the maximum strata of a local
invariant,} Invent. Math. {\bf 
128} (1997), no.~2, 207--302

\bibitem{bp}   
F. Bogomolov and T. Pantev, {\em Weak Hironaka theorem},
Math. Res. Lett. {\bf 3} (1996), no.~3, 299--307.
  
\bibitem{Brylinski}  
J-L. Brylinski, {\em  Propri\'et\'es de ramification \`a l'infini du
groupe modulaire de Teichm\"uller.} Ann. Sci. Ecole Norm. Sup.
(4) {\bf 12} (1979), 295-333.

\bibitem{C}  
C. Chevalley, {\em   Une d\'emonstration d'un th\'eor\`eme sur les
groups alg\'ebriques,} Journ. Math. Pures Appl. (9) {\bf 39} (1960),
307--317. 

\bibitem{cox}   
D. Cox, {\em  Toric varieties and toric resolutions,} this
volume. 

\bibitem{Deligne}  
P. Deligne, {\em  Le lemme de Gabber.} In: S\'em. sur les pinceaux
arithm\'etiques: la conjecture de Mordell (Ed. L.. Szpiro).
Ast\'erisque {\bf 127} (1985), 131-150.

\bibitem{D.M}  
P. Deligne and  D. Mumford, {\em  The irreducibility of the space of 
curves of given genus,} Inst. Hautes \'Etudes
Sci. Publ. Math. No. {\bf 36} (1969), 75--109. 

\bibitem{E.H}   
D. Eisenbud and J. Harris, {\em Schemes: the language of
modern algebraic geometry}. Wadsworth \& Brooks/Cole Adv. Books
Software, Pacific Grove, CA, 1992. Forthcoming edition as: {\em Why
Schemes?}, Springer. 

\bibitem{Enc-Vil}  
S. Encinas and O. Villamayor,{\em Constructive desingularization and
equivariance}, this volume. 

\bibitem{Fulton}   
W. Fulton,  {\it Introduction to toric varieties}, Annals
of 
Math. Studies 131, Princeton Univ. Press, Princeton, NJ, 1993.

\bibitem{vG.FO}  
B. van Geemen  and F. Oort, {\em A compactification of a fine moduli 
spaces of curves.} This volume. 



\bibitem{Gies}  
D. Gieseker, {\em Lectures on moduli for curves.} Published for the
Tata Institute of Fundamental Research, Bombay; Springer-Verlag,
Berlin-New York, 1982 



\bibitem{Groebner}  
W. Gr\"{o}bner, {\em Moderne algebraische Geometrie. Die
idealtheoretischen Grundlagen}, Springer, Wien und Innsbruck, 1949. 



\bibitem{FGA}  
A. Grothendieck, {\em Fondements de la g\'eom\'etrie alg\'ebrique.} Extraits du
S\'em. Bourbaki, 1957 - 1962. Secr\'etariat Math., Paris, 1962 

\bibitem{g}   
A. Grothendieck, {\it Rev\^{e}tements \'etales et groupe 
fondamental}, Lecture Notes in Math., 224, Springer, Berlin, 1971.

\bibitem{gr-mu}   
A. Grothendieck and J. P. Murre, {\it The tame fundamental group of a
formal neighbourhood of a divisor with normal crossings on a scheme},
Lecture Notes in Math., 208, Springer, Berlin, 1971.

\bibitem{H.M}  
J. Harris  and D. Mumford, {\em  On the Kodaira dimension of the moduli space
of curves.} Invent. Math. {\bf 67} (1982), 43-70. 

\bibitem{harris}   
J. Harris, {\em Algebraic Geometry - a first course}.
Grad. texts in Math. {\bf  133}, Springer - Verlag, 1992.

\bibitem{hironaka}    
H. Hironaka, {\it Resolution of singularities of an
algebraic variety over a field of characteristic zero: I, II,\/} Ann. of Math.
(2) {\bf 79} (1964), 109-326.

\bibitem{Ig}  
J.-I. Igusa, {\it Arithmetic variety of moduli for genus two.}
Ann. Math. {\bf 72} (1960), 612-649. 


\bibitem{dj2}   
A. J. de Jong, {\em Families of curves and
alterations},  Ann. Inst. Fourier (Grenoble) {\bf 47} (1997), no.~2, 599--621.

\bibitem{dJ.FO}   
A. J. de Jong and F. Oort, {\em On extending families of
curves}, Journ. Algebr. Geom. {\bf 6} (19997),  545-562. 

\bibitem{Karu1}  
K. Karu, {\em 
 Semistable reduction in characteristic 0 for families of surfaces and
 three-folds,} preprint. {\tt alg-geom/9711020}

\bibitem{Karu2}  
K. Karu, {\em 
 Minimal models and boundedness of stable varieties,} 
preprint. {\tt math.AG/9804049}

\bibitem{kato}    
K. Kato, {\em logarithmic structures of
Fontaine-Illusie}, in {\it Algebraic analysis, geometry, and number
theory (Baltimore, MD, 1988)}, 191--224, Johns Hopkins 
Univ. Press, Baltimore, MD.

\bibitem{kato1}   
K. Kato, {\em Toric singularities.} Amer. J. Math. {\bf 116} (1994),
no.~5, 1073--1099. 

\bibitem{Keel-Mori}   
S. Keel and S. Mori,  {\em Quotients by groupoids,} Ann. of
Math. (2) {\bf 145} (1997), no.~1, 193--213. 
 
\bibitem{kkms}   
G. Kempf, F. Knudsen, D. Mumford and B. Saint-Donat,
{\em Toroidal Embeddings I}, Lecture Notes in Math., 339, Springer,
Berlin, 1973. 

\bibitem{Knudsen}  
F. F. Knudsen, {\em  The projectivity of the moduli space of stable curves, 
II: the stacks $M_{g,n}$. III: The line bundles on $M_{g,n}$ and a 
proof of the projectivity of $\overline{M_{g,n}}$ in characteristic 
zero.} Math. Scand. {\bf 52} (1983). 161-199, 200-221.



\bibitem{Kollar}   
J. Koll\'ar,  {\em Projectivity of complete moduli,}
J. Differential Geom. {\bf 
32} (1990), no.~1,
235--268.

\bibitem{Lichtenbaum}
S. Lichtenbaum, {\em Curves over discrete valuation rings.} Amer. Journ. Math. {\bf 90} (1968), 380--405.


\bibitem{Loo}  
E. Looijenga, {\em  Smooth Deligne-Mumford compactifications by means
of Prym level structures.}  Jour. Algebr. Geom. {\bf 3} (1994),
283-293. 

\bibitem{Matsumura}  
H. Matsumura,  {\it Commutative ring
theory}, Translated from the Japanese by M. Reid, Second edition,  
Cambridge Univ. Press, Cambridge, 1989; MR 90i:13001 

\bibitem{Mochizuki}  
S. Mochizuki, {\em Extending Families of Curves I, II.} RIMS
preprints 1189 and 1188, 1998.

\bibitem{Mostafa}  
M. Mostafa, {\em  Die Singularit\"aten de Modulmannigfaltigkeiten
$\overline{M}_g(n)$ der stabilen Kurven vom Geslecht $g \geq 2$ mit
$n$-Teilungspunkstruktur.} Journ reine angew. Math. (Crelle),
{\bf 343} (1983), 81--98.


\bibitem{Mumford-WH}  
D. Mumford, {\em The boundary of moduli problems.} In: Lect. Notes
Summ. Inst. Algebraic Geometry, Woods Hole 1964; 8 pp., {\em unpublished.} 

\bibitem{Curves-on-Surfaces}  
D. Mumford, {\it Lectures on curves on an algebraic surface},
Princeton Univ. Press, Princeton, N.J., 1966.


\bibitem{Curves+Jacobians}   
D. Mumford, {\em Curves and their Jacobians}, The University of
Michigan Press, Ann Arbor, Mich., 1975. 

\bibitem{FO}  
F. Oort, {\em  Good and stable reduction of abelian
varieties.} Manuscr. Math. {\bf 
11} (1974), 171 - 197.


\bibitem{FO2}  
F. Oort, {\em Singularities of the moduli scheme for curves of genus
three.} Nederl. Akad. Wetensch. Proc. Ser. A {\bf
78}=Indag. Math. {\bf 37} (1975), 170--174. 

\bibitem{FO3}  
F. Oort,  {\em Singularities of coarse moduli schemes.}
in {\it S\'eminaire d'Alg\`ebre Paul Dubreil, 29\`eme ann\'ee
(Paris,1975--1976)}, 61--76, Lecture Notes in 
Math., 586, Springer, Berlin, 1977. 



\bibitem{FO4}  
F. Oort, {\em The algebraic fundamental group.} in {\it Geometric
Galois actions, 1}, 67--83, Cambridge Univ. Press, Cambridge 1997. 

 


\bibitem{Par}  
 K. Paranjape, {\em Bogomolov-Pantev resolution - An expository
 account.} To appear in: proceedings of the Warwick Algebraic
Geometry Conference July/August 1996. 

\bibitem{P.dJ}  
M. Pikaart and A. J. de Jong, {\em Moduli of curves with non-abelian
level structure.} In: The moduli space of curves (Ed. R. Dijkgraaf,
C. Faber and G. van der Geer), Proc. 1994 Conference Texel, PM {\bf
129}, Birkh\"auser, 1995, pp. 483 - 510. 

\bibitem{Popp1}  
H. Popp, {\em   Moduli theory and classification theory of algebraic
varieties.} 
Lecture Notes in Math., 620, Springer, Berlin, 1977.



\bibitem{Popp2}  
H. Popp, {\em    On the moduli of algebraic varietes III, Fine moduli spaces.}
Compos. Math.  {\bf 31} (1975), 237 - 258.

\bibitem{Popp3}  
H. Popp, {\em    The singularities of moduli schemes of curves.}
Journ. Number theory {\bf 1} (199), 90 - 107.

\bibitem{Rauch}  
H. E. Rauch, {\em The singularities of the modulus space.} Bull. Amer. Math. Soc. {\bf 68} (1962), 390 - 394.


\bibitem{Raynaud}  
M. Raynaud, {\em  Compactification du module des courbes.} in {\it
S\'eminaire Bourbaki (23\`eme ann\'ee, 1970/1971), Exp. No. 385},
47--61. Lecture Notes in 
Math., 244, Springer, Berlin, 1971. 

\bibitem{R.G}  
M. Raynaud and L. Gruson,  {\it Crit\`eres de platitude et de
projectivit\'e, Technique de ``platification" d'un module,}
Invent. Math. {\bf 13} (1971), 1-89. 


\bibitem{R}
B. Riemann,    {\it Theorie der Abel'schen Funktionen.} Journ. reine angew. Math. {\bf 54}  (1857), pp. 115 - 155.

\bibitem{Roberts}  
P. C. Roberts, {\em Intersection theory and the homological
conjectures in commutative algebra.} Proceed. ICM, Kyoto
1990. Math. Soc. Japan \& Springer - Verlag, 1991; Vol. I, pp. 361 -
368. 




\bibitem{Sam}  
P. Samuel, {\em Invariants arithm\'etiques des courbes de genre 2}
(d'apr\`es Jun Ichi Igusa), in {\it S\'eminaire Bourbaki, Vol.\ 7}, (1961/62)
Exp.\ \ 228, 81--93, Soc. Math. 
France, Paris.


\bibitem{Schlessinger}  
M. Schlessinger, {\em  Functors of Artin rings,}
Trans. Amer. Math. Soc. {\bf 130} (1968), 208--222 
 
\bibitem{S}  
J.-P. Serre, {\em Rigidit\'e du foncteur de Jacobi d'echelon $n \geq
3$.}  Appendix of Exp. 17 of S\'eminaire Cartan 1960/61.  



\bibitem{Serre}  
J.-P. Serre,  {\em Alg\`ebre locale - multiplicit\'es.}  Cours au
Coll\`ege de France, 1957--1958, r\'edig\'e par Pierre
Gabriel. Seconde \'edition, 1965. Lecture 
Notes in Math., 11, Springer, Berlin, 1965.



\bibitem{S.T}  
J.-P. Serre and J. Tate, {\em   Good reduction of abelian
varieties.} Ann. of Math. (2) {\bf 88} (1968), 492--517.  [Serre \OE II, 79.] 


\bibitem{Shimura}  
G. Shimura,    {\it On the field of rationality for an abelian
variety.} Nagoya Math. Journ. {\bf 45}, (1971), 167 -- 178. 



\bibitem{Silv}  
J. Silverman, {\em   The arithmetic of elliptic curves.} Grad. texts in
Math. {\bf 
106}, Springer - Verlag, 1986.

\bibitem{Villa}  
O. Villamayor, {\it Constructiveness of Hironaka's resolution,}
Ann. Sci. \'Ecole Norm. Sup. (4) {\bf 22} (1989), no.~1, 1--32. 

%
 \end{thebibliography}
\end{document}